\documentclass[aos,preprint]{imsart}

\usepackage{latexsym}
\usepackage{amssymb}
\usepackage{amsmath}
\usepackage{amsthm}
\usepackage{mathrsfs}
\usepackage[numbers]{natbib}
\usepackage{stmaryrd}
\usepackage{url}
\usepackage{longtable}
\usepackage{color}
\usepackage[resetlabels]{multibib}
\usepackage{dsfont}
\usepackage[T1]{fontenc}
\usepackage{lmodern}
\RequirePackage[colorlinks,citecolor=blue,urlcolor=blue]{hyperref}
\usepackage{mathtools}

\usepackage{graphicx}
\usepackage[latin1]{inputenc}
\usepackage{enumitem}

\usepackage{amsopn}
\usepackage{comment}

\newcommand{\bfi}{\begin{fig}}
\newcommand{\efi}{\end{fig}}
\newcommand{\btab}{\begin{tab}}
\newcommand{\etab}{\end{tab}}
\newcommand{\barr}{\begin{array}}
\newcommand{\earr}{\end{array}}
\newcommand{\beq}{\begin{equation}}
\newcommand{\eeq}{\end{equation}}
\newcommand{\bdis}{\begin{displaymath}}
\newcommand{\edis}{\end{displaymath}\noindent}

\newcommand{\bbn}{\mathbb{N}}

\newcommand{\bbr}{\mathbb{R}}

\newcommand{\bbe}{\mathbb{E}}

\newcommand{\bbp}{\mathbb{P}}

\newcommand{\bone}{\mathds 1}

\newcommand{\R}{\mathbb{R}}
\newcommand{\N}{\mathbb{N}}

\newcommand{\limw}{\stackrel{\rm w}{\longrightarrow}}
\newcommand{\limst}{\stackrel{\mathrm{st}}{\Longrightarrow}}
\newcommand{\lims}{\stackrel{\mathrm{st}}{\longrightarrow}}
\newcommand{\limL}{\stackrel{L^1}{\Longrightarrow}}
\newcommand{\limp}{\stackrel{\mathbb{P}}{\longrightarrow}}

\newcommand{\ucp}{\stackrel{\rm u.c.p.}{\Longrightarrow}}

\newcommand{\lec}{\lesssim}

\newcommand{\cals}{{\cal S}}
\newcommand{\calc}{{\cal C}}

\newcommand{\calf}{{\cal F}}

\newcommand{\calb}{{\cal B}}
\newcommand{\call}{{\cal L}}

\newcommand{\calm}{{\cal M}}

\newcommand{\calv}{{\cal V}}

\newcommand{\calz}{{\cal Z}}

\newcommand{\al}{{\alpha}}
\newcommand{\la}{{\lambda}}
\newcommand{\La}{{\Lambda}}
\newcommand{\eps}{{\epsilon}}

\newcommand{\ga}{{\gamma}}
\newcommand{\Ga}{{\Gamma}}

\newcommand{\si}{{\sigma}}
\newcommand{\om}{{\omega}}
\newcommand{\Om}{{\Omega}}

\newcommand{\var}{{\mathrm{Var}}}
\newcommand{\cov}{{\mathrm{Cov}}}

\newcommand{\ov}{\overline}
\newcommand{\un}{\underline}
\newcommand{\wh}{\widehat}
\newcommand{\wt}{\widetilde}

\newcommand{\dd}{\mathrm{d}}
\newcommand{\ee}{\mathrm{e}}
\newcommand{\ii}{\mathrm{i}}

\newcommand{\bb}{\mathrm{b}}

\newcommand{\frc}{\mathfrak{c}}
\newcommand{\frr}{\mathfrak{r}}

\newcommand{\lb}{\left[}
\newcommand{\rb}{\right]}
\newcommand{\lp}{\left(}
\newcommand{\rp}{\right)}
\newcommand{\lv}{\left|}
\newcommand{\rv}{\right|}
\newcommand{\bv}{\bigg|}
\newcommand{\as}{\qquad\text{as}\qquad n\to\infty}
\newcommand{\sumt}{\sum_{i=1}^{t^*_n}}
\newcommand{\sumT}{\sum_{i=1}^{T^*_n}}
\newcommand{\del}{\Delta_n}
\newcommand{\Del}{\Delta_n}
\newcommand{\Delh}{\Delta_n^{\frac12}}
\newcommand{\delh}{\Delta_n^{\frac12}}

\newcommand{\supt}{\sup_{t\in[0,T]}}
\newcommand{\pk}{\pi^2\kappa |\xi|^2}

\theoremstyle{plain}
\newtheorem{Theorem}{Theorem}[section]
\newtheorem{Corollary}[Theorem]{Corollary}
\newtheorem{Lemma}[Theorem]{Lemma}
\newtheorem{Proposition}[Theorem]{Proposition}
\theoremstyle{definition}
\newtheorem{Definition}[Theorem]{Definition}
\newtheorem{Example}[Theorem]{Example}
\newtheorem{Assumption}{Assumption}

\theoremstyle{remark}
\newtheorem{Remark}[Theorem]{Remark}

\newcommand{\bthm}{\begin{Theorem}}
\newcommand{\ethm}{\end{Theorem}}

\newcommand{\bcor}{\begin{Corollary}}
\newcommand{\ecor}{\end{Corollary}}

\newcommand{\blem}{\begin{Lemma}}
\newcommand{\elem}{\end{Lemma}}

\newcommand{\bprop}{\begin{Proposition}}
\newcommand{\eprop}{\end{Proposition}}

\newcommand{\bdf}{\begin{Definition}}
\newcommand{\edf}{\end{Definition}}

\newcommand{\bex}{\begin{Example}}
\newcommand{\eex}{\end{Example}}

\newcommand{\brem}{\begin{Remark}}
\newcommand{\erem}{\end{Remark}}

\newcommand{\bass}{\begin{Assumption}}
\newcommand{\eass}{\end{Assumption}}

\newcommand{\bpr}{\begin{proof}}
\newcommand{\epr}{\end{proof}}

\newcommand{\benu}{\begin{enumerate}}
\newcommand{\eenu}{\end{enumerate}}

\newcommand{\bit}{\begin{itemize}}
\newcommand{\eit}{\end{itemize}}

\newcommand{\itt}{\textit}

\numberwithin{equation}{section}

\newcites{A,B}{References,References}

\allowdisplaybreaks[2]

\begin{document}

\begin{frontmatter}
	
	\title{High-frequency analysis of parabolic stochastic PDEs}
	\runtitle{High-frequency analysis of parabolic stochastic PDEs}
	
	\begin{aug}
		\author{\fnms{Carsten} \snm{Chong}\ead[label=e1]{carsten.chong@epfl.ch}}
		
		\runauthor{C. Chong}
		
		\affiliation{École Polytechnique Fédérale de Lausanne}
		
		\address{Institut de mathématiques\\ École Polytechnique Fédérale de Lausanne\\ Station 8\\ CH-1015 Lausanne\\ \printead{e1}}
	\end{aug}

	\begin{abstract}
		We consider the problem of estimating stochastic volatility for a class of second-order parabolic stochastic PDEs. Assuming that the solution is observed at high temporal frequency, we use limit theorems for multipower variations and related functionals to construct consistent nonparametric estimators and asymptotic confidence bounds for the integrated volatility process. As a byproduct of our analysis, we also obtain feasible estimators for the regularity of the spatial covariance function of the noise. 
	\end{abstract}
	
	
	
\end{frontmatter}

\section{Introduction}\label{intro}

A central objective of stochastic modeling is to capture the fluctuations of a system evolving under the influence of random noise. Being able to quantify the degree of variability and uncertainty in such a system is inevitable for the control and prediction of its future behavior. 

In the financial and econometrics literature, a key concept designed to measure and describe the amount of randomness present in the evolution of asset prices, interest rates, or other financial indices is that of \emph{stochastic volatility}. Over the past decades, a huge amount of work has been devoted to building stochastic volatility models that are able to reproduce stylized features found in empirical financial data. We only refer to \citeA{Brockwell09} for a comprehensive overview.

Of course, the notion of stochastic volatility is not only limited  to mathematical finance. For example, in the literature of turbulence, it is commonly referred to as \emph{intermittency}; see \citeA{BN04,Pereira16,Robert08} for various models of stochastic intermittency. In a related context, the phenomenon of intermittency has also been intensively studied in the theory of \emph{stochastic partial differential equations} (stochastic PDEs). To be more precise, let us consider a parabolic stochastic PDE of the form
\beq\label{SHE} \partial_t Y(t,x)=\frac{\kappa}{2} \Delta Y(t,x) - \la Y(t,x) +\si(t,x)\dot W(t,x), \eeq
where $\kappa>0$ is a diffusion or viscosity constant, $\la>0$ is a damping rate, $\si$ is a predictable random field, and $\dot W$ is a Gaussian noise. When we consider \eqref{SHE} for $t\geq0$ and $x\in\R$ with an initial condition at $t=0$ that is bounded away from $0$, it is known from \citeA{Foondun09} that if $\si$ is a linear function of the solution with sufficient growth, then the solution exhibits a strong mass concentration at large times by forming exponentially large peaks on exponentially small areas. On the other hand, when $\si$ is a bounded function of the solution, this kind of intermittent behavior does not occur. Hence, the knowledge of the form of $\si$ is essential for determining the behavior of the solution $Y$ to \eqref{SHE}. 

Furthermore, in many applications of \eqref{SHE}, or stochastic PDEs of a similar form, the random field $\sigma$ models the level of noise that acts on a process described by an otherwise deterministic PDE. Examples include \citeA{Cont05} on term structure models,  \citeA{Denaro13,ElSaadi07} on plankton distribution,  \citeA{Jones99} on the motion of particles in gravitational fields, \citeA{Sigrist12,Sigrist15} on precipitation models, and \citeA{Tuckwell13,Walsh81} on neuron spikes. In these applications, the knowledge of $\sigma$ is essential for assessing to which degree the solution to \eqref{SHE} deviates from the solution to the deterministic PDE.

\subsection{Objective and related literature}

Motivated from the applications mentioned above,
the purpose of the present article is to establish consistent estimators and asymptotic confidence bounds for the random field $\si$ in \eqref{SHE}, which we henceforth call the \emph{(stochastic) volatility process} (even outside the financial context). To this end, we assume that we are given observations of a single path of the solution $Y(t,x)$ at a finite number of spatial points $x_1,\ldots,x_N \in\R^d$ and at a high number of time points $t=\del, 2\del,\ldots, [T/\del]\del$ within a finite interval $[0,T]$ with $T<\infty$ ($[\,\cdot\,]$ stands for the integer part). Here, $\del$ is a small time step, and we seek estimators of $\si$ with the properties mentioned above when $\del\to0$. Hence, our observation scheme has \emph{high frequency} in time and \emph{low resolution} in space. This is a realistic framework for many of the applications mentioned above, where high-frequency recordings are only available at a small number of measuring sites.

The high-frequency analysis of Itô semimartingales has been fully accomplished in the past ten years; see the treatises \citeA{AitSahalia14,Jacod12} for a complete account. For example, given a continuous semimartingale $X(t)=\int_0^t \si(s)\,\dd B(s)$ where $B$ is a Brownian motion and $\si$ a predictable process, the basic idea to estimate $\si$ is to consider \emph{(normalized) power variations} of $X$, i.e.,
\beq\label{pvar} V^n_p(X,t) = \del \sum_{i=1}^{[t/\del]} \lv\frac{\Delta^n_i X}{\sqrt{\del}}\rv^p,\qquad t\in[0,T], \qquad p>0, \eeq
where $\Delta^n_i X = X(i\del)-X((i-1)\del)$ is an \emph{increment} of $X$ from $(i-1)\del$ to $i\del$. Under minimal assumptions on $\si$, one can show that
\beq\label{LLN-ito} V^n_p(X,t) \ucp V_p(X,t)=\bbe[|Z|^p]\int_0^t |\si(s)|^p\,\dd s, \eeq
where $Z\sim N(0,1)$ and $\ucp$ denotes uniform convergence in probability on compacts; see Theorem~3.4.1 in \citeA{Jacod12}. Under further regularity assumptions on $\si$, we have an associated central limit theorem of the form
\beq\label{CLT-ito} \del^{-\frac12}(V^n_p(X,t)-V_p(X,t)) \limst \calz, \eeq
where $\calz$ is a Gaussian process with independent increments and explicitly known variance, conditionally on $\si$; see Theorem~5.3.6 in \citeA{Jacod12}. In \eqref{CLT-ito}, $\limst$ denotes functional stable convergence in law with respect to the uniform topology; see Section~3.2 in \citeA{AitSahalia14}, Section~2.2 in \citeA{Jacod12}, or \citeA{Podolskij10}. The two results \eqref{LLN-ito} and \eqref{CLT-ito} can then be used to construct asymptotic confidence bounds for the \emph{integrated volatility process} $\int_0^t |\si(s)|^p\,\dd s$, which is what we understand by ``estimating $\si$'' throughout this article. 

When we leave the class of semimartingales and consider a moving average process of the form $X_t=\int_{-\infty}^t g(t-s)\si(s)\,\dd B_s$, where $g$ is a kernel that is smooth except at the origin, then the functionals 
\beq\label{pvar-ma} V^n_p(X,t) = \del \sum_{i=1}^{[t/\del]} \lv\frac{\Delta^n_i X}{\tau_n}\rv^p,\qquad t\in[0,T], \qquad p>0, \eeq
with a normalizing factor $\tau_n$ depending on the singularity of $g$ at the origin, still satisfy \eqref{LLN-ito} and \eqref{CLT-ito} (with  a slightly larger  variance for $\calz$) if $g$ and $\si$ are sufficiently regular; see \citeA{BN11,BN13,Corcuera13}. In particular, this applies to fractional Brownian motion with Hurst parameter $H<\frac 34$; see  \citeA{Corcuera06,Istas97}.

 In the context of stochastic PDEs, estimation problems have been considered by many authors; see \citeA{Cialenco18} for a recent survey.  The majority of literature in this respect focuses on the estimation of $\kappa$, assuming that $\si$ is constant and known, and that the solution to \eqref{SHE}, or certain transformations thereof, is observed continuously in time and/or space; see \citeA{Cialenco18,Lototsky17} and the references therein. In practice, of course, measurements are \emph{discrete}, and the amount of literature is much smaller when it comes to estimating $\si$ based on discrete observation schemes. 
 
 When $\si$ is a deterministic function of $t$ only, and \eqref{SHE} is considered in one spatial dimension on an interval, \citeA{PR02} constructs an estimator for $\si$ based on high-frequency observations in time of a fixed number of Fourier coefficients of $Y$. But these are solutions to certain stochastic differential equations and hence semimartingales, so the estimation problem can be fully solved by the techniques of \citeA{AitSahalia14,Jacod12}. By contrast, the solution $Y$ itself at fixed spatial positions is \emph{not} a semimartingale (if $\dot W$ is a Gaussian space-time white noise and $d=1$, it has a nontrivial finite quartic variation in time; see \citeA{Swanson07,Walsh81}). Assuming high-frequency observations of $Y$ in time (as we do in this work), but still with deterministic $\si$ that only depends on $t$ (but not on $x$), the papers \citeA{Bibinger17,Cialenco17} use power variations as in \eqref{pvar-ma} with $p=4$ and $p=2$, respectively, to establish asymptotic confidence bounds for $\si$. We also refer to \citeA{Liu16} for related results from a probabilistic point of view and to \citeA{Bibinger19} for some extensions of \citeA{Bibinger17}.
 
With stochastic $\si$, \citeA{Pospisil07} shows a variant of \eqref{LLN-ito} with $p=4$ for the solution $Y$ to \eqref{SHE} if $d=1$ and $\dot W$ is a space-time white noise; see also \citeA{Cialenco17,Foondun15}. The papers \citeA{BN11b,Pakkanen14} contain certain limit theorems for space-time moving averages when $\si$ is independent of the noise $W$ (so by conditioning on $\si$, this reduces to the case of deterministic $\si$). Apart from these particular cases, to our best knowledge, no further results are available for \eqref{SHE}, and in particular, no central limit theorems as in \eqref{CLT-ito} exist in the case of stochastic $\si$, and not even a law of large numbers as in \eqref{LLN-ito} if $\dot W$ is a spatially colored noise. 

It is therefore the main objective of this paper to fill this gap and to derive consistent estimators and confidence bounds for the integrated volatility process (with respect to time and for fixed values of $x$) if $\si$ is a random field. In fact, we consider much more general functionals than \eqref{pvar-ma}. Given a sufficiently regular evaluation function $f\colon \R^{N\times L}\to\R^M$ with $L, M \in\N$, we consider \emph{(normalized) variation functionals} of the form $V^n_f(Y_x,t) = (V^n_f(Y_x,t)_1,\ldots,V^n_f(Y_x,t)_M)'$ where
\beq\label{Vnf} V^n_f(Y_x,t)_m = \Delta_n \sum_{i=1}^{[t/\Delta_n]-L+1} f_m\left( \frac{\Delta^n_i Y_{x}}{\tau_n},\ldots,\frac{\Delta^n_{i+L-1} Y_{x}}{\tau_n} \right) \eeq
for $t\in[0,T]$ and $m=1,\ldots,M$.
Here, $Y_x$ is the $N$-dimensional process whose $j$th component is $Y(\cdot,x_j)$, $\Delta^n_i Y_x=Y_x(i\Delta_n)-Y_x((i-1)\Delta_n)=(\Delta^n_i Y(\cdot,x_1),\ldots,\Delta^n_i Y(\cdot, x_N))'$, and $\tau_n$ is a normalizing factor to be introduced in Section~\ref{model}. The main examples for $f$ are multipowers, which we will study in detail in Section~\ref{multipower}. The reader may consult \citeA{AitSahalia14,Jacod12} and \citeA{BN11,BN13,Corcuera13} for multipower variations of semimartingales and moving averages, respectively.

\subsection{Results and methodology}

After a short introduction to stochastic PDEs, the two main limit theorems are formulated in Section~\ref{limit}. Theorem~\ref{LLN} gives a law of large numbers for the functionals in \eqref{Vnf}, while Theorem~\ref{CLT} gives the associated central limit theorem at a rate of $\sqrt{\del}$. Section~\ref{multipower} shows how these limit theorems apply to the important example of multipower variations (see Corollary~\ref{multipower-result}), which will then be used in Section~\ref{est} to construct feasible estimators for $\si$. It turns out that we can even estimate the spatial correlation structure of the noise $\dot W$ in \eqref{SHE}, which we assume to be parametrized by an exponent $\al$. Theorem~\ref{known-alpha} addresses the case of estimating $\si$ when $\al$ is known, while Theorems~\ref{est-alpha} and \ref{correst} propose two estimation procedures for $\al$ and Theorem~\ref{unknown-alpha} one for $\si$ when $\al$ is unknown. Our results indicate that this spatial correlation index $\al$ plays a very similar role to the kernel smoothness parameter in \citeA{BN11,BN13,Corcuera13}.

The proofs will be given in Section~\ref{proofs}. As we already know from the semimartingale case considered in \citeA{AitSahalia14,Jacod12}, having a stochastic instead of a deterministic volatility process complicates the proofs considerably, in particular for the central limit theorem. For instance, the proofs of \citeA{Bibinger17,Cialenco17} do not apply in our context as they make heavy use of the Gaussian distribution of the solution $Y$ to \eqref{SHE} when $\si$ is nonrandom. Conceptually, $Y$ can be viewed as a moving average process in space and time; see formula \eqref{Yx} below, so it seems natural to transfer techniques from \citeA{BN11,BN13,Corcuera13} to the stochastic PDE setting.  A crucial step in their proofs is to factorize the volatility process out of the stochastic integral by discretizing $\si$ along a subgrid of $\del, 2\del,\ldots,[T/\del]\del$ \emph{before} showing the actual central limit theorem, and then to prove, using fractional calculus methods from \citeA{Corcuera14}, that this discretization only induces an asymptotically negligible error. If one wishes to apply this method to stochastic PDEs, one would have to discretize the volatility process $\si$ both in time and space. Although the heat kernel in \eqref{Yx} is concentrated around the origin, in general, this localization is simply not strong enough \emph{in space} on a $\sqrt{\del}$ rate, which would be needed for the central limit theorem. Thus, we see no way to apply the methods of \citeA{BN11,BN13,Corcuera13} to \eqref{SHE}; cf.\ part (4) of Remark~\ref{rem-discr}.

Instead, we will show that a \emph{combination} of the \emph{martingale methods} of \citeA{AitSahalia14,Jacod12} (for the discretization part and the identification of the limit law) with \emph{analysis on the Wiener space} as in \citeA{BN11,BN13,Corcuera13} (for tightness) will give the desired central limit theorem for \eqref{Vnf}. The advantage of this strategy is that we only need to make spatial approximations of $\si$ \emph{after} the actual central limit theorem, where we can use \emph{symmetry} properties of certain measures related to the heat kernel to compensate for its bad spatial concentration properties. Since $Y$ is not a semimartingale, for this method to work, we have to use a complex procedure to approximate $V^n_f(Y_x,t)$ by martingale-type sums in a first step. 

Our  ``martingale proof'' also provides an interesting alternative to proving limit theorems for moving average processes in the purely temporal case. For instance, with the new method, the results of \citeA{BN11,BN13,Corcuera13} can be extended to allow for volatility processes that are semimartingales, which include the majority of stochastic volatility models available in the literature. Moreover, we believe that the martingale techniques we develop in this paper will pave the way for further statistical procedures to estimate spot volatility, to handle measurement errors, or to detect and estimate the density of jumps for stochastic PDEs (and moving average processes). We refer to Chapter~III.8 in \citeA{AitSahalia14}, \citeA{Jacod10}, and \citeA{AitSahalia09, AitSahalia09b,Jacod14}, respectively, for the corresponding results in the semimartingale framework, which are all proved with martingale techniques. 
 
This paper is accompanied by some supplementary material in \citeA{Chong18a}. All references and numberings starting with a letter, like (A.1) or Lemma~B.1, refer to \citeA{Chong18a}, except for Assumptions~\ref{AssLLN} and \ref{AssCLT}, which are stated in Section~\ref{model}.

In what follows, we often write $\iint_a^b =\int_a^b \int_{\R^d}$, $\iint=\iint_{-\infty}^\infty$, $\iiint_a^b = \int_a^b\int_{\R^d}\int_{\R^d}$, and $\iiint = \iiint_{-\infty}^\infty$. Moreover, $\bbn=\{1,2,\ldots\}$ and $\bbn_0=\{0,1,2,\ldots\}$.

\section{Model and main results}\label{model}

On a given filtered probability space $(\Om,\calf,(\calf_t)_{t\in\R},\bbp)$ satisfying the usual conditions, we consider the stochastic PDE \eqref{SHE} for $t\in\R$ and $x\in\R^d$  driven by a zero-mean $(\calf_t)_{t\in\R}$-Gaussian noise $\dot W$ which is white in time but possibly colored in space. More precisely, we have an $L^2$-valued centered Gaussian measure $W(A)$, indexed by bounded Borel sets $A\in\calb_\bb(\R^{d+1})$, such that $W(A)$ is independent of $\calf_t$ if $A\cap((-\infty,t]\times\R^d)=\emptyset$, and such that for any $A_1,A_2\in\calb_\bb(\R^{d+1})$,
\[ \bbe\lb W(A_1)W(A_2) \rb = \iiint\bone_{A_1}(s,y)\bone_{A_2}(s,z)\,\La(\dd y,\dd z)\,\dd s.   \]
In this paper, we assume that $\La(\dd y,\dd z) = F(z-y)\,\dd y\,\dd z$ where $F$ is the \emph{Riesz kernel}  $F(x)=c_\al|x|^{-\al}$ for some $\al\in(0,d\wedge 2)$ and the normalizing constant is given by $c_\al = \pi^{d/ 2-\al} \Ga(\frac \al 2)/\Ga(\frac{d-\al}{2})$. Here and throughout the paper, $|x|$ denotes the Euclidean norm. In dimension $1$, we also allow for the case where $\dot W$ is a \emph{Gaussian space-time white noise}, which corresponds to $\La(\dd y,\dd z) = \delta_z(\dd y)\,\dd z$ and $F(x)=\delta_0(x)$. We set $\al=1$ and $c_\al=1$ in this case. In dimensions $d\geq2$, it is well known that no function-valued solution to \eqref{SHE} exists if $\dot W$ is a space-time white noise, or if $2\leq \al < d$; see \citeA{Dalang99}.

By the classical integration theory of \citeA{Walsh86} (see also \citeA{Dalang99,Nualart07} for extensions), an Itô integral against $W$ can be constructed for integrands from the space $\call(W)$ of predictable random fields $\phi\colon \Om\times\R\times\R^d\to\R$ satisfying
\[ \bbe\lb \iiint |\phi(s,y)\phi(s,z)|\,\La(\dd y,\dd z)\,\dd s\rb<\infty. \]
In particular, as soon as
the predictable random field $\si$ satisfies
\beq\label{bounded-moments} \sup_{(t,x)\in\R\times\R^d} \bbe[\si^2(t,x)]<\infty, \eeq
the stochastic PDE \eqref{SHE} for $(t,x)\in \R\times \R^d$ admits a \emph{mild solution} given by
\beq\label{Yx} Y(t,x) = \iint_{-\infty}^t G(t-s,x-y)\si(s,y) \,W(\dd s,\dd y),\qquad (t,x)\in\R\times\R^d, \eeq
where
\beq\label{heatkernel} G(t,x)=G_x(t)=(2\pi\kappa t)^{-\frac d2}\ee^{-\frac{|x|^2}{2\kappa t}-\la t} \bone_{t>0},\qquad (t,x)\in\R\times\R^d, \eeq
is the heat kernel for \eqref{SHE}. We remark that although the integration theory in \citeA{Dalang99,Nualart07,Walsh86} is developed for $t\geq0$, their results extend without any change to the  case  $t\in\R$. Also, while all our results are formulated for \eqref{SHE} with $t\in\R$, they remain valid if \eqref{SHE} is considered for $t\geq0$ as soon as the initial condition at $t=0$ is sufficiently regular; see Remark~\ref{non-stat} below.

As soon as $\si$ is jointly stationary with the increments of $W$, the mild solution in \eqref{Yx} is stationary in space and time. In particular, if $\si\equiv 1$, all components of  $\Delta^n_i Y_x$ are normally distributed with mean $0$ and variance
\beq\label{taun}\begin{split} \qquad\tau^2_n&=\bbe\left[\lv \iint_{-\infty}^{i\del} (G_{x_j-y}(i\del-s)-G_{x_j-y}((i-1)\del-s))\,W(\dd s,\dd y)\rv^2\right]\\
		& =\iiint_0^\infty  (G_y(s)-G_y(s-\del))(G_z(s)-G_z(s-\del))\,\La(\dd y,\dd z)\,\dd s,
\end{split}\raisetag{-2.75\baselineskip}\eeq 
which depends neither on $i$ nor on $j$. This will be the normalizing factor we choose in \eqref{Vnf} so that $V^n_f(Y_x,t)$ is typically ``of order $1$'' and we may hope for convergence as $n\to\infty$. An explicit formula for $\tau_n$ can be found in Lemma~\ref{Ga} in the supplementary material \citeA{Chong18a}.

\subsection{Limit theorems for normalized variation functionals}\label{limit}

A first-order limit theorem for the normalized variation functionals \eqref{Vnf} can be shown under mild assumptions on $f$ and $\si$. In what follows, the Euclidean norm $|z|$ for some matrix $z\in\R^{N\times L}$ is defined by viewing $z$ as an element of $\R^{NL}$. 

\bass\label{AssLLN} There exists $p\geq2$ with the following properties:
\benu
\item[A1.] The function $f\colon \R^{N\times L}\to\R^M$ is continuous with $f(z)=o(|z|^p)$  as $|z|\to\infty$. 
\item[A2.] For some $\eps>0$, we have
\[ \sup_{(t,x)\in\R\times\R^d} \bbe[|\si(t,x)|^{p+\eps}]<\infty. \]
\item[A3.] The random field $\si$ is uniformly $L^2$-continuous on $\bbr\times\bbr^d$: as $\eps\to0$,
\[ w(\eps)= \sup\left\{\bbe[|\si(t,x)-\si(s,y)|^2]\colon |t-s|+|x-y|<\eps\right\} \to 0. \]
\eenu
\eass

The next theorem is our first main result. If $(X^n(t))_{t\geq0}$ and $(X(t))_{t\geq0}$ are stochastic processes, we write $X^n\limL X$ or $X^n(t)\limL X(t)$ if for every $T>0$, we have $\bbe[\supt |X^n(t)-X(t)|]\to0$ as $n\to\infty$.
\bthm[Law of large numbers] \label{LLN} Under Assumption~\ref{AssLLN}, we have
\beq\label{LLN-conv}V^n_f(Y_x,t) \limL V_f(Y_x,t)= \int_0^t \mu_f\left(\si^2(s,x_1),\ldots,\si^2(s,x_N) \right)\,\dd s.\eeq
Here, $\mu_f\colon \R^N\to \R^M$ is the function given by 
$\mu_f(v_1,\ldots,v_N)=\bbe[f(Z)]$, where $Z=(Z_{jk})_{j,k=1}^{N,L}$ is multivariate normal with mean $0$ and 
\beq\label{cov-Z}\cov(Z_{j_1k_1},Z_{j_2 k_2}) = \Ga_{|k_1-k_2|}v_j \bone_{j_1=j_2=j},\eeq
and
\beq\label{Ga-formula} \Ga_0=1,\qquad \Ga_r=\frac12 \Big( (r+1)^{1-\frac\al 2}-2r^{1-\frac\al 2}+(r-1)^{1-\frac\al 2}\Big),\qquad r\geq1.  \eeq 
\ethm

\brem\label{RemLLN} Fix some $m=1,\ldots,M$.
If $f_m$ only depends on the variables $(z_{jk}\colon j\in J,~ k=1,\ldots,L)$, where $J\subseteq \{1,\ldots,N\}$, that is, if $V^n_f(Y_x,t)_m$ only uses the increments observed at the points $(x_j\colon j\in J)$, then its limit in \eqref{LLN-conv} will only depend on $(\si(\cdot,x_j)\colon j\in J)$. In other words, by taking measurements at $x_j$, one can obtain isolated information about  $\si(\cdot,x_j)$, independently from the values of $\si$ at all other positions.
\erem

In order to obtain a central limit theorem for \eqref{LLN-conv}, we need to put stronger regularity assumptions on $f$ and $\si$, which is already necessary for semimartingales (cf.\ \citeA{AitSahalia14,Jacod12}) and moving averages (cf.\ \citeA{BN11,BN13,Corcuera13}). In the first two references, $\si$ itself has to be a semimartingale, while in the next three references, $\si$ has to be (essentially) Hölder continuous with exponent $>\frac12$. We will assume that $\si$ has one of these two properties plus additional regularity in space.

\bass\label{AssCLT}~
\benu
\item[B1.] The function $f\colon \R^{N\times L}\to\R^M$ is even [i.e., we have $f(z)=f(-z)$ for all $z\in\R^{N\times L}$] and four times continuously differentiable. Moreover, there are $p\geq2$ and $C>0$ such that 
\begin{align*} &|f_m(z)|\leq C(1+|z|^p), \qquad \left|\textstyle\frac{\partial}{\partial z_\al} f_m(z)\right|\leq C(1+|z|^{p-1}),\\
&\left|\textstyle\frac{\partial^2}{\partial z_\al\,\partial z_\beta }f_m(z)\right|+\left|\textstyle\frac{\partial^3}{\partial z_\al\,\partial z_\beta\,\partial z_\ga }f_m(z)\right|+\left|\textstyle\frac{\partial^4}{\partial z_\al\,\partial z_\beta\,\partial z_\ga \,\partial z_\delta}f_m(z)\right|\\
&\quad\leq C(1+|z|^{p-2})
\end{align*} for all $m\in\{1,\ldots,M\}$ and $\al,\beta,\ga,\delta\in \{1,\ldots, N\}\times\{1,\dots, L\}$. 
\item[B2.] If $F$ is the Riesz kernel with $0<\al<1$, or $\al=1$ and $\dot W$ is \emph{not} a space-time white noise, we assume for each $m=1,\ldots, M$ that $f_m(z)$ only depends on $z_{j1},\ldots,z_{jL}$ for some $j =j(m)\in\{1,\ldots,N\}$.
\item[B3.] The volatility process $\si$ takes the form 
\beq\label{si-semi} \si(t,x) = \si^{(0)}(t,x)+\iint_{-\infty}^t K(t-s,x-y)\rho(s,y)\,W'(\dd s,\dd y)\eeq
for $(t,x)\in\R\times\R^d$,
with the following specifications:
\bit
\item $\si^{(0)}$ is a predictable process satisfying 
\beq\label{mom-si} \sup_{(t,x)\in\R\times\R^d} \bbe[|\si^{(0)}(t,x)|^{2p+\eps}]<\infty\eeq
for some $\eps>0$, and 
\beq\label{Holder} \sup_{x\in\R^d} \bbe[|\si^{(0)}(t,x)-\si^{(0)}(s,x)|^{2p}]^{\frac 1 {2p}} \leq C' |t-s|^\ga \eeq
for some $\ga\in(\frac 12,1]$ and $C'>0$.
In addition, for every $t\in\R$, the mapping $x\mapsto \si^{(0)}(t,x)$ is almost surely twice differentiable such that for $ j,k=1,\ldots,d$, we have
\beq\label{si-der} \sup_{(t,x)\in\R\times\R^d} \bbe\lb\lv \textstyle\frac{\partial}{\partial x_j} \si^{(0)}(t,x)\rv^p+\lv \textstyle\frac{\partial^2}{\partial x_j\,\partial x_k} \si^{(0)}(t,x)\rv^p\rb<\infty. \eeq 
\item $W'$ is an $(\calf_t)_{t\in\R}$-Gaussian noise that is white in time and possibly colored in space [such that $(W,W')$ is bivariate Gaussian] with $\La'(\dd y,\dd z)=F'(z-y)\,\dd y\,\dd z$, where $F'$ is the Riesz kernel with some $\al'\in(0,2)\cap(0,d]$.
\item $K\colon [0,\infty)\times\R^d\to\R$ is a kernel such that the partial derivatives $\frac{\partial}{\partial t}K$, $\frac{\partial}{\partial x_j}K$, $\frac{\partial^2}{\partial x_j\,\partial x_k}K$, and $\frac{\partial^3}{\partial x_j\,\partial x_k\,\partial x_l} K$ exist and belong to $\call(W')$ for all $j,k,l=1,\ldots,d$.
\item $\rho$ is a predictable process satisfying the same moment condition \eqref{mom-si} as $\si^{(0)}$, and furthermore, for some $\eps'>0$ and $C''>0$,
\beq\label{Holder-2} \sup_{x\in\R^d} \bbe[|\rho(t,x)-\rho(s,x)|^{2p+\eps}]^{\frac1 {2p+\eps}}\leq C'' |t-s|^{\eps'}. \eeq
\eit
\eenu
\eass

For our second main result, we use $\limst$ to denote functional stable convergence in law in the space of càdlàg functions $[0,\infty)\to\R^M$, equipped with the local uniform topology, while stable convergence in law between finite-dimensional random variables will be denoted by $\lims$. We refer the reader to \citeA{AitSahalia14,Jacod12} for a definition of this mode of convergence and also for the definition of a \emph{very good filtered extension} of $(\Om,\calf,(\calf_t)_{t\geq0},\bbp)$. The only property of stable convergence in law we need is the following (see, for example, Proposition~2(i) in \citeA{Podolskij10}):
\beq\label{stable} X_n\lims X,\qquad Y_n\limp Y \quad\implies \quad(X_n,Y_n)\lims (X,Y).  \eeq
Since the limiting objects in \eqref{LLN-conv} are random, this will allow us to studentize \eqref{CLT-statement} below and obtain feasible confidence bounds for $\si$. Just convergence in law, of course, will not suffice for this purpose.
\bthm\label{CLT} Under Assumption~\ref{AssCLT}, we have as $n\to\infty$,
\beq\label{CLT-statement} \del^{-\frac12} \lp V^n_f(Y_x,t)-V_f(Y,t)\rp \limst \calz, \eeq
where $(\calz(t)= (\calz_1(t),\ldots,\calz_M(t))')_{t\geq0}$ is a continuous process defined on a very good filtered extension $(\ov\Om,\ov\calf, (\ov\calf_t)_{t\geq0},\ov\bbp)$ of the original probability space $(\Om,\calf,(\calf_t)_{t\geq0},\bbp)$, which, conditionally on the $\si$-field $\calf$, is a centered Gaussian process with independent increments such that the covariance function $\calc_{m_1m_2}(t) = \ov\bbe[\calz_{m_1}(t) \calz_{m_2}(t)\,|\,\calf]$, for $m_1,m_2=1,\ldots,M$, is given by 
\beq\label{calc} \begin{split}\calc_{m_1m_2}(t) &= \int_0^t \rho_{f_{m_1},f_{m_2}}(0;\si^2(s,x_1),\ldots,\si^2(s,x_N))\,\dd s\\
	&\quad+\sum_{r=1}^\infty \int_0^t \rho_{f_{m_1},f_{m_2}}(r;\si^2(s,x_1),\ldots,\si^2(s,x_N))\,\dd s\\
&\quad+\sum_{r=1}^\infty \int_0^t \rho_{f_{m_2},f_{m_1}}(r;\si^2(s,x_1),\ldots,\si^2(s,x_N))\,\dd s. \end{split}\eeq 
In the last line, for $r\in\N_0$, we define 
\beq\rho_{f_{m_1},f_{m_2}}(r;v_1,\ldots,v_N)=\cov(f_{m_1}(Z^{(1)}),f_{m_2}(Z^{(2)})),\eeq 
where $Z^{(1)}=(Z^{(1)}_{jk})_{j,k=1}^{N,L}$ and $Z^{(2)}=(Z^{(2)}_{jk})_{j,k=1}^{N,L}$ are jointly Gaussian, both with the same law as the matrix $Z$ in Theorem~\ref{LLN} 
and cross-covariances  
\beq\label{cross-cov} \cov(Z^{(1)}_{j_1k_1}, Z^{(2)}_{j_2k_2})=\Ga_{|k_1-k_2+r|}v_j\bone_{j_1=j_2=j}. \eeq
Part of the statement is that the series in \eqref{calc} converge in the $L^1$-sense.
\ethm

\brem\label{rem-discr} Let us comment on the assumptions of Theorem~\ref{CLT}.
\benu
	\item Assumption~\ref{AssCLT}1 can be relaxed by allowing, for example, $f$ to be continuous but not differentiable at $z=0$, very similar to \citeA{BN11} or Chapter~11.2 in \citeA{Jacod12}. Due to the technical proofs already needed under the stronger assumptions, we refrain from doing so in this paper.
	\item Increments at different measurements sites $x_j \neq x_{j'}$ contribute in the limit $n\to\infty$ independently to the right-hand side of \eqref{LLN-conv}; see  Remark~\ref{RemLLN}. However, in the cases specified in Assumption~\ref{AssCLT}2, the correlation between two such increments at different locations decays in general at a slower rate than $\sqrt{\del}$. So for Theorem~\ref{CLT} to hold, we must assume in these cases that each coordinate of $f$ uses increments at no more than one measurement site. The symmetry assumption on $f$ is standard and already needed in the semimartingale context in order to avoid an asymptotic bias; see Theorem~5.3.6 in \citeA{Jacod12}. 
	\item Assumption~\ref{AssCLT}3 on the temporal regularity of $\si$ is the ``union'' of two typical cases considered in the literature. The part $\si^{(0)}$ is (essentially) Hölder continuous of order strictly larger than $\frac12$, as considered, for instance, in \citeA{BN11,BN13,Corcuera13}. If one wants to include volatility processes that are of the roughness of Brownian motion, one has to make further structural assumptions as in \citeA{AitSahalia14,Jacod12},  namely, that $\si^{(1)}=\si-\si^{(0)}$ is a semimartingale. As we will show in Lemma~\ref{takesform}, \eqref{si-semi} is one possibility to obtain such a semimartingale structure, jointly in $x$.
	\item The volatility process must also have nice regularity in space. In fact, we assume that $\si$ is pathwise twice differentiable in space. By Theorem~\ref{LLN}, the limit of variation functionals taken at one measurement site only depends on the volatility at this site. However, this spatial concentration at the origin, which is due to the properties of the heat kernel, is very weak, so we must use the differentiability assumption \emph{and} the symmetry of the heat kernel (notice, however, Remark~\ref{symm}) to obtain a localization at a faster rate than $\sqrt{\del}$. For more details, we refer the reader to Remark~\ref{rem-discr-2} in the supplement.
	Examples of volatility models for $\si^{(1)}$ include  Ornstein--Uhlenbeck processes in space and time (see \citeA{BN04,Nguyen15}) and their generalizations (see \citeA{Brockwell17,Pham17}). If $\si$ does not depend on $x$ as in \citeA{Bibinger17}, then \eqref{si-der} is clearly satisfied.
\eenu
\erem

\brem\label{gen} In principle, the results and techniques developed in this paper apply to more general equations than \eqref{SHE}, or more general kernels $G$ in \eqref{Yx} and spatial covariance functions $F$ of the noise (for the existence of a Gaussian noise with $\La(\dd y,\dd z) = F(z-y)\,\dd y\,\dd z$, it is necessary and sufficient that the function $F$ be the Fourier transform of a nonnegative tempered measure on $\R^d$; see Section~2 of \citeA{Dalang99} for more details). 

As the proof shows, the law of large numbers (Theorem~\ref{LLN}) continues to hold as long as $G\in \call(W)$ has a dominating singularity at the origin such that, for some $\Ga_r\in\R$, the measures defined in \eqref{Pin} with the new $G$ and $F$ satisfy $\Pi^n_{r,0}\limw \Ga_r\delta_0$ and $|\Pi^n_{r,h}|\limw 0$ for all $r=0,1,\ldots$ and $h\neq0$.

For the central limit theorem (Theorem~\ref{CLT}), we additionally need the following assumptions [expressed in terms of $\Pi^n_{r,h}$ and $|\Pi^n_{r,h}|$ from \eqref{Pin}]:
\benu
\item $G$ is symmetric in the sense that $G(t,x)=G(t,-x)$ for any $t>0$ and $x\in\R^d$. In addition, for every $r\in\N_0$, one has, as $n\to\infty$,
\[ \iiint_0^\infty (|y|^2+|z|^2)\,|\Pi^n_{r,0}|(\dd s,\dd y,\dd z) = o(\delh). \]
\item For every $r\in\N_0$, we have, as $n\to\infty$, $$|\Pi^n_{r,0}([0,\infty)\times\R^d\times\R^d)-\Ga_r| = o(\sqrt{\del}).$$ Moreover, either $f$ satisfies Assumption~\ref{AssCLT}2, or for all $r\in\N_0$ and $h\neq0$, $$|\Pi^n_{r,h}|([0,\infty)\times\R^d\times\R^d)=o(\sqrt{\del}).$$ 
\item  There is some decreasing square-summable sequence $(\ov\Ga_r\colon r\in\N_0)$ such that for all $n\in\N$, $r\in\N_0$, and $h\in\R^d$,  $$|\Pi^n_{r,h}|([0,\infty)\times\R^d\times\R^d) \leq \ov \Ga_r.$$ 
\item There is $\nu>1$ such that for all $\theta\in(0,1)$,  $$|\Pi^n_{0,0}|((\Delta_n^{1-\theta},\infty)\times\R^d\times\R^d) = O(\Delta_n^{\nu \theta}).$$
\eenu
Condition (1) is needed for the reasons explained in part (4) of Remark~\ref{rem-discr} but can be relaxed; see Remark~\ref{symm}. Condition~(2) is used for the terms $K^{n,i}_3$, $K^{n,i}_5$, and $K^{n,i}_6$ in the proof of Lemma~\ref{remove-s}. Condition (3) is crucial for the actual central limit theorem in Proposition~\ref{CLT-core} (if the asymptotic covariances of the increments fail to be square-summable, there is no hope to see a central limit theorem; cf.\ \citeA{Nourdin10}). Finally, condition (4) is needed for nearly all approximations in the proof. 

If the kernel has singularity fronts as, for example, in the case of the wave equation, the limits in \eqref{Vnf} will have a different shape, which we shall discuss in a separate work. 
\erem

\brem\label{symm} In the setting of Remark~\ref{gen}, the symmetry assumption on $G$ can be weakened. Suppose that $G(t,x)=\wt G(t,x)H(x)$ where $\wt G$ satisfies condition (1) in Remark~\ref{gen}, and $H\colon \R^d \to \R$ is differentiable such that 
\beq\label{H0} \sup_{i=1,\ldots,d} \sup_{u\in[0,1]} \lv\frac{\partial}{\partial x_i}H(ux)\rv+ \sup_{u\in[0,1]} |H(ux)|  \leq H_0(x) \eeq
for some function $H_0\colon \R^d \to [0,\infty)$.

 In the proof of Theorem~\ref{CLT}, the symmetry of $G$ is only used to show the identity \eqref{important}. In the asymmetric case, we observe from \eqref{Pin} that $\Pi^n_{r,0}(\dd s,\dd y,\dd z) = H(y)H(z)\,\wt\Pi^n_{r,0}(\dd s,\dd y,\dd z)$ for all $r\in\N_0$, where $\wt\Pi^n_{r,h}$ is the measure that arises from the first equation in \eqref{Pin} when we replace $G$ by $\wt G$. By the mean value theorem, applied to $(y,z)\mapsto H(y)H(z)$, and property \eqref{H0}, the left-hand side of \eqref{important} is bounded by a constant times
 \begin{align*}
&H(0)^2\lv\iiint_0^{(\la_n+(k\vee k'))\del} (y_l+z_l) \,\wt \Pi^n_{|k'-k|,0}(\dd s,\dd y,\dd z)\rv\\
&\quad + \iiint_0^\infty |y_l+z_l|H_0(y)H_0(z)(|y|+|z|)\,|\wt \Pi^n_{|k'-k|,0}|(\dd s,\dd y,\dd z)  .
 \end{align*}
 The first term is zero by symmetry. Thus, if we impose the condition
 \beq\label{newcond} \iiint_0^\infty (|y|^2+|z|^2)H_0(y)H_0(z)\,|\wt \Pi^n_{r,0}|(\dd s,\dd y,\dd z)=o(\sqrt{\Del}) \eeq
as $n\to\infty$ for all $r\in\N_0$, we no longer need $G$ to be symmetric in Remark~\ref{gen}.
 
 Let us apply the previous discussion to the important example where $\wt G$ is the heat kernel \eqref{heatkernel}, and $H(x)=\ee^{\theta\cdot x}$ for some $\theta \in\R^d$, which corresponds to \eqref{SHE} with an additional gradient term (considered, in similar forms, by \citeA{Bibinger17} and many of the applications mentioned in the introduction). It is easily verified that the resulting kernel $G$ belongs to $\call(W)$ (and hence, a stationary mild solution exists in the case of constant $\si$) if and only if $\la>\frac\kappa 2 |\theta|^2$. 
 Under this additional constraint, one can then show that \eqref{newcond} holds with $H_0(x)=(1+|\theta|)(\ee^{\theta\cdot x}\vee1)$ in \eqref{H0}. Indeed, by symmetry considerations, the left-hand side of \eqref{newcond} is bounded by a constant times
 \[  \iiint_0^\infty (|y|^2+|z|^2)\ee^{\theta\cdot y}\ee^{\theta\cdot z}\,|\wt \Pi^n_{r,0}|(\dd s,\dd y,\dd z), \]
 which equals the left-hand side of \eqref{pin2}.
 Using the identities $G(t,x)=\wt G(t,x) \ee^{\theta\cdot x}$ and $G(t,x)=\wt G(t,x-\kappa \theta t)\ee^{\kappa |\theta|t/2}$, it is not difficult to see that \eqref{cross-bound-3} and \eqref{dens} remain valid, as well as \eqref{int-bound-2}, \eqref{convolution}, and \eqref{der} if we replace $\la$ by $\la_0=\la-\frac\kappa2|\theta|$. One can now follow the arguments given in the proof of Lemma~\ref{y2} in order to show that \eqref{pin2} and hence \eqref{newcond} hold true. 
\erem

\brem\label{non-stat} All results in this section remain valid if we consider \eqref{SHE} for $t>0$ and $x\in\R^d$, subject to some bounded and sufficiently regular initial condition $y_0$ at time $t=0$. Indeed, the mild solution is then given by
\beq\label{mild-2} Y(t,x)=\int_{\R^d} G(t,x-y)y_0(y)\,\dd y + \iint_0^t G(t-s,x-y)\si(s,y)\,W(\dd s,\dd y) \eeq
for $(t,x)\in(0,\infty)\times\R^d$.
Let us denote the first and the second term by $Y^{(0)}$ and $Y^{(1)}$, respectively, and fix $T>0$. Under the hypothesis that $y_0$ is  Hölder continuous with some exponent $>1-\frac\al 2$ (resp., differentiable with a derivative that is  Hölder continuous with some exponent $>1-\frac\al 2$), we know from classical PDE theory (see Theorem~5.1.2(ii) in \citeA{Lunardi95}) that $t\mapsto Y^{(0)}(t,x)$ is Hölder continuous on $[0,T]$ with some exponent $\eta>\frac12-\frac\al 4$ (resp., $\eta>1-\frac\al 4$), uniformly for $x\in\R^d$. In particular, by \eqref{taun-asymp}, we have $|\Delta^n_i Y^{(0)}_x/\tau_n| \lec \Delta_n^{\eta}/\tau_n\lec \Delta_n^{\eta- 1/2 + \al/4}$, where the last exponent is strictly positive (resp., larger than $\frac12$). From this, it is straightforward to deduce that the contribution of $Y^{(0)}$ to \eqref{Vnf}  is asymptotically negligible in Theorem~\ref{LLN} (resp., Theorem~\ref{CLT}). 

We are left to show that $Y^{(1)}$ has the same asymptotic behavior as the expression in \eqref{Yx}. For the law of large numbers, this is straightforward, while for the central limit theorem, it can be proved analogously to Step 1 in Section~\ref{overview}. The details are omitted at this point. We further remark
that the assumption $\la>0$ is superfluous when \eqref{SHE} is considered for $t\geq0$, and it is sufficient to formulate  Assumptions~\ref{AssLLN}2 and \ref{AssLLN}3 as well as Assumption~\ref{AssCLT}3 with $t\in[0,T]$ for any $T>0$ instead of $t\in\R$, and we may replace $-\infty$ in \eqref{si-semi} by $0$. 
\erem

\brem In the literature of stochastic PDEs, one often considers equations where the random field $\si$ is an explicit functional of the solution $Y$, that is, equation \eqref{mild-2} where $\si=B(Y)$ and $B$ is an operator satisfying certain regularity and growth conditions such that \eqref{mild-2}  admits a mild solution; see \citeA{DaPrato92}. While Assumption~\ref{AssLLN} is relatively weak in this situation (for example, it is satisfied if
\beq\label{3cases} \begin{split} B(Y)(t,x)&=b(Y(t,x)),\\
	 B(Y)(t,x)&=\iint_0^t  H(t-s,x-y)Y(s,y)\,\dd y\,\dd s\\ \text{or}\qquad B(Y)(t,x)&=\iint_0^t  K(t-s,x-y)Y(s,y)\,W(\dd s,\dd y),\end{split}\eeq
and $b$ is Lipschitz continuous, $H\in L^1([0,T]\times\R^d)$, and $K\bone_{[0,T]}\in \call(W)$ for all $T>0$), Assumption~\ref{AssCLT} is more restrictive. In fact, for the functionals in \eqref{3cases}, it is only satisfied if $b$ is constant, or if $B(Y)$ is the second or third expression with functions $H$ and $K$ that are sufficiently smooth.
\erem

\subsection{Multipower variations} \label{multipower} We apply Theorems~\ref{LLN} and \ref{CLT} to an important class of functionals, namely to so-called \emph{multipower variations} $V^n_\Phi(Y_x,t)$ or \emph{signed multipower variations} $V^n_\Psi(Y_x,t)$, where $\Phi,\Psi\colon \R^{N\times L}\to \R^N$ (note that $N=M$) are given by
\begin{align} \label{mpv}  \Phi_m(z) &=\Phi_m\lp (z_{jk})_{j,k=1}^{N,L}\rp= \prod_{k=1}^L |z_{mk}|^{w_{mk}},\\
\label{smpv} \Psi_m(z)&=\Psi_m\lp (z_{jk})_{j,k=1}^{N,L}\rp= \prod_{k=1}^L (z_{mk})^{w_{mk}},\qquad m=1,\ldots,N,\end{align}
with $w_{mk}\geq0$ in \eqref{mpv} and $w_{mk}\in \N_0$ in \eqref{smpv}. We shall write $w=(w_{mk})_{m,k=1}^{N,L}$, $w_m = w_{m1}+\cdots+w_{mL}$, and $\ov w = \max\{w_1,\ldots,w_N\}$. If we want to emphasize the dependence on $w$, we write $\Phi(z)=\Phi(w;z)$ and $\Psi(z)=\Psi(w;z)$. If $w_{mk} = p_k$ for all $m$ and $k$, we write
\beq\label{equal-m} \Phi(p_1,\ldots,p_L;z)=\Phi(w;z)\qquad\text{and}\qquad  \Psi(p_1,\ldots,p_L;x) = \Psi(w;z). \eeq

For multipowers, Theorems~\ref{LLN} and \ref{CLT} take the following form: 
\bcor\label{multipower-result} Assume that $\al\in(0,2)\cap (0,d]$.
\benu
\item If Assumptions~\ref{AssLLN}2 and \ref{AssLLN}3 hold with $p=\ov w\vee 2$, then we have for all $m=1,\ldots,N$,
\begin{align}\label{LLN-mpv} V^n_{\Phi_m|\Psi_m}(Y_x,t) &\limL \mu_{\Phi_m|\Psi_m} \int_0^t |\si(s,x_m)|^{w_m}\,\dd s,\end{align}
where $\mu_{f}=\mu_{f}(1,\ldots,1)$ (as defined in Theorem~\ref{LLN}) and $\Phi_m|\Psi_m$ means that we can either take $\Phi_m$ or $\Psi_m$ in \eqref{LLN-mpv}. Note that $\mu_{\Psi_m}=0$ if $w_m$ is odd.
\item Suppose that $w_{mk}\in\{0,2\}$ or $w_{mk}\geq 4$ in the case of \eqref{mpv}, and that all $w_m$ are even in the case of \eqref{smpv}. Further assume that Assumption~\ref{AssCLT}3 holds with $p=\ov w$. Then \eqref{CLT-statement} holds for $f=\Phi|\Psi$, and the $\calf$-conditional covariance processes in \eqref{calc} are given by
\beq\label{cov-mpv}
\calc_{m_1 m_2}(t)=\begin{cases}  \rho_{\Phi_m|\Psi_m}\displaystyle\int_0^t |\si(s,x_m)|^{2w_m}\,\dd s, &m_1=m_2=m,\\ 0,&m_1\neq m_2,\end{cases}\eeq
where $\rho_{f} = \rho_{f,f}(0;1,\ldots,1)+2\sum_{r=1}^\infty \rho_{f,f}(r;1,\ldots,1)$ (as defined in Theorem~\ref{CLT}).
\eenu
\ecor

Because of their particular importance in high-frequency statistics, we further specialize Corollary~\ref{multipower-result} to the normalized power variations 
\beq V^n_p(Y_x,t)=(V^n_p(Y_x,t)_m)_{m=1}^N = \lp \Del\sum_{i=1}^{[t/\Del]} \lv \frac{\Delta^n_i Y(\cdot,x_m)}{\tau_n}\rv^p \rp_{m=1}^N, \eeq
where $p>0$ and which corresponds to the special case $L=1$ and the function $\Phi(p,\cdot)$ in \eqref{equal-m}.
\bcor\label{pvar-cor} Assume that $\al\in(0,2)\cap(0,d]$.
\benu
	\item Let $p>0$ and $\si$ be a predictable random field satisfying Assumptions~\ref{AssLLN}2 (with $p\vee 2$ instead of $p$) and \ref{AssLLN}3. Then for every $m=1,\ldots,N$,
	\beq V^n_p(Y_x,t)_m \limL \mu_p\int_0^t |\si(s,x_m)|^p\,\dd s, \eeq
	where $\mu_p=\bbe[|Z|^p]$ with $Z\sim N(0,1)$. 
	\item Let $p=2$ or $p\geq4$ and suppose that Assumption~\ref{AssCLT}3 holds. Then
	\beq \lp\Del^{-\frac12}\lp V^n_p(Y_x,t)_m-\mu_p\int_0^t |\si(s,x_m)|^p\,\dd s\rp\rp_{m=1}^N \limst \calz, \eeq
	where $\calz$ is a process as described after \eqref{CLT-statement} and
	\beq \calc_{m_1m_2}(t) = \begin{cases} R_p \displaystyle \int_0^t |\si(s,x_m)|^{2p}\,\dd s, &m_1=m_2=m,\\ 0,&m_1\neq m_2. \end{cases} \eeq
	In the previous line,
$R_p=\rho_p(1)+2\sum_{r=1}^\infty\rho_p\left(\Ga_r\right)$, where $\Ga_r$ is defined in \eqref{Ga-formula},
	and $\rho_p(r)= \cov(|X|^p,|Y|^p)$ for $(X,Y)\sim N\big(0,\big(\begin{smallmatrix} 1 & r\\ r & 1\end{smallmatrix}\big)\big)$.
\eenu
\ecor

\bex If $\dot W$ is a space-time white noise and $p=2$ or $p=4$, by expressing $x^2$ and $x^4$ in terms of Hermite polynomials and then using Lemma~1.1.1 in \citeA{Nualart06}, we obtain
\beq \begin{split} R_2 &= 2+4\sum_{r=1}^\infty \left(\frac12 \sqrt{r+1}-\sqrt{r}+\frac12 \sqrt{r-1}\right)^2 = 2.357487..., \\
R_4 &= 96+144\sum_{r=1}^\infty \left(\frac12 \sqrt{r+1}-\sqrt{r}+\frac12 \sqrt{r-1}\right)^2\\
&\quad+48 \sum_{r=1}^\infty \left(\frac12 \sqrt{r+1}-\sqrt{r}+\frac12 \sqrt{r-1}\right)^4 = 109.223069...,  \end{split}\eeq
which are larger than the corresponding constants $2$ and $96$ in the semimartingale framework (cf.\ Theorem~6.1 and Example~6.5 in \citeA{AitSahalia14}) and are the same as in the setting of moving average processes (cf.\ Theorem~4 
in \citeA{BN11}). The reason for larger constants compared to the semimartingale case is the nonvanishing asymptotic correlation between increments of $Y(\cdot,x_m)$.
Let us also remark that $R_2=\pi\Ga$ for the constant $\Ga$ in Theorem~4.2 of \citeA{Bibinger17}, and that $R_4=\check \si^2$ for the constant $\check\si^2$ in Equation~(A.2) of \citeA{Cialenco17}. 
\eex

\subsection{Estimation of volatility and spatial noise correlation index}\label{est}

In this section, we will explain how Theorems~\ref{LLN} and \ref{CLT} can be applied to estimate the volatility process $\si$ and the spatial correlation index $\al$ of the noise. 

For both problems, the knowledge of the parameter $\la$ is irrelevant as we shall see. This is important because there is no way to estimate $\la$ consistently under our observation scheme. Indeed, a Girsanov argument (see Proposition~1.6 in \citeA{Nualart94}) shows that for constant $\si$, the laws of the solution $Y$ on a compact space-time set are equivalent for different values of $\lambda$. 

 Furthermore, for the estimation of $\si$, we will assume that the parameter $\kappa$ is known. In fact, if $N=1$ and $\si$ is a constant, then $Y(\cdot,x_1)$ is a stationary Gaussian process whose distribution is completely determined by its covariance function. By \eqref{Rt} and the scaling properties of the normal distribution, this only depends on the ratio $\si^2/\kappa^{\al/ 2}$, so there is no way to identify the pair $(\si,\kappa)$ based on observations of $Y(\cdot,x_1)$. If measurements are recorded at $N\geq2$ spatial positions, then by the second statement of Lemma~\ref{delta} (1), the normalized increments at different locations are asymptotically uncorrelated, and hence independent. Thus, it is impossible to consistently estimate both $\kappa$ and $\si$ based on observations at finitely many space points. 

Despite this restriction, even if $\kappa$ is unknown, the subsequent results can be easily modified to yield consistent and asymptotically mixed normal estimators for the \emph{viscosity-adjusted volatility} $\kappa^{-\al w_m/4}\int_0^t |\si(s,x_m)|^{w_m}\, \dd s$ or the \emph{relative volatility} $\int_0^t |\si(s,x_m)|^{w_m}\, \dd s/\int_0^T |\si(s,x_m)|^{w_m}\, \dd s$ with $m=1,\ldots,N$, $w_m$ as below, and $t\in[0,T]$. Both quantities are constant multiples of the integrated volatility and thus completely describe the shape of the temporal fluctuations of $\si$, which is sufficient for many applications; see, for example, \citeA{BN14}, where the concept of relative volatility was introduced and further applied to turbulence data.

We first consider the situation when the spatial correlation index $\al$ is \emph{known}. Then Corollary~\ref{multipower-result} immediately yields consistent estimators and asymptotic confidence bounds for the integrated volatility process at the measurement sites $x_1,\ldots,x_N$.
In the theorems of this section, we will often divide by asymptotic $\calf$-conditional variances during studentization procedures which may be zero in some degenerate situations. In these cases, convergence in probability and stable convergence in law should be understood in restriction to the set where all involved realized variation functionals are strictly positive. For the theoretical background of this concept for stable convergence in law, we refer the reader to Chapter~3.2 in \citeA{AitSahalia14} and to \citeA{Podolskij10}.

\bthm\label{known-alpha} Assume that $\al\in(0,2)\cap (0,d]$ and $\kappa>0$ are known. Define $\wt V^n_{\Phi|\Psi}(Y_x,t)$ in the same way as $V^n_{\Phi|\Psi}(Y_x,t)$ but with $\tau_n$ replaced by
\beq\label{tildetau} \wt\tau_n^2 =  \displaystyle\frac{\pi^{\frac d 2-\al}\Ga(\frac\al 2)}{(2\kappa)^{\frac \al 2}(1-\frac \al 2)\Ga(\frac d 2)} \del^{1-\frac \al 2}. \eeq 
Then, under the hypotheses of Corollary~\ref{multipower-result} (1) and (2), we have 
\begin{align}\label{LLN-mpv-2} \wt V^n_{\Phi_m|\Psi_m}(Y_x,t) &\limL \mu_{\Phi_m|\Psi_m} \int_0^t |\si(s,x_m)|^{w_m}\,\dd s,\qquad m=1,\ldots,N,\end{align}
and, for every $T>0$,
\beq\label{stud-1-2}\begin{split} &\left\{\del^{-\frac12}\frac{\mu_{\Phi_m|\Psi_m}}{\sqrt{\rho_{\Phi_m|\Psi_m}}}\sqrt{\frac{\mu_{\Phi_m|\Psi_m(2w;\cdot)}}{\wt V^n_{\Phi_m|\Psi_m(2w;\cdot)}(Y_x,T)}} \vphantom{\lp \frac{\wt V^{n}_{\Phi_m|\Psi_m}(Y_x,T)}{\mu_{\Phi_m|\Psi_m}} \rp} \right.\\
	&\left.\quad \times \lp \frac{\wt V^{n}_{\Phi_m|\Psi_m}(Y_x,T)}{\mu_{\Phi_m|\Psi_m}} - \int_0^T |\si(s,x_m)|^{w_m}\,\dd s\rp\right\}_{m=1}^N \lims N(0,\mathrm{Id}_N) \end{split}\eeq
as $n\to\infty$. The left-hand sides of \eqref{LLN-mpv-2} and \eqref{stud-1-2} are independent of $\la$.
\ethm

If $\al$ is unknown, we first have to find a consistent estimator for $\al$, for which we propose two solutions. The first estimator is a regression-type estimator similar to the change-of-frequency estimator in \citeA{BN13,Corcuera13} for the kernel singularity of a moving average process and similar to the estimator for the Hölder index of a Gaussian process proposed in \citeA{Istas97}. Define the function $\Phi^{(2)}(p;\cdot)\colon \R^{N\times 2}\to\R^N$ by
\begin{align*} \Phi^{(2)}_m(p;x)&=\Phi^{(2)}_m\lp p; (x_{j1},x_{j2})_{j=1}^{N}\rp = |x_{m1}+x_{m2}|^p, \qquad m=1,\ldots,N.\end{align*}
Furthermore, recalling $\rho_p(r)$ from Corollary~\ref{pvar-cor}, define
\beq\label{C0alpha} \calc_0(\al)= \lp\frac{4}{p\log 2}\rp^2 \lp \calc_{11}-2 \frac{\calc_{12}}{(2+2\Ga_1)^{\frac p2}}+\frac{\calc_{22}}{(2+2\Ga_1)^p} \rp, \eeq
where 
\beq \label{Cs}
\begin{split}
	\calc_{11}&= \rho_p(1)+2\sum_{r=1}^\infty \rho_p(\Ga_r),\\
	\calc_{22}&=(2+2\Ga_1)^p\lp \rho_p(1)+2\sum_{r=1}^\infty \rho_p\lp\frac{2\Ga_r+\Ga_{r-1}+\Ga_{r+1}}{2+2\Ga_1}\rp\rp,\\
	\calc_{12}&=(2+2\Ga_1)^{\frac p 2}\lp \rho_p\lp \sqrt{\frac{1+\Ga_1}{2}}\rp + \sum_{r=1}^\infty \rho_p\lp \frac{\Ga_r+\Ga_{r-1}}{\sqrt{2+2\Ga_1}}\rp\right.\\
	&\quad\left.+\vphantom{\rho_p\lp \sqrt{\frac{1+\Ga_1}{2}}\rp} \sum_{r=1}^\infty \rho_p\lp \frac{\Ga_r+\Ga_{r+1}}{\sqrt{2+2\Ga_1}}\rp\rp .
\end{split}\eeq
Note that $\calc_0$ depends on $\alpha$ via $\Ga_r$; see \eqref{Ga-formula}.
\bthm\label{est-alpha}
Let $\al\in(0,2)\cap (0,d]$.
\benu
\item Suppose that $p>0$ and Assumptions~\ref{AssLLN}2 and \ref{AssLLN}3 hold with exponent $p\vee 2$. Then as $n\to\infty$,
\beq\label{alpha-conv}\begin{split}\wh \al^{(p)}_n&=\frac1N \sum_{m=1}^{N} \wh \al^{(p),m}_n =\frac1N \sum_{m=1}^N\lp 2-\frac{4}{p}\log_2 \lp 
\frac{V^n_{\Phi^{(2)}_m(p;\cdot)}(Y_x,T)}{V^n_{\Phi_m(p;\cdot)}(Y_x,T)}
\rp\rp\\ &= 2-\frac{4}{pN}\sum_{m=1}^N\log_2 \lp 
\frac{\sum_{i=1}^{[T/\del]-1} {|\Delta^n_i Y_{x_m} +\Delta^n_{i+1} Y_{x_m}|}^p}{\sum_{i=1}^{[T/\del]} {|\Delta^n_i Y_{x_m}|}^p}
\rp \limp \al,\end{split}\raisetag{-3\baselineskip}
\eeq 
\item Suppose that $p=2$ or $p\geq 4$ and that Assumption~\ref{AssCLT}3 holds for this value of $p$. Then
\beq\label{conf-int}\begin{split} &\frac{N}{\del^{\frac12}} \sqrt{\frac{2^p\Ga(\frac{2p+1}{2})}{\pi^{\frac12}\calc_0(\wh\al^{(p)}_n) }} \lp \sum_{m=1}^N \frac{V^n_{\Phi_m(2p;\cdot)}(Y_x,T)}{(V^n_{\Phi_m(p;\cdot)}(Y_x,T))^2}\rp^{-\frac12} (\wh\al^{(p)}_n -\al)\\
	&\quad \lims N(0,1).\end{split} \eeq
\eenu
Note that the left-hand sides of \eqref{alpha-conv} and \eqref{conf-int} do not depend on the parameters $\kappa$ and $\lambda$.
\ethm

The second estimator is a correlation estimator (compare with the modified realized variation ratio of \citeA{BN13}). To this end, we 
define
\beq\label{C0tilde} \wt \calc_0(\al)=\lp\frac{2}{\log 2}\rp^2(\wt C_{11}-2\wt C_{12}\Ga_1+\wt C_{22}\Ga_1^2), \eeq
where
\begin{align*}
		\wt\calc_{11}&= 1+\Ga_1^2+2\sum_{r=1}^\infty (\Ga_r^2+\Ga_{r+1}\Ga_{r-1}),\qquad	\wt\calc_{22}=2+4\sum_{r=1}^\infty \Ga_r^2,\\
		\wt\calc_{12}&=2\Ga_1+2\sum_{r=1}^\infty \Ga_r(\Ga_{r+1}+\Ga_{r-1}).
\end{align*}

\bthm\label{correst}
Let $\al\in(0,2)\cap(0,d]$ and $F(x)=-2\log_2(1+x)$. 
\benu
\item Under  Assumptions~\ref{AssLLN}2 and \ref{AssLLN}3 with $p=2$, we have as $n\to\infty$,
\beq\label{tilde-alpha} \wt\al_n = \frac1N\sum_{m=1}^N \wt\al^{m}_n =  \frac1N \sum_{m=1}^N F\lp\frac{V^n_{\Psi_m(1,1;\cdot)}(Y_x,T)}{V^n_{\Phi_m(2;\cdot)}(Y_x,T)}\rp\limp \al.  \eeq
\item Under Assumption~\ref{AssCLT}3 with $p=2$, we have as $n\to\infty$,
\beq\label{tilde-alpha-CLT}\begin{split} &\frac{N}{\del^{\frac12}}\sqrt{\frac{3}{\wt \calc_0(\wt \al_n)}} \lp  \sum_{m=1}^N \frac{V^n_{\Phi_m(4;\cdot)}(Y_x,T)}{(V^n_{\Psi_m(1,1;\cdot)+\Phi_m(2;\cdot)}(Y_x,T))^2}\rp^{-\frac12}(\wt\al_n-\al)\\ &\quad\lims N(0,1).\end{split}  \eeq
\eenu
Both quantities on the left-hand side of \eqref{tilde-alpha} and \eqref{tilde-alpha-CLT} do not depend on $\kappa$ and $\lambda$.
\ethm

\brem Let us compare the asymptotic variances of the two estimators $\wh\al^{(p)}_n$ and $\wt \al_n$ in the case where $\si(t,x) \equiv \si$ is constant (but nonzero). If $p=2$, then under the assumptions of Theorems~\ref{est-alpha} and \ref{correst},
\beq\begin{aligned}
 \lim_{n\to\infty}	\Del^{-\frac12}\var[\wh\al^{(2)}_n-\al] &=\frac{1}{NT} \calc_0(\al),\\
\lim_{n\to\infty}	\Del^{-\frac12}\var[\wt\al_n-\al] &=\frac{1}{NT} \frac{\wt \calc_0(\al)}{(1+\Ga_1)^2}.
\end{aligned}\eeq
A straightforward calculation shows that $\calc_0(\al)=\wt \calc_0(\al)/(1+\Ga_1)^2$ for all $\al\in(0,2)$. In other words, from the viewpoint of asymptotic variance, the two estimators $\wh\al^{(2)}_n$ and $\wt \al_n$ are equivalent. 
By varying the value of $p$ for the estimator $\wh\al^{(p)}_n$, we can further check how reliable the estimates for $\al$ are.
\erem

With (either of) the two estimators for $\al$ at hand, we can now proceed to the estimation of $\si$. The rate of convergence is slower by a logarithmic factor compared to the case where $\al$ is known; see Theorem~\ref{known-alpha}. This is the same phenomenon that occurs when the smoothness and the variance of a Gaussian process are to be estimated at the same time; see \citeA{Brouste18,Fukasawa18,Istas97}.

\bthm\label{unknown-alpha} Assume that $\al\in(0,2)\cap(0,d]$ and that $\kappa$ is known. Let $\al_n = \wh\al^{(p_0)}_n$ for some $p_0>0$ or $\al_n=\wt \al_n$, in which case we set $p_0=2$. 
Define
\beq\label{taunhat} \wh\tau_n^2 =  \displaystyle\frac{\pi^{\frac d 2-\al_n}\Ga(\frac{\al_n}  2)}{(2\kappa)^{\frac {\al_n}  2}(1-\frac {\al_n}  2)\Ga(\frac d 2)} \del^{1-\frac {\al_n} 2},\qquad n\in\N,\eeq
and the functionals $\wh V^n_{\Phi|\Psi}(Y_x,t)$ in the same way as $V^n_{\Phi|\Psi}(Y_x,t)$ but with $\wh\tau_n$ instead of $\tau_n$.
\benu
	\item If Assumptions~\ref{AssLLN}2 and \ref{AssLLN}3 hold with $p_0\vee \ov w \vee 2$, then for every $T\geq0$ and $m=1,\ldots,N$,
	\[ \wh V^n_{\Phi_m|\Psi_m}(Y_x,T) \limp \mu_{\Phi_m|\Psi_m}\int_0^T |\si(s,x_m)|^{w_m}\,\dd s\as. \]
	\item Assume that all $w_{mk}\in \{0,2\}\cup[4,\infty)$ in the case of $\Phi$, and that all $w_m$ are even in the case of $\Psi$. Also assume that $p_0=2$ or $\geq 4$, and that Assumption~\ref{AssCLT}3 holds with exponent $p_0\vee \ov w$. Then, if $Z\sim N(0,1)$, we have as $n\to\infty$, in the case $\al_n=\wh\al_n^{(p_0)}$,
	\begin{align*} &\Bigg\{\frac{\del^{-\frac12}}{|\log\del|} \frac{4N}{w_m \wh V^n_{\Phi_m|\Psi_m}(Y_x,T)}\\
	&\quad\times \sqrt{\frac{2^{p_0}\Ga(\frac{2p_0+1}{2})}{\pi^{\frac12}\calc_0(\wh\al^{(p_0)}_n) }}\lp \sum_{j=1}^N  \frac{V^n_{\Phi_j(2p_0;\cdot)}(Y_x,T)}{( V^n_{\Phi_j(p_0;\cdot)}(Y_x,T))^2}  \rp^{-\frac12}\\
		&\quad\times \Bigg(\wh V^n_{\Phi_m|\Psi_m}(Y_x,T)-\mu_{\Phi_m|\Psi_m}\int_0^T |\si(s,x_m)|^{w_m}\,\dd s\Bigg)\Bigg\}_{m=1}^N \lims \begin{pmatrix} Z\\ \vdots\\ Z\end{pmatrix},
		\end{align*}
	while in the case $\al_n=\wt\al_n$,
	\begin{align*} &\Bigg\{\frac{\del^{-\frac12}}{|\log\del|} \frac{4N}{w_m \wh V^n_{\Phi_m|\Psi_m}(Y_x,T)}\\
		&\quad\times\sqrt{\frac{3}{\wt \calc_0(\wt \al_n)}}  \lp \sum_{j=1}^N  \frac{V^n_{\Phi_j(4;\cdot)}(Y_x,T)}{( V^n_{\Psi_j(1,1;\cdot)+\Phi_j(2,\cdot)}(Y_x,T))^2}  \rp^{-\frac12}\\
		&\quad\times \Bigg(\wh V^n_{\Phi_m|\Psi_m}(Y_x,T)-\mu_{\Phi_m|\Psi_m}\int_0^T |\si(s,x_m)|^{w_m}\,\dd s\Bigg)\Bigg\}_{m=1}^N \lims \left(\begin{matrix} Z\\ \vdots\\ Z\end{matrix}\right).
	\end{align*}
\eenu
\ethm

\brem Whereas different coordinates are asymptotically independent in \eqref{stud-1-2}, they are identical in the limit in Theorem~\ref{unknown-alpha} (2). The former is a consequence of \eqref{cov-mpv}, while the latter is due the fact that the dominating term in the case of unknown $\al$ comes from the difference $\al_n-\al$ (which is independent of $m$); see the proof of Theorem~\ref{unknown-alpha}.
\erem

\section{Overview of proofs}\label{proofs}

The main ideas for the proof of Theorem~\ref{LLN} and \ref{CLT} are sketched in this section, while the details will be given in Sections~\ref{details-LLN} and \ref{details}. The results of Sections~\ref{multipower} and \ref{est} are shown in Section~\ref{proof-rest}.

Without risk of confusion, we shall use the notation $x=(x_1,\ldots,x_N)'\in(\R^d)^N$, $\un\Delta^n_i Y_x = (\Delta^n_i Y_{x},\ldots,\Delta^n_{i+L-1} Y_x)\in\R^{N\times L}$, and 
\begin{align*}
\Delta^n_i G_y(s)&=G(i\del-s,y)-G((i-1)\del-s,y) \in \R,\\
\Delta^n_i G_{x,y}(s) &= (\Delta^n_i G_{x_1-y}(s),\ldots,\Delta^n_i G_{x_N-y}(s))'\in \R^N,\\
\un\Delta^n_i G_{x,y}(s)&=(\Delta^n_i G_{x,y}(s),\ldots,\Delta^n_{i+L-1} G_{x,y}(s))\in\R^{N\times L}
\end{align*}
  for $s\in\R$ and $y\in\R^d$. Similarly, $\un\si(s,x) = (\si(s,x_1),\ldots,\si(s,x_N))'$, and $\un\si^2(s,x) = (\si^2(s,x_1),\ldots,\si^2(s,x_N))'$. Moreover, we write $t^*_n = [t/\Del]-L+1$ for $t\in[0,\infty)$ and $A\lec B$ if $A\leq CB$ for some finite constant $C>0$ that does not depend on any important parameter.

The following measures will play an important role in identifying the limit behavior of \eqref{Vnf}:
\beq\label{Pin}\begin{split} \Pi^n_{r,h}(A)
	&=\iiint_A\frac{ G_y(s)-G_y(s-\Delta_n)}{\tau_n}\\
	&\quad\times\frac{G_{z+h}(s+r\del)-G_{z+h}(s+(r-1)\del)}{\tau_n}\,\dd s\,\La(\dd y,\dd z),\\
	|\Pi^n_{r,h}|(A)
	&=\iiint_A\frac{ |G_y(s)-G_y(s-\Delta_n)|}{\tau_n}\\
	&\quad\times\frac{|G_{z+h}(s+r\del)-G_{z+h}(s+(r-1)\del)|}{\tau_n}\,\dd s\,\La(\dd y,\dd z),
\end{split}\eeq
where $r\in\N_0$, $h\in\R^d$, and $A\in \calb([0,\infty)\times\R^d\times \R^d)$. By \eqref{taun}, we have $\Pi^n_{0,0}([0,\infty)\times\R^d\times \R^d)=1$, so $\Pi^n_{0,0}$ is a probability measure. In fact, if we consider an arbitrary, say, the first component of the increment $\Delta^n_i Y_x$ with $\si\equiv1$, then for $A_1 \in \calb([0,\infty))$ and $A_2\in\calb(\R^d)$, $\Pi^n_{0,0}(A_1\times A_2\times A_2)$ is the proportion of the variance of $Y(i\del,x_1)-Y((i-1)\del,x_1)$ that is explained by the integral in \eqref{increment} on the set $\{(s,y)\colon (i\del-s,x_1-y) \in A_1\times A_2\}$. 

In general,  $\Pi^n_{r,h}([0,\infty)\times\R^d\times \R^d) \in [0,1]$ is the correlation between two increments $\Delta^n_i Y(\cdot,x_1)$ and $\Delta^n_{i+r}Y(\cdot,x_1+h)$, taken at a temporal distance of $r\del$ and a spatial distance of $h$. The value $\Pi^n_{r,h}(A_1\times A_2\times A_3)$, with $A_1$ and $A_2$ as above and $A_3\in\calb(\R^d)$, then quantifies how much this correlation is caused by the restrictions of the corresponding integrals in \eqref{increment} to the domains $\{(s,y)\colon (i\del-s,x_1-y) \in A_1\times A_2\}$ and $\{(s,z)\colon ((i+r)\del-s,x_1+h-z) \in A_1\times A_3\}$, respectively. Some important properties of these measures are proved in Section~\ref{heat-est}.

\subsection{Overview of the proof of Theorem~\ref{LLN}}

By arguing componentwise, we may assume without loss of generality that $M=1$.
As a first step, we show that we may further assume that $\sigma$ is a bounded random field.
\blem\label{sibounded} In order to prove Theorem~\ref{LLN}, one may additionally assume that $\si$ is uniformly bounded in $(\om,t,x)\in\Om\times\bbr\times\bbr^d$.
\elem

For the remaining analysis, by writing $f$ as the difference of its positive and negative part, which still satisfy Assumption~\ref{AssLLN}1, we may assume that $f$ is nonnegative. Then both $V^n_f(Y_x,t)$ and $V_f(Y_x,t)$ are increasing processes in $t$, so it suffices to prove $\bbe[|V^n_f(Y_x,t) - V_f(Y_x,t)|]\to0$ for every $t\geq0$. 
By definition, we have
\beq\label{increment} \un \Delta^n_i Y_x = \iint \un\Delta^n_i G_{x,y}(s) \si(s,y)\,W(\dd s,\dd y). \eeq
As we shall see, asymptotically as $\Delta_n\to0$, only the portion of the integral where $s$ is close to $i\del$ contributes to the size of $\un\Delta^n_i Y_x$. More precisely, we have the following result: 
\blem\label{lemma-eps} 
For $\eps>0$, if we define
\beq\label{alpha} \alpha^{n,i,\eps}_x = \iint  \un\Delta^n_i G_{x,y}(s)\si(s,y)\bone_{s>i\del-\eps}\,W(\dd s,\dd y), \eeq
then
$V^n_f(Y_x,t)-\del\sumt f\lp \frac{\alpha^{n,i,\eps}_x}{\tau_n}\rp \limL 0$ as $n\to\infty$.
\elem

As a next step, we discretize the volatility process in \eqref{alpha} along the points $i\del-\eps$. By the following lemma, this only introduces an asymptotically negligible error.
\blem\label{LLN-disc} If 
$$\wh\alpha^{n,i,\eps}_x = \iint \un \Delta^n_i G_{x,y}(s)\si(i\del-\eps,y)\bone_{s>i\del-\eps}\,W(\dd s,\dd y),$$
then
\[\lim_{\eps\to0}\limsup_{n\to\infty}\bbe\left[\supt\left|\del\sumt \left\{f\lp \frac{\alpha^{n,i,\eps}_x}{\tau_n}\rp-f\lp \frac{\wh\alpha^{n,i,\eps}_x}{\tau_n}\rp\right\}\right|\right] =0. \]
\elem

For small $\eps$, many of the integrals in \eqref{alpha} are taken over disjoint intervals. By exploiting this kind of conditional independence, we shall prove the following result:
\blem\label{centering} If $\si$ is bounded, then we have for every $t\geq0$,
\beq\label{expr-2} \sup_{n\in\N}\bbe\left[\left|\del\sumt \left\{f\lp \frac{\wh\alpha^{n,i,\eps}_x}{\tau_n}\rp-\bbe\lb f\lp \frac{\wh\alpha^{n,i,\eps}_x}{\tau_n} \rp\,\bigg|\, \calf_{i\del-\eps} \rb \right\}\right|^2\right]\to0  \eeq
as $\eps\to0$.
\elem

Finally, we show that the conditional expectations in \eqref{expr-2} converge to the correct limit.
\blem\label{cond-exp-conv} If $\si$ is bounded, then we have for every $t\geq0$,
\beq\label{conv} \begin{split} & \limsup_{n\to\infty}\bbe \lb \left|\del\sumt \bbe\lb f\lp \frac{\wh\alpha^{n,i,\eps}_x}{\tau_n} \rp\,\bigg|\, \calf_{i\del-\eps} \rb - \int_0^t \mu_f(\un\si^2(r,x))\,\dd r\right|\rb\\
&\quad \to 0\qquad\text{as}\qquad\eps\to0.\end{split}  \eeq
\elem

\subsection{Overview of the proof of Theorem~\ref{CLT}}\label{overview}

At the heart of our proof, we use a martingale central limit theorem for triangular arrays (see Theorem~2.2.15 in \citeA{Jacod12})  to obtain the stable convergence in law in \eqref{CLT-statement}. However, since the solution process \eqref{Yx} to \eqref{SHE} at fixed spatial points lacks the semimartingale property, many approximations are needed---before and after---to turn the left-hand side of \eqref{CLT-statement} into a term with a martingale structure. 

\subsubsection*{Step 1: Martingalization}
 We want to truncate the increments $\un\Delta^n_i Y_x$ in a similar way as in \eqref{alpha} to make sure that a large portion of the truncated increments are stochastic integrals over disjoint intervals (and hence have some sort of conditional independence) for different values of $i$. However, if we take a fixed level of $\eps$ as in \eqref{alpha}, the number of overlapping increments will still be of order $\eps/\del$. Consequently, the total approximation error for $V^n_f(Y_x,t)$ will be of order $\del (\eps/\del) = \eps$, which is not sufficient due to the $1/\sqrt{\del}$ prefactor. Of course, the best truncation that one can hope for is to only keep the integral in \eqref{alpha} on the set $\{s>(i-1)\del\}$, so that a given increment will only overlap with a finite number of other increments. But since $|\Pi^n_{0,0}|((\del,\infty)\times\R^d\times\R^d)\not\to 0$, this approximation is simply not valid (even without dividing by $\sqrt{\del}$). The idea is therefore to consider a truncation in between, that is, to take the integral in \eqref{alpha} only on the set $\{s>i\del-\la_n\del\}$ where $\la_n$ is a sequence increasing to $\infty$ with $\la_n\del\to0$. As it turns out,  the best (i.e., smallest) choice for $\la_n$ is achieved when we carry out the truncation iteratively.
 Hence, in a first step, we consider truncations of the form
 \beq\label{Ykn} \ga^{n,i,0}_x= \iint \un \Delta^n_i G_{x,y}(s)\si(s,y)\bone_{s>i\del-\la^0_n\del}\,W(\dd s,\dd y), \eeq
for which we have the following result:
\blem\label{int-trunc} If $\la_n^{0} = [\Delta_n^{-a_0}]$ for some $a_0>\frac1 \nu$, where $\nu=1+\frac\al 2$ as in Lemma~\ref{nu}, then
\beq\label{int-trunc-res}  \Delh \sumt \left\{ f\bigg( \frac{\un\Delta^n_i Y_x}{\tau_n}\bigg) - f\bigg( \frac{\ga^{n,i,0}_x}{\tau_n}\bigg) \right\} \limL 0. \eeq
\elem
Since two truncated increments $\ga^{n,i,0}_x$ and $\ga^{n,j,0}_x$ are defined on disjoint intervals as soon as $|i-j|>\la_n^0+L-1$, we can employ martingale techniques to improve (i.e., decrease) the order of $\la_n^0$. As mentioned above, we use an iterative truncation procedure and consider numbers $a_1>\cdots> a_R>\frac{1}{2\nu}$ satisfying $a_r>\frac{a_{r-1}}{\nu}$ for all $r=1,\ldots,R$. We define $\la_n^{r}=[\del^{-a_r}]$ and $\ga^{n,i,r}_x$ as in \eqref{Ykn} but with $\la_n^0$ replaced by $\la_n^r$. Furthermore, with the notation
\beq\label{not-filt} \calf^n_i = \calf_{i\del}, \eeq
we introduce the variables 
\beq\label{gamma} \begin{split} \delta^{n,r}_i &=\Delh \lp f\bigg( \frac{\ga^{n,i,r-1}_x}{\tau_n}\bigg) - f\bigg(\frac{\ga^{n,i,r}_x}{\tau_n}\bigg)\rp,\\ \ov\delta^{n,r}_i &= \delta^{n,r}_i - \bbe[\delta^{n,r}_i \,|\, \calf^n_{i-\la_n^{r-1}}]\end{split} \eeq
for $r=1,\ldots, R$. 
As $n\to\infty$, we can now shrink the domain of integration to the set $\{s>i\del-\la_n^R\del\}$ as the following two lemmas show:
\blem\label{knknprime}
For every $r=1,\ldots, R$,  $\sumt \ov\delta^{n,r}_i\limL 0$.
\elem

\blem\label{condex-trunc}
For every $r=1,\ldots,R$,  $\sumt \bbe[\delta^{n,r}_i\,|\, \calf^n_{i-\la_n^{r-1}}] \limL 0$.
\elem

In what follows, we define $a=a_R$, $\la_n = \la_n^R$, and $\ga^{n,i}_x = \ga^{n,i,R}_x$. Since $\nu>1$, after possibly increasing $R$, we may assume that $a$ is larger but arbitrarily close to $\frac1{2\nu}$. Although the iterative truncation procedure above has greatly reduced the number of overlapping increments (one increment now overlaps with roughly $\la_n$ instead of $\del^{-1}$ increments), this number $\la_n$ is still increasing in $n$, and hence, the increments are still far from having a martingale structure.
A classical block splitting technique, similar to \citeA{Jacod10} (see also Chapter~12.2.4 in \citeA{Jacod12}), will now help us to finally obtain a martingale array. To this end, 
define 
\[ V^n(t)=\sumt \psi^n_i, \qquad \psi^n_i =\delh\lp f\left(\frac{\ga^{n,i}_x}{\tau_n}\right)-\bbe\lb f\lp \frac{\ga^{n,i}_x }{\tau_n}\rp \,\Big|\, \calf^n_{i-\la_n}\rb\rp. \]
We now arrange the summands $\psi^n_i$ into blocks of length $m\la_n$ (where $m\in\N$), leaving out $\la_n+L-1$ terms between two consecutive blocks. More precisely, we decompose $V^n(t)$ into $V^n(t) =  V^{n,m,1}(t)+  V^{n,m,2}(t)+ V^{n,m,3}(t)$ with
\begin{align*}
 V^{n,m,1}(t) &= \sum_{j=1}^{J^{n,m}(t)} V^{n,m}_j, \\
 V^{n,m,2}(t)&=\sum_{j=1}^{J^{n,m}(t)}\sum_{k=1}^{\la_n+L-1} \psi^n_{(j-1)((m+1)\la_n+L-1)+m\la_n+k},\\
  V^{n,m,3}(t)&=\sum_{i=J^{n,m}(t)((m+1)\la_n+L-1)+1}^{t^*_n} \psi^n_i.
\end{align*}	
Here,  
\[ V^{n,m}_j=\sum_{k=1}^{m\la_n} \psi^n_{(j-1)((m+1)\la_n+L-1)+k},\qquad j=1,\ldots,J^{n,m}(t),   \]
are the blocks that we have up to time $t$. There are $J^{n,m}(t)=\lb\frac{t^*_n}{(m+1)\la_n+L-1}\rb$ complete blocks and possibly a boundary term $V^{n,m,3}(t)$, while $V^{n,m,2}(t)$ contains all summands that have been left out between blocks.
\blem\label{approx} If $a>\frac 1 {2\nu}$ is sufficiently small, we have for $i=2,3$,
\[ \lim_{m\to\infty}\limsup_{n\to\infty} \bbe\lb \supt | V^{n,m,i}(t)|\rb = 0. \]
\elem
For different values of $j$, the terms $V^{n,m}_j$ in the definition of $V^{n,m,1}$ comprise stochastic integrals over disjoint time domains. But the volatility process is evaluated continuously in time, so in order to finally obtain a martingale structure, we fix $\si$ at the beginning of each block $V^{n,m}_j$. To this end, we 
define
\beq\label{hatpsi}
\begin{split}
	\wh\xi^n_{i,k} &= \iint \frac{\un\Delta^n_i G_{x,y}(s)}{\tau_n} \si( (i-\la_n-k)\del,y)\bone_{s>(i-\la_n)\del} \,W(\dd s,\dd y),\\
	\wh\psi^n_{i,k} &=\Delh  \Big( f(\wh \xi^n_{i,k})-\bbe[f(\wh \xi^n_{i,k})\,|\, \calf^n_{i-\la_n-k}]\Big),\\
\end{split}
\eeq
and
\beq\label{hatV}\begin{split}
\wh V^{n,m,1}(t)&= \sum_{j=1}^{J^{n,m}(t)} \wh V^{n,m}_j,\\
 \wh V^{n,m}_j&= \sum_{k=1}^{m\la_n} \wh \psi^n_{(j-1)((m+1)\la_n+L-1)+k,k},\qquad j=1,\ldots,J^{n,m}(t).\end{split}
\eeq

\blem\label{Vhat} For every $m\in\N$,  we have $V^{n,m,1}(t)-\wh V^{n,m,1}(t) \limL 0$ as $n\to\infty$ if $a>\frac 1 {2\nu}$ is small.
\elem

\subsubsection*{Step 2: The martingale central limit theorem}
From \eqref{hatpsi} and \eqref{hatV}, it is easy to see that for every $m\in\N$, the variables $((\wh V^{n,m}_j)_{j=1,\ldots,J^{n,m}(t)}\colon n\in\N)$ form a \emph{triangular array}, in the sense of Chapter~2.2.4 in \citeA{Jacod12}, with respect to the filtrations $((\calf^n_{j((m+1)\la_n+L-1)-\la_n})_{j=0,\ldots,J^{n,m}(t)} \colon n\in\N)$. We use Proposition~2.2.4 and Theorem~2.2.15 in \citeA{Jacod12} to establish the asymptotic distribution of $\wh V^{n,m,1}(t)$.
\bprop\label{CLT-core} For small $a>\frac 1 {2\nu}$,
we have $\wh V^{n,m,1}\limst \calz^{(m)}$, where $\calz^{(m)}$ is a process characterized by the same properties as $\calz$ in Theorem~\ref{CLT} but with $\calc(t)$ replaced by $\calc^m(t) = \frac{m}{m+1} \calc(t)$.
\eprop

We need two approximation results, whose proofs are given after that of Proposition~\ref{CLT-core}.
\blem\label{si-bounded}
In order to prove Proposition~\ref{CLT-core}, we may assume that $\si$ is bounded.
\elem
\blem\label{f-poly}
In order to prove Proposition~\ref{CLT-core}, we may assume that $f$ is an even polynomial.
\elem

 Because $\bbe[\wh V^{n,m}_j\,|\, \calf^n_{(j-1)((m+1)\la_n+L-1)-\la_n}]=0$ by construction, it remains to prove the following properties. In point (3) below, we consider a sequence  $((W^{\iota}(t))_{t\in\R}\colon \iota\in\N)$ of independent two-sided standard Brownian motions such that the stochastic integral of any $H\in\call(W)$ against $W$ can be expressed as an $L^2$-series of stochastic integrals against $W^{\iota}$. The existence of such a sequence follows as in Section~2.3 of \citeA{Dalang11}; see their Proposition~2.6(b) in particular.
\benu
\item For all $t>0$ and $m_1,m_2=1,\ldots,M$, we have as $n\to\infty$,
\beq\label{Cmt} \begin{split} \calc^{n,m}_{m_1 m_2}(t)&=\sum_{j=1}^{J^{n,m}(t)} \bbe\lb (\wh V^{n,m}_j)_{m_1}(\wh V^{n,m}_j)_{m_2}  \,|\, \calf^n_{(j-1)((m+1)\la_n+L-1)-\la_n}\rb\\
	&\stackrel{\bbp}{\longrightarrow} \calc^m_{m_1m_2}(t). \raisetag{-3\baselineskip} \end{split} \eeq
\item There is $q>2$ such that for all $t>0$ and $m\in\N$, we have as $n\to\infty$,
\beq\label{conv-r} \sum_{j=1}^{J^{n,m}(t)} \bbe\left[|\wh V^{n,m}_j|^q \,|\, \calf^n_{(j-1)((m+1)\la_n+L-1)-\la_n}\right] \stackrel{\bbp}{\longrightarrow}  0. \eeq
\item Let $(M(t))_{t\geq0}$ be either the restriction to $[0,\infty)$ of $W^\iota$ for some $\iota\in\N$ or a bounded $(\calf_t)_{t\geq0}$-martingale that is orthogonal (in the martingale sense) to $W^\iota$ for all $\iota\in\N$. Then, given $m\in\N$ and $t>0$, we have
\beq\label{conv-mart} \begin{split} \qquad&\sum_{j=1}^{J^{n,m}(t)} \bbe\lb \wh V^{n,m}_j\left(M({\tau^{n}_j})-M({\tau^{n}_{j-1}})\right) \,|\, \calf^n_{(j-1)((m+1)\la_n+L-1)-\la_n}\rb\\ &\quad\stackrel{\bbp}{\longrightarrow} 0 \as,\end{split}\raisetag{-2\baselineskip}\eeq
where $\tau^{n}_j = \inf\{t\geq0\colon J^{n,m}(t)\geq j\} = (j((m+1)\la_n+L-1)+L-1)\Del$.
\eenu
Let us remark that in contrast to Theorem~2.2.15 in \citeA{Jacod12}, we fix a countable, and not just a finite number of Brownian motions and then consider martingales orthogonal to all these Brownian motions. The proof of the mentioned theorem, which is based on \citeA{Jacod97}, extends to this situation with no change.

\subsubsection*{Step 3: Computing the conditional expectation}
Since $\calc^m \limL \calc$, the results up to now and Proposition~2.2.4 in \citeA{Jacod12} imply that in the limit $n\to\infty$,
\beq\label{intermediate} \del^{-\frac12} \lp V^n_f(Y_x,t)-\del\sumt \bbe\lb f\lp \frac{\ga^{n,i}_x}{\tau_n}\rp \,\Big|\, \calf^n_{i-\la_n}\rb \rp \limst \calz. \eeq
The conditional expectation cannot be computed explicitly as it involves the volatility process sampled continuously. The purpose of the next lemma is therefore to discretize $\si$ at a fixed number of intermediate points (in the same spirit as Lemma~3 in \citeA{BN11}), where we use the following notation:
	\beq\label{abbr} t^{n,i}_{k} = i\Delta_n-k\Delta_n,\qquad I^{n,i}_{k,l} = (t^{n,i}_{k},t^{n,i}_{l}),\qquad \bone^{n,i}_{k,l}(s) = \bone_{I^{n,i}_{k,l}}(s).\eeq
\blem\label{discretization} If $a>\frac{1}{2\nu}$ is small enough, there exist numbers $a>a^{(1)}>\cdots>a^{(Q-1)}$ such that 
\beq\label{term-2} \Delh \sumt \left\{ \bbe\lb f\lp \frac{\ga^{n,i}_x}{\tau_n}\rp \,\Big|\, \calf^n_{i-\la_n}\rb-\bbe\lb f( \theta^n_i ) \,|\, \calf^n_{i-\la_n}\rb\right\} \limL 0, \eeq
where
\beq\label{thetani} \theta^n_i=\iint \frac{\un \Delta^n_i G_{x,y}(s)}{\tau_n} \sum_{q=1}^Q \si\lp t^{n,i}_{\la_n^{(q-1)}},y\rp\bone^{n,i}_{\la_n^{(q-1)},\la_n^{(q)}}(s) \,W(\dd s,\dd y),\eeq
and $\la_n^{(0)}=\la_n=[\del^{-a}]$, $\la_n^{(q)}=[\del^{-a^{(q)}}]$ for $q=1,\ldots,Q-1$, and $\la_n^{(Q)}=0$.
\elem
Notice the difference between $\la_n^r$ and $a_r$ as defined before Lemma~\ref{knknprime}  and $\la_n^{(q)}$ and $a^{(q)}$ as defined in the previous lemma. In fact, we have 
\begin{align*}
0=\la_n^{(Q)}<\la^{(Q-1)}_n<\cdots<\la_n^{(0)}&=\la_n=\la_n^R<\cdots< \la_n^0,\\
 a^{(Q-1)}<\cdots<a^{(0)}&=a=a_R<\cdots<a_0.
 \end{align*}
 
After suitable approximations, the conditional expectation of $f(\theta^n_i)$ can now be evaluated along these intermediate time points. To this end, define $m^{n,i}_r \in \R^{N\times L}$ and $v^{n,i}_r\in(\R^{N\times L})^2$ by
\beq\label{mv}\begin{split}
(m^{n,i}_r)_{jk} &=\iint  \frac{\Delta^n_{i+k-1} G_{x_j-y}(s)}{\tau_n} \sum_{q=1}^{r} \si( t^{n,i}_{\la_n^{(q-1)}},y)\bone^{n,i}_{\la_n^{(q-1)},\la_n^{(q)}}(s)\,W(\dd s,\dd y),\\
(v^{n,i}_r)_{jk,j'k'} &= \iiint   \frac{\Delta^n_{i+k-1} G_{x_j-y}(s)\Delta^n_{i+k'-1} G_{x_{j'}-z}(s)}{\tau^2_n}\\
&\quad\times \sum_{q=r+1}^{Q} \si( t^{n,i}_{\la_n^{(q-1)}},y)\si( t^{n,i}_{\la_n^{(q-1)}},z)\bone^{n,i}_{\la_n^{(q-1)},\la_n^{(q)}}(s)\,\dd s\,\La(\dd y,\dd z)
\end{split}\eeq
for $r=0,\ldots,Q$ and $j,j'=1,\ldots,N$ and $k,k'=1,\ldots,L$. In particular, $m^{n,i}_0 = 0$, $m^{n,i}_Q = \theta^n_i$, and $v^{n,i}_Q=0$. Then, recalling the definition of $\un\mu_f$ after \eqref{muf}, the following results hold if  $a$ is  sufficiently close to $\frac{1}{2\nu}$:
\blem\label{iter-condexp} 
$\Delh \sumt \left\{ \bbe[f(\theta^n_i)\,|\, \calf^n_{i-\la_n}] -  \un\mu_f(\bbe[v^{n,i}_0\,|\, \calf^n_{i-\la_n}])\right\} \limL 0.$
\elem
\blem\label{remove-condexp} 
$\Delh\sumt \left\{\un\mu_f(\bbe[v^{n,i}_0 \,|\, \calf^n_{i-\la_n}]) - \un\mu_f(v^{n,i}_0)\right\} \limL 0.$
\elem

\subsubsection*{Step 4: Approximation of the Lebesgue integral}
 We complete the proof of Theorem~\ref{CLT} by showing in two steps that $\del\sum_{i=1}^{[t/\del]} \un\mu_f(v^{n,i}_0)$ approximates the integral $V_f(Y_x,t)=\int_0^t \mu_f(\un\si^2(s,x))\,\dd s$ at a rate faster than $\sqrt{\del}$. 
\blem\label{remove-s} 
$ \delh\sumt \left\{\un\mu_f(v^{n,i}_0) -\mu_f(\un\si^2((i-1)\del,x))\right\}\limL 0.$
\elem
\blem\label{lebesgue} 
$\Delta_n^{-\frac12}\lp \del\sumt \mu_f(\un\si^2((i-1)\del,x)) -V_f(Y_x,t)\rp\limL 0.$
\elem
As expected from the semimartingale literature (cf.\ \citeA{Jacod12}), the two previous lemmas only hold true under strong regularity assumptions on $\si$. In fact, Assumption~\ref{AssCLT}2 and the spatial differentiability assumptions on $\si$ in \ref{AssCLT}3 are only needed for this step.

\section*{Acknowledgments}
I would like to thank Robert Dalang, Jean Jacod, and Mark Podolskij for their valuable advice. Moreover, I am grateful to two anonymous referees for their careful reading of the manuscript, which led to many improvements of the presentation.

\appendix
\begin{supplement}\label{suppA}
	\stitle{Supplement to ``High-frequency analysis of parabolic stochastic PDEs.''}
	\sdescription{This paper is accompanied by supplementary material in \citeA{Chong18a}. Section~\ref{SectA} in \citeA{Chong18a} gives some auxiliary results needed for the proofs in this paper. In Section~\ref{heat-est}, some important estimates related to the heat kernel are derived. Sections~\ref{details-LLN} and \ref{details} provide the details for the proof of Theorems~\ref{LLN} and \ref{CLT}, respectively. And finally, Section~\ref{proof-rest} contains the proofs for Sections~\ref{multipower} and \ref{est}.}
\end{supplement}

\bibliographystyleA{plainnat}
\bibliographyA{bib-VolaEstimation}

\printaddresses

\newpage
\renewcommand*\footnoterule{}
\setcounter{page}{1}
\setattribute{journal}{name}{}
\allowdisplaybreaks[3]
\begin{frontmatter}
	
	\title{Supplement to ``High-frequency analysis of parabolic stochastic PDEs''}
	\runtitle{High-frequency analysis of parabolic stochastic PDEs}
	
	\begin{aug}
		\author{\fnms{Carsten} \snm{Chong}\ead[label=e2]{carsten.chong@epfl.ch}}
		
		\runauthor{C. Chong}
		
		\affiliation{École Polytechnique Fédérale de Lausanne}
		
		\address{Institut de mathématiques\\ École Polytechnique Fédérale de Lausanne\\ Station 8\\ CH-1015 Lausanne\\ \printead{e2}}
	\end{aug}

	\begin{abstract}
		In this supplement, we prove, in Section~\ref{SectA}, some auxiliary results needed for the proofs in the main text \citeB{Chong18}. In Section~\ref{heat-est}, some important estimates related to the heat kernel are derived. Sections~\ref{details-LLN} and \ref{details} provide the details for the proof of Theorems~\ref{LLN} and \ref{CLT} in \citeB{Chong18}, respectively. And finally, Section~\ref{proof-rest} contains the proofs for Sections~\ref{multipower} and \ref{est} in \citeB{Chong18}. 
	\end{abstract}

\end{frontmatter}

	We use the numberings and notations from the main text \citeB{Chong18}.
	
	\appendix

	\section{Auxiliary results}\label{SectA}
	
	\blem\label{takesform} Let $\si^{(1)} = \si -\si^{(0)}$ be the stochastic integral process in \eqref{si-semi}. Under Assumption~\ref{AssCLT}, the following statements hold:
	\benu
	\item For fixed $x\in\R^d$, we have the representation
	\beq\label{si-fixed-x}\begin{split} \si^{(1)}(t,x) &= \si^{(1)}(-1,x)+\int_{-1}^t \si^{(11)}_x(s)\,\dd s \\
		&\quad+ \iint_{-1}^t \si^{(12)}_x(s,y)\,W'(\dd s,\dd y), \qquad t>-1,\end{split} \eeq
	where
	\beq\label{si11-12} \begin{split} \si^{(11)}_x(s)&=\iint_{-\infty}^s \frac{\partial}{\partial t} K(s-r,x-y)\rho(r,y)\,W'(\dd r,\dd y),\\
		\si^{(12)}_x(s,y)&=K(0,x-y)\rho(s,y), \end{split}\eeq  
	and both integrals in \eqref{si-fixed-x} are well defined. 
	\item For any $t\in\R$, the map $x\mapsto \si^{(1)}(t,x)$ is twice differentiable, and \eqref{si-der} also holds for $\si^{(1)}$.
	\eenu\elem
	\bpr
	Since $K$ is differentiable in $t$, by the fundamental theorem of calculus and the stochastic Fubini theorem (see Theorem~2.6 in \citeB{Walsh86-2}), we can rewrite $\si^{(1)}(t,x)-\si^{(1)}(-1,x)$ as
	\begin{align*}
	&\iint_{-\infty}^{-1} (K(t-s,x-y)-K(-1-s,x-y))\rho(s,y)\,W'(\dd s,\dd y) \\
	&\qquad+\iint_{-1}^{t} K(t-s,x-y)\rho(s,y)\,W'(\dd s,\dd y)\\
	&\quad=\int_{-1}^t \iint_{-\infty}^{r} \frac{\partial}{\partial t} K(r-s,x-y)\rho(s,y)\,W'(\dd s,\dd y)\,\dd r\\
	&\quad\quad+\iint_{-1}^t K(0,x-y)\rho(s,y)\,W'(\dd s,\dd y), \end{align*}
	from which \eqref{si-fixed-x} follows. 
	With similar reasoning, one can show that the partial derivatives $\frac{\partial}{\partial x_j} \si^{(1)}$, $\frac{\partial^2}{\partial x_j\,\partial x_k} \si^{(1)}$, and $\frac{\partial^3}{\partial x_j\,\partial x_k\,\partial x_l} \si^{(1)}$ for $j,k,l=1,\ldots,d$ exist and are given by the same stochastic integral as in \eqref{si-semi} but with $K$ replaced by the corresponding partial derivative of $K$. Moreover, since the third derivatives exist, all first- and second-order partial derivatives of $\si^{(1)}$ are continuous, which implies that $\si^{(1)}$ is twice (totally) differentiable. Property \eqref{si-der} for $\si^{(1)}$ now follows from the integral representation of the partial derivatives together with the assumptions on $K$ and $\rho$.
	\epr
	
	\blem\label{mom-ex} Under Assumption~\ref{AssCLT}, also the whole volatility process $\si$ satisfies \eqref{mom-si} as well as \eqref{Holder} with $\ga=\frac12$. 
	\elem
	\bpr For $\si^{(0)}$, both statements hold by assumption. Thanks to the Burkholder--Davis--Gundy (BDG) inequality, the absolute moment of order $2p+\eps$ of the integral against $W'$ is bounded by a constant times
	\begin{align*}
	\sup_{(t,x)\in\R\times\R^d} \bbe[|\rho(t,x)|^{2p+\eps}] \left( \iiint_0^\infty |K(s,y)K(s,z)|\,\dd s\,\La(\dd y,\dd z)\right)^{\frac {2p+\eps} 2}.
	\end{align*}
	This is finite and independent of $t$ and $x$ since $\rho$ satisfies \eqref{mom-si} and $K\in \call(W')$. Hence, we have \eqref{mom-si} for $\si$. With a similar moment estimate, \eqref{Holder} follows from the representation \eqref{si-fixed-x}. 
	\epr

	\blem\label{cov-add}  Suppose that $X=(X_1,\ldots,X_d)'$ and $Y=(Y_1,\ldots,Y_d)'$ are jointly Gaussian and identically distributed random vectors with zero mean. If $f_1, f_2\colon \R^d\to\R$ are even and twice differentiable such that for some $p\geq0$, we have $|\frac{\partial^2}{\partial x_i \, \partial x_j} f_1(x)|+|\frac{\partial^2}{\partial x_i \, \partial x_j} f_2(x)|\lec 1+|x|^p$ for all $i,j=1,\ldots,d$, then
	\beq\label{toprove}\begin{split}
		|\cov(f_1(X),f_2(Y))|&\lec \lp \sum_{i,j=1}^d \frac{c_{ij}^2}{2} +\sum_{(i_1,j_1)\neq (i_2,j_2)} |c_{i_1j_1}c_{i_2 j_2}|\rp\\
		&\quad\times \left(1+\sum_{i=1}^d |v_i|^{\frac p2}\right),\end{split}
	\eeq
	where $v_i=\bbe[X_i^2]=\bbe[Y_i^2]$ and $c_{ij} = \bbe[X_iY_j]$.
	\elem
	\bpr By an approximation argument, we may assume that the covariance matrix of $X$ (which coincides with that of $Y$) is symmetric positive definite. Then it is well known that there exists an orthogonal matrix $A\in\R^{d\times d}$ such that $AX$ and $AY$ are Gaussian vectors with independent nondegenerate entries, respectively. Writing $f_{1|2}(x)=f^A_{1|2}(Ax)$, where $f^A_{1|2}(x)=f_{1|2}(A^{-1}x)$ are again even functions such that $|\frac{\partial^2}{\partial x_i \, \partial x_j} f^A_{1|2}(x)|\lec 1+|x|^p$ with the same constant as for $f_{1|2}$, we may restrict ourselves to the case where $X$ and $Y$ consist of independent entries, respectively, with $v_i>0$ for all $i=1,\ldots,d$. 
	
	Now consider  the Hermite polynomials $H_n(x)=\frac{(-1)^n}{n!}\ee^{x^2/2} \frac{\dd^n}{\dd x^n} \ee^{-x^2/2}$ for $n\in\N_0$ (with $H_0(x)=1$) and let $H_n(v,x) = \prod_{i=1}^d v_i^{ n_i/2}H_{n_i}(\frac{x_i}{\sqrt{v_i}})$ for $n\in\N^d_0$ and $x\in\R^d$ be  multivariate generalizations thereof. 
	Since $\bbe[f^2_{1|2}(X)]<\infty$ and $f_{1|2}$ is even, we have (cf.\ Proposition~1.1.1 and Lemma~1.1.1 in \citeB{Nualart06-2}) $$f_{1|2}(x) = \sum_{n \in \N^d_0} a_{n}(v,f_{1|2})H_{n}(v,x) =  \sum_{|n|\, \text{even}} a_{n}(v,f_{1|2})H_{n}(v,x)$$
	with $a_{n}(v,f)=\frac{n!}{v^{n}}\bbe[f(X)H_n(v,X)]$  
	and multi-index notation $|n|=n_1+\cdots+n_d$, $n!=n_1!\cdots n_d!$, and $v^n = v_1^{n_1}\cdots v_d^{n_d}$. If all entries of $v$ are equal to $1$, we simply write $H_n(x)$ and $a_n(f)$. 
	
	Next, with a similar reasoning as in Lemma~1.1.1 in  \citeB{Nualart06-2} and using Equation~(35) in \citeB{Slepian72}, we have
	\beq\label{cross} \bbe[H_n(v,X)H_m(v,Y)] = v^{\frac{n+m}{2}} \sum_{\begin{smallmatrix} \mu\colon \frr(\mu)=n,\\ \frc(\mu)=m \end{smallmatrix}} \frac{c^\mu}{v^{\frac{\frr(\mu)+\frc(\mu)}{2}} \mu!} = \sum_{\begin{smallmatrix} \mu\colon \frr(\mu)=n,\\ \frc(\mu)=m \end{smallmatrix}} \frac{c^\mu}{\mu!} \eeq
	if $|n|=|m|$, and $\bbe[H_n(v,X)H_m(v,Y)] =0$ otherwise.
	In the last line, the sum is taken over all $(d\times d)$-matrices $\mu$ such that $\frr(\mu)_i = \sum_{j=1}^d \mu_{ij} = n_i$ and $\frc(\mu)_j = \sum_{i=1}^d \mu_{ij} = m_j$ for all $i,j=1,\ldots,d$. Moreover, $\mu!=\prod_{i,j=1}^d \mu_{ij}!$ and $c^{\mu}= \prod_{i,j=1}^d c_{ij}^{\mu_{ij}}$.
	
	As a consequence of \eqref{cross}, we obtain
	\begin{align*}
	\cov(f_1(X),f_2(Y)) = \sum_{\begin{smallmatrix} |n|=|m|\geq 2\\ \text{and even}\end{smallmatrix}} a_n(v,f_1)a_m(v,f_2) \sum_{\begin{smallmatrix} \mu\colon \frr(\mu)=n,\\ \frc(\mu)=m \end{smallmatrix}} \frac{c^\mu}{\mu!}.
	\end{align*}
	Expanding this expression via Taylor's formula at $c=0$, and observing that the zeroth- and the first-order terms vanish because $\mu$ contains at least one entry $\geq2$ or two entries $\geq1$, we derive
	\begin{align*}
	&\cov(f_1(X),f_2(Y))\\
	&\quad = \sum_{i,j=1}^d \frac{c_{ij}^2}{2} \sum_{\begin{smallmatrix} |n|=|m|\geq 2\\ \text{and even}\end{smallmatrix}} a_n(v,f_1)a_m(v,f_2) \sum_{\begin{smallmatrix} \mu\colon \frr(\mu)=n,\\ \frc(\mu)=m, \, \mu_{ij}\geq2 \end{smallmatrix}} \frac{(hc)^{\mu-2e_{ij}}}{(\mu-2e_{ij})!}\\
	&\qquad +\sum_{(i_1,j_1)\neq (i_2,j_2)} c_{i_1j_1}c_{i_2 j_2} \sum_{\begin{smallmatrix}  |n|=|m|\geq 2\\ \text{and even}\end{smallmatrix}} a_n(v,f_1)a_m(v,f_2)\\
	&\qquad\quad\times\sum_{\begin{smallmatrix} \mu\colon \frr(\mu)=n,\, \frc(\mu)=m, \\ \mu_{i_1j_1}\geq1,\, \mu_{i_2 j_2} \geq1 \end{smallmatrix}} \frac{(hc)^{\mu-e_{i_1j_1}-e_{i_2j_2}}}{(\mu-e_{i_1j_1}-e_{i_2j_2})!}.
	\end{align*}
	Here, $h$ lies between $0$ and $1$, and $e_{ij}$ is the matrix with entry $1$ at the position $(i,j)$ and $0$ elsewhere. Next, we make the substitution $\mu-2e_{ij} \mapsto \mu$ (resp., $\mu-e_{i_1j_1}-e_{i_2j_2}\mapsto \mu$) as well as $(n-2e_i,m-2e_j)\mapsto (n,m)$ [resp., $(n-e_{i_1}-e_{i_2},m-e_{j_1}-e_{j_2})\mapsto (n,m)$], which yields
	\begin{align*}
	&\cov(f_1(X),f_2(Y))\\
	&\quad = \sum_{i,j=1}^d \frac{c_{ij}^2}{2} \sum_{|n|=|m|\,\text{even}} a_{n+2e_i}(v,f_1)a_{m+2e_j}(v,f_2) \sum_{\begin{smallmatrix} \mu\colon \frr(\mu)=n,\\ \frc(\mu)=m \end{smallmatrix}} \frac{(hc)^{\mu}}{\mu!}\\
	&\qquad +\sum_{(i_1,j_1)\neq (i_2,j_2)} c_{i_1j_1}c_{i_2 j_2} \sum_{|n|=|m|\,\text{even}} a_{n+e_{i_1}+e_{i_2}}(v,f_1)a_{m+e_{j_1}+e_{j_2}}(v,f_2)\\
	&\qquad\quad\times \sum_{\begin{smallmatrix} \mu\colon \frr(\mu)=n,\\ \frc(\mu)=m \end{smallmatrix}} \frac{(hc)^{\mu}}{\mu!}.
	\end{align*}
	Using Equations~(1.2) and (1.3) in \citeB{Nualart06-2} and integrating by parts, one may verify by induction that
	\beq\label{a-calc}\begin{split}
		a_{n+m}(v,f) = 	
		a_{n}(v,\partial^m f)
	\end{split}\eeq
	for all $m\in\N^d_0$ and functions $f$ where $\partial^m f= \frac{\partial^{|m|}}{\partial x_1^{m_1}\cdots\partial x_d^{m_d}} f$ exists. As a result, if we recall \eqref{cross}, $\cov(f_1(X),f_1(Y))$ is given by
	\begin{align*}
	& \sum_{i,j=1}^d \frac{c_{ij}^2}{2} \bbe\lb \frac{\partial^2}{\partial x_i^2} f_1(\sqrt{h} X) \frac{\partial^2}{\partial x_j^2} f_2(\sqrt{h} Y)\rb \\
	&\quad+\sum_{(i_1,j_1)\neq (i_2,j_2)} c_{i_1j_1}c_{i_2 j_2} \bbe\lb \frac{\partial^2}{\partial x_{i_1}\,\partial x_{i_2}} f_1(\sqrt{h} X) \frac{\partial^2}{\partial x_{j_1}\,\partial x_{j_2}} f_2(\sqrt{h} Y)\rb, 
	\end{align*}
	so \eqref{toprove} follows from the Cauchy--Schwarz inequality and the growth properties of $\frac{\partial^2}{\partial x_i\,\partial x_j} f_{1|2}$.
	\epr
	
	\section{Heat kernel estimates}\label{heat-est}
	
	According to \citeB{Dalang99-2}, if we denote by
	\beq\label{FT} \calf \phi(\xi)=\int_{\R^d} \ee^{-2\pi \ii \xi\cdot x}\phi(x)\,\dd x,\qquad \xi\in\R^d, \eeq
	the \emph{Fourier transform} of a compactly supported smooth function $\phi\colon\R^d\to\R$, then there exists a measure $\mu$ called the \emph{spectral measure} of $W$ such that
	for all compactly supported smooth functions $\psi_1,\psi_2\colon \R^{d+1}\to\R$, we have
	\beq\label{isometry} \begin{split} &\iiint \psi_1(s,y)\psi_2(s,z)F(y-z)\,\dd y\,\dd z\,\dd s\\ &\quad = \iint \calf \psi_1(s,\cdot)(\xi)\ov{\calf\psi_2(s,\cdot)(\xi)}\,\mu(\dd \xi)\,\dd s. \end{split} \eeq
	With our choice of $F$ (including the case $d=1$ and $F=\delta_0$), it is well known (see Chapter~V, \S 1, Lemma~1(a) in \citeB{Stein70}) that the spectral measure is given by $\mu(\dd \xi) =|\xi|^{\al-d}\,\dd \xi$.
	
	\blem\label{Ga} Let $\al\in(0,2)\cap(0,d]$.
	\benu
	\item We have
	\beq\label{taun-asymp} \begin{split}\tau_n^2 &= \frac{\pi^{\frac d 2-\al}\Ga(\frac\al 2)}{(2\kappa)^{\frac \al 2}\Ga(\frac d 2)\la^{1-\frac \al 2}}\ga(1-\textstyle \frac\al 2 , \la\del)\\ &= \displaystyle\frac{\pi^{\frac d 2-\al}\Ga(\frac\al 2)}{(2\kappa)^{\frac \al 2}(1-\frac \al 2)\Ga(\frac d 2)} \del^{1-\frac \al 2} +o(\del^{1-\frac\al 2})\end{split}\eeq
	as $n\to\infty$,
	where $\ga(a,x)=\int_0^x \ee^{-u}u^{a-1}\,\dd u$ is the lower incomplete gamma function.
	\item Let $\Ga^n_r=\Pi^n_{r,0}([0,\infty)\times\R^d\times\R^d)$ for $r\in\N_0$,
	where $\Pi^n_{r,0}$ is given in \eqref{Pin}. Then 
	\beq\label{Ganr}\Ga^n_r =\frac{\ga(1-\textstyle\frac\al 2 ,\la(r+1)\del)-2\ga(1-\frac\al 2, \la r \del)+\ga(1-\frac\al 2, \la(r-1)\del)}{2\ga(1-\textstyle \frac\al 2, \la\del)},\eeq
	and recalling the definition of $\Ga_r$ in \eqref{Ga-formula}, we have as $n\to\infty$,
	\beq\label{Ganr-speed} |\Ga^n_r-\Ga_r| = O(\del).\eeq
	\eenu
	\elem
	\bpr
	Let $Y_0$ be given by \eqref{Yx} with $\si$ replaced by the constant $1$. We claim that for any $t\in\R$ and $\tau\geq0$ we have
	\beq\label{Rt} R(\tau)=\bbe[(Y_0(t+\tau,x)-Y_0(t,x))^2] = \frac{\pi^{\frac d 2-\al}\Ga(\frac\al 2)}{(2\kappa)^{\frac \al 2}\Ga(\frac d 2)\la^{1-\frac\al 2}} \ga(1-\textstyle\frac\al 2, \la \tau). \eeq
	By \eqref{isometry}, we have
	\begin{align*}
	R(\tau)&= \iint_{-\infty}^{t} (\calf G(t+\tau-s,\cdot)(\xi)-\calf G(t-s,\cdot)(\xi))^2\,\dd s\,\mu(\dd \xi)\\ &\quad+ \iint_t^{t+\tau} (\calf G(t+\tau-s,\cdot)(\xi))^2 \,\dd s\,\mu(\dd \xi)\\
	&=\iint_{0}^{\infty} (\calf G(s+\tau,\cdot)(\xi)-\calf G(s,\cdot)(\xi))^2\,\dd s\,\mu(\dd \xi)\\ &\quad + \iint_0^{\tau} (\calf G(s,\cdot)(\xi))^2 \,\dd s\,\mu(\dd \xi).
	\end{align*}
	Since \beq\label{fourierG} (2\pi\kappa t)^{-\frac d 2}\calf \ee^{-\frac{|\cdot|^2}{2\kappa t}}(\xi)=\ee^{-2\pi^2\kappa t |\xi|^2},\eeq
	we obtain by first integrating with respect to time,
	\beq\label{termR}\begin{split}
		R(\tau) &= \iint_0^\infty \lp\ee^{-(s+\tau)(\la+2\pi^2\kappa |\xi|^2)}-\ee^{-s(\la+2\pk)}\rp^2\,\dd s\,\mu(\dd \xi)\\ &\quad+\iint_0^\tau \ee^{-2s(\la+2\pk)}\,\dd s\,\mu(\dd \xi)\\
		&=\int_{\R^d} \frac{(1-\ee^{-(\la+2\pk)\tau})^2}{2\la+4\pk}\,\mu(\dd\xi)+\int_{\R^d} \frac{1-\ee^{-2(\la  + 2\pi^2 \kappa |\xi|^2)\tau}}{2\la+4\pi^2\kappa |\xi|^2}\,\mu(\dd\xi)\\
		&=\int_{\R^d} \frac{1-\ee^{-(\la+2\pk)\tau}}{\la+2\pk}\,\mu(\dd\xi).
	\end{split}\raisetag{-4\baselineskip}\eeq
	By changing to polar coordinates $z=|\xi|$, and substituting $\sqrt{\frac{2\pi^2\kappa}{\la}}z \mapsto z$, we derive
	\begin{align*}
	R(\tau) &=|\cals^{d-1}|\int_0^\infty \frac{1-\ee^{-(\la+2\pi^2\kappa z^2)\tau}}{\la+2\pi^2\kappa z^2}z^{\al-1}\,\dd z\\ &=\frac{|\cals^{d-1}|}{(2\pi^2\kappa)^{\frac \al 2}\la^{1-\frac\al 2}}\int_0^\infty \frac{1-\ee^{-\la(1+z^2)\tau}}{1+z^2}z^{\al-1}\,\dd z,
	\end{align*}
	where $|\cals^{d-1}|={2\pi^{d/2}}{\Ga(\frac d 2)^{-1}}$ is the surface measure of the $d$-dimensional unit ball. As $1-\ee^{-ax}=a\int_0^x \ee^{-au}\,\dd u$ and $\int_0^\infty \ee^{-(1+z^2)u}z^{\al-1}\,\dd z= \frac12 \ee^{-u}u^{-\al/2}\Ga(\frac\al 2)$, we further obtain
	\begin{align*}
	R(\tau) &=\frac{|\cals^{d-1}|}{(2\pi^2\kappa)^{\frac \al 2}\la^{1-\frac\al 2}}  \int_0^\infty \lp\int_0^{\la \tau} \ee^{-(1+z^2)u}\,\dd u \rp z^{\al-1}\,\dd z\\ &= \frac{\pi^{\frac d 2-\al}\Ga(\frac\al 2)}{(2\kappa)^{\frac \al 2}\Ga(\frac d 2)\la^{1-\frac\al 2}} \int_0^{\la\tau} \ee^{-u}u^{-\frac\al 2}\,\dd u,
	\end{align*}
	from which \eqref{Rt} follows. Because $L_0(x)=\ga(1-\frac \al 2,x)x^{-(1-\al/2)}\to (1-\frac\al2)^{-1}$ as $x\to0$ and $\frac{\partial^2}{\partial x^2} \ga(1-\frac\al2,x)=\ee^{-x}(-\frac\al2-x)x^{-1-\al/2}$, one can easily show that conditions (A1), (A2), and (A3) in \citeB{BN09} are all satisfied with their $\beta$ equal to $1-\frac\al 2$. All assertions, except for \eqref{Ganr-speed}, now follow from Equation~(A.2) and Lemma~1 in \citeB{BN09}. Finally, \eqref{Ganr-speed} holds because the function $x\mapsto \Ga_r(x)$, defined by \eqref{Ganr} with $\del$ replaced by $x$, is differentiable. This can be easily seen from the identities $\Ga_0(x)=1$ and
	\[ \Ga_r(x)=\frac{(r+1)^{1-\frac\al2}L_0((r+1)x)-2r^{1-\frac\al2}L_0(rx)+(r+1)^{1-\frac\al2}L_0((r-1)x)}{2L_0(x)}\]
	for $r\geq1$ 
	and the fact that $L_0$ is differentiable on $[0,\infty)$ (see Equation~6.5.4 in \citeB{Abramowitz72}) and $L_0(0)>0$. 
	\epr
	
	\blem\label{delta} Let $\al\in(0,2)\cap(0,d]$ and consider the measures $\Pi^n_{r,h}$ and $|\Pi^n_{r,h}|$ defined in \eqref{Pin}.
	\benu
	\item As $n\to\infty$, we have for all $r\in\N_0$ that $\Pi^n_{r,0} \limw \Ga_r \delta_0$, where $\Ga_r$ is defined in \eqref{Ga-formula} and $\limw$ denotes weak convergence of measures. Moreover, for $h\neq0$, we have
	\beq\label{speed} |\Pi^n_{r,h}|([0,\infty)\times\R^d\times\R^d) \lec \begin{cases} \del^{\frac\al2}|\log \del|&\text{if } \al\in(0,d\wedge2),\\ \del^{\frac32}|\log\del|&\text{if } d=1,~F=\delta_0. \end{cases} \eeq
	In particular, we have $|\Pi^n_{r,h}|\limw 0$ for all $h\neq 0$. 
	\item There exists a decreasing sequence $(\ov\Ga_r\colon r\in\N_0)$ such that $\sum_{r=0}^\infty \ov\Ga_r^2< \infty$, and for all $n\in\N$, $r\in\N_0$, and $h\in\R^d$,
	\beq\label{cross-bound-3} |\Pi^n_{r,h}|([0,\infty)\times\R^d\times\R^d) \leq \ov \Ga_r. \eeq
	\eenu
	\elem
	\bpr Let us first consider
	$h\neq0$. The cases $r=0$ and $r=1$ can be treated  similarly to  below, so we only consider $r\geq2$. Recalling that $\La(\dd y,\dd z)=c_\al|z-y|^{-\al}\,\dd y\,\dd z$ (in the case $F=\delta_0$, the following quantity is simply $0$), we obtain
	\beq\label{Pi-first}\begin{split}
		&|\Pi^n_{r,h}|([0,\infty)\times\{(y,z)\colon |z-y|>\textstyle\frac{|h|}{4}\})\\
		&\quad\lec\int_0^\infty  \int_{\R^d} \frac{|G_y(s)-G_y(s-\del)|}{\tau_n}\,\dd y\\
		&\quad\quad\times \int_{\R^d} \frac{|G_{z+h}(s+r\del)-G_{z+h}(s+(r-1)\del)|}{\tau_n}\,\dd z\,\dd s\\
		&\quad =  \int_0^{\del} \ee^{-\la s} \int_{\R^d} \frac{|G_{z+h}(s+r\del)-G_{z+h}(s+(r-1)\del)|}{\tau^2_n}\,\dd z \,\dd s \\
		&\quad\quad + \int_0^\infty  \int_{\R^d} \frac{|G_y(s+\del)-G_y(s)|}{\tau_n}\,\dd y\\
		&\quad\qquad\times \int_{\R^d} \frac{|G_{z+h}(s+(r+1)\del)-G_{z+h}(s+r\del)|}{\tau_n}\,\dd z\,\dd s.
	\end{split}\eeq
	Let $r_n(s)=\sqrt{2\kappa s (1+\frac{s}{\del})(\frac d 2 \log(1+\frac{\del}{s})+\la\del)}$ be the value of $|y|$ where $G_y(s+\del)=G_y(s)$. In the following calculation, if we substitute $y/{\sqrt{\kappa s}} \mapsto y$ for $G_y(s)$ and $y/{\sqrt{\kappa (s+\del)}} \mapsto y$ for $G_y(s+\del)$ in the last step, we derive
	\begin{align*}
	&\int_{\R^d} |G_y(s+\del)-G_y(s)|\,\dd y\\
	&\quad=  \int_{|y|\leq r_n(s)} (G_y(s)-G_y(s+\del))\,\dd y + \int_{|y|>r_n(s)} (G_y(s+\del)-G_y(s))\,\dd y\\
	&\quad = \frac{\ee^{-\la s}}{(2\pi)^{\frac d 2}} \Bigg( (1-\ee^{-\la\del})\int_{|y|\leq \frac{r_n(s)}{\sqrt{\kappa (s+\del)}}} {\ee^{-\frac{|y|^2}{2}}}\,\dd y\\
	&\quad\quad + (1+\ee^{-\la\del}) \int_{\frac{r_n(s)}{\sqrt{\kappa (s+\del)}}<|y|\leq \frac{r_n(s)}{\sqrt{\kappa s}}} \ee^{-\frac{|y|^2}{2}}\,\dd y\\
	&\quad \quad -(1-\ee^{-\la\del})\int_{|y|> \frac{r_n(s)}{\sqrt{\kappa s}}} {\ee^{-\frac{|y|^2}{2}}}\,\dd y \Bigg).
	\end{align*}
	Since $r_n(s)\lec \sqrt{s}$ for $\del \leq s \leq 1$ and $r_n(s)\lec s$ for $s>1$, if we change to polar coordinates in the second integral of the previous display, we obtain for $\del\leq s \leq 1$,
	\begin{align*}
	&\int_{\R^d} |G_y(s+\del)-G_y(s)|\,\dd y\\
	&\quad\lec \ee^{-\la s} \lp\del + \left(\frac{ r_n(s)}{\sqrt{s}}\right)^{d-1} r_n(s) \lp \frac{1}{\sqrt{s}} -\frac{1}{\sqrt{s+\del}}\rp\rp\\
	&\quad\lec\ee^{-\la s} \lp\del+\sqrt{\frac{s}{{\del}}}\lp\frac{1}{\sqrt{s/\del}} -\frac{1}{\sqrt{s/\del+1}} \rp\rp\lec \ee^{-\la s} \frac{\del} s.
	\end{align*}
	For $s>1$, the left-hand side is bounded by $\del \ee^{-\la s}s^{(d/2 -1)\vee0}$. It is also bounded by $2$, regardless of the value of $s$, because $G$ is a density up to an exponential factor. Hence, we obtain
	\beq\label{int-bound-2} \int_{\R^d} |G_y(s+\del)-G_y(s)|\,\dd y\lec \begin{cases} 2,&0 \leq s \leq \del,\\ \frac{\del}s \ee^{-\frac\la 2 s}, &s> \del. \end{cases}\eeq
	Inserting this into \eqref{Pi-first}, we derive
	\beq\label{intermediate-1}\begin{split}
		&|\Pi^n_{r,h}|([0,\infty)\times\{(y,z)\colon |z-y|>\textstyle\frac{|h|}{4}\})\\
		&\quad \lec \frac{1}{\tau_n^2}\int_0^{\del} \left(\ee^{-\frac\la 2 (s+(r-1)\del)} \frac{\del}{s+(r-1)\del}+  \ee^{-\frac\la 2 (s+r\del)}\frac{\del}{s+r\del}\right) \,\dd s  \\
		&\qquad + \frac{1}{\tau_n^2}\int_{\del}^\infty \ee^{-\frac \la 2 s}\frac{\del}{s}\ee^{-\frac\la 2 (s+r\del)}\frac{\del}{s+r\del}\,\dd s\\
		&\quad \lec \frac{1}{\tau_n^2(r-1)} \int_0^{\del} \ee^{-\frac\la 2 s}\,\dd s+\frac{\del\ee^{-\frac\la 2 r\del}}{\tau_n^2 r}\int_{\del}^\infty \frac{\ee^{-\la s}}{s}\,\dd s \lec \frac{\del|\log \del|}{\tau_n^2(r-1)},
	\end{split}\raisetag{-4\baselineskip}\eeq
	where  we used $\frac{\del}{s+r\del}\leq \frac 1 r$ in the penultimate step and
	$\del\ee^{- \la  r\del/2} \lec \frac 1 r$ and $\int_{\del}^\infty \frac{\ee^{-\la s}}{s}\,\dd s \lec |\log \del|$ in the final step.

	Next, we consider $|\Pi^n_{r,h}|$ on the set $[0,\infty)\times\{(y,z)\colon |y|\leq \textstyle\frac{|h|}{4},~ |z-y|\leq\textstyle\frac{|h|}{4}\}$. In particular, we have $|z| = |z-y+y|\leq \frac{|h|}{2}$, and hence, $|z+h|\geq |h|-|z|\geq \frac{|h|}{2}$. An elementary calculation shows that $|\frac{\partial}{\partial t} G_{z+h}(t)|\lec 1\wedge t^{-d/ 2 -1}$, uniformly for $|z+h|\geq \frac{|h|}{2}$, so we have $|G_{z+h}(s+r\del)-G_{z+h}(s+(r-1)\del)|\lec (1\wedge ((r-1)\del)^{- d /2 -1})\del$ on the set considered above by the mean value theorem. Thus, with obvious modifications if $F=\delta_0$, we deduce from \eqref{int-bound-2},
	\beq\label{intermediate-2}\begin{split}
		&|\Pi^n_{r,h}|([0,\infty)\times\{(y,z)\colon |y|\leq \textstyle\frac{|h|}{4},~ |z-y|\leq\textstyle\frac{|h|}{4}\}) \\
		&\quad\lec \displaystyle \frac{(1\wedge((r-1)\del)^{-\frac d 2-1})\del}{\tau_n^2}\\
		&\qquad\times\int_0^\infty \int_{\R^d} |G_y(s)-G_y(s-\del)| \int_{|z-y|\leq \frac{|h|}{4}} |z-y|^{-\al}\,\dd z \,\dd y\,\dd s\\
		&\quad\lec \displaystyle \frac{(1\wedge((r-1)\del)^{-\frac d 2-1})\del}{\tau_n^2}\left(\del + \del \int_{\del}^\infty \frac{\ee^{-\frac \la 2 s}}{s}  \,\dd s\rp\\
		&\quad\lec \frac{\del|\log \del|}{(r-1)\tau_n^2},
	\end{split}\eeq
	where for the last inequality, we used the estimate $r\del(1\wedge (r\del)^{- d /2-1})\leq 1$.
	
	Furthermore, since the heat kernel is smooth outside the origin, we have $G_y(s)\lec 1$ and $|G_y(s+\del)-G_y(s)|\lec \del$ uniformly for all $|y|>\frac{|h|}{4}$. Thus, again with obvious changes if $d=1$ and $F=\delta_0$,
	\begin{align*}
	&|\Pi^n_{r,h}|([0,\infty)\times\{(y,z)\colon |y|> \textstyle\frac{|h|}{4},~ |z-y|\leq\textstyle\frac{|h|}{4}\})\\
	&\quad\lec\frac{1}{\tau_n^2}\int_0^{\del}\int_{\R^d} \lp\int_{|z-y|\leq\frac{|h|}{4}} |z-y|^{-\al}\,\dd y \rp\\
	&\qquad\quad\times|G_{z+h}(s+r\del)-G_{z+h}(s+(r-1)\del)|\,\dd z\,\dd s\\
	&\qquad+\frac{\del}{\tau_n^2}\int_0^{\infty} \int_{\R^d} \lp \int_{|z-y|\leq \frac{|h|}{4}} |z-y|^{-\al} \,\dd y\rp\\
	&\qquad\quad\times |G_{z+h}(s+(r+1)\del)-G_{z+h}(s+r\del)|\,\dd s\,\dd z.
	\end{align*}
	With similar techniques as in \eqref{intermediate-1}, this can be shown to be bounded by a constant times $\frac{\del}{(r-1)\tau_n^2}$. Together with \eqref{intermediate-1} and \eqref{intermediate-2}, this shows that
	\beq\label{hnot0} |\Pi^n_{r,h}|([0,\infty)\times\R^d\times\R^d) \lec \frac{\del|\log \del|}{(r+1)\tau_n^2}, \eeq
	which is actually true for all $r\in\N_0$ and, by \eqref{taun-asymp}, shows all assertions of the lemma for $h\neq0$, except for \eqref{speed} when $F=\delta_0$. For this, we decompose $|\Pi^n_{r,h}|([0,\infty)\times\R^d\times\R^d)=A_1+A_2$, where
	\begin{align*}
	A_1&= \iint \frac{|G_y(s)-G_y(s-\del)|}{\tau_n}\\
	&\quad\times\frac{|G_{y+h}(s+r\del)-G_{y+h}(s+(r-1)\del)|}{\tau_n}\bone_{|y|\leq \frac{h}{2}}\,\dd s\,\dd y,\\
	A_2&=\iint \frac{|G_y(s)-G_y(s-\del)|}{\tau_n}\\
	&\quad\times\frac{|G_{y+h}(s+r\del)-G_{y+h}(s+(r-1)\del)|}{\tau_n}\bone_{|y|> \frac{h}{2}}\,\dd s\,\dd y.
	\end{align*}
	If $|y|\leq\frac{|h|}2$, then $|y+h|\geq|h|-|y|\geq\frac{|h|}2$, so $A_1$ can be handled in the same way as $A_2$. For the latter, we can now use \eqref{ylarge} and \eqref{int-bound-2} to obtain the desired bound:
	\begin{align*}
	A_2&\lec \frac{\del}{\tau_n^2} \lp \int_0^{2\del} 2\,\dd s + \int_{2\del}^\infty \frac{\del}{s-\del}\ee^{-\frac\la 2 (s-\del)}\,\dd s\rp \lec \frac{\del^2}{\tau_n^2} (1+|\log \del|)\\ &\lec \del^{\frac32}|\log\del|.
	\end{align*}
	
	Next, we examine the case $h=0$. Regarding the first statement, since $\Pi^n_{r,0}([0,\infty)\times\R^d\times\R^d) = \Ga^n_r$ by definition and $\Ga^n_r\to\Ga_r$ by \eqref{Ganr-speed}, the claim follows from the Portmanteau theorem if we can show that $|\Pi^n_{r,0}|(A)\to 0$ for all closed sets $A\subseteq [0,\infty)\times\R^d\times\R^d$ that do not contain the origin. In this case, $A$ is bounded away from $y=0$, $z=0$ or $t=0$. As $G$ is smooth outside the origin,  we have $|G_y(s)-G_y(s-\del)|/\del \leq \ov G(s,y)$, $|G_z(s+r\del)-G_z(s+(r-1)\del)|/\del \leq \ov G(s,z)$, or both on $A$, where $\ov G$ is a bounded and rapidly decreasing function that is independent of $n$. By symmetry, let us suppose that the second bound holds. Then
	\beq\label{z-away}\begin{split}  \frac{\tau_n^2}{\del}|\Pi^n_{r,0}|(A)&=\iiint_A |G_y(s)- G_y(s-\del)|\\
		&\quad\times \frac{|G_z(s+r\del)-G_z(s+(r-1)\del)|}{\del}\,\dd s\,\La(\dd y,\dd z) \\
		&\leq\iiint |G_y(s)-G_y(s-\del)| \ov G(s,z)\,\dd s\,\La(\dd y,\dd z) \lec 1, \end{split}\eeq
	and (1) follows because $\del/\tau_n^2\to0$ by \eqref{taun-asymp}.  
	
	For statement (2), we again leave the cases $r=0$ and $r=1$ to the reader and assume $r\geq2$. Also, we exclude the case $F=\delta_0$, which can be handled along the following lines and is much simpler in fact. Then we decompose $|\Pi^n_{r,0}|([0,\infty)\times\R^d\times\R^d) = B_1+B_2+B_3$ where
	\begin{align*}
	B_1&=\frac{c_\al}{\tau^2_n}\iiint_0^{\del} \frac{G_y(s)|G_z(s+r\del)-G_z(s+(r-1)\del)|}{|z-y|^\al}\,\dd s\,\dd y\,\dd z,\\
	B_2&=\frac{c_\al}{\tau^2_n}\iiint \frac{|G_y(s+\del)-G_y(s)||G_z(s+(r+1)\del)-G_z(s+r\del)|}{|z-y|^{\al}}\\
	&\quad\times\bone_{|z-y|\leq \sqrt{\del}}\,\dd s\,\dd y\,\dd z,\\
	B_3&=\frac{c_\al}{\tau^2_n}\iiint \frac{|G_y(s+\del)-G_y(s)||G_z(s+(r+1)\del)-G_z(s+r\del)|}{|z-y|^{\al}}\\
	&\quad\times\bone_{|z-y|> \sqrt{\del}}\,\dd s\,\dd y\,\dd z.
	\end{align*}
	For $B_1$, we first use Parseval's identity to compute
	\beq\label{convolution} \begin{split}
		\quad	&c_\al \int_{\R^d} G_y(s)|z-y|^{-\al}\,\dd y = c_\al \lv \int_{\R^d} \calf G(s,\cdot)(\xi)\ee^{2\pi \ii \xi\cdot z} \calf |\cdot|^{-\al}(\xi)\,\dd \xi\rv\\
		&\quad=\lv\int_{\R^d} \ee^{-(\la+2\pi^2\kappa |\xi|^2)s}\ee^{2\pi \ii \xi\cdot z}|\xi|^{\al-d}\,\dd \xi\rv \leq \int_{\R^d} \ee^{-(\la+2\pi^2\kappa |\xi|^2)s}|\xi|^{\al-d}\,\dd \xi\\
		&\quad=|\cals^{d-1}|\ee^{-\la s}\int_0^\infty \ee^{-2\pi^2\kappa r^2 s}r^{\al-1}\,\dd r= \frac{|\cals^{d-1}|\Ga(\frac\al2)\ee^{-\la s}}{2(2\pi^2\kappa s)^{\frac\al 2}}. 
	\end{split}\raisetag{-3\baselineskip}\eeq
	From this, together with \eqref{int-bound-2} and \eqref{taun-asymp}, we obtain
	\begin{align*}
	B_1&\lec \frac{1}{\tau_n^2} \int_0^{\del} s^{-\frac\al 2} \frac{\del}{s+(r-1)\del}\ee^{-\frac \la 2 (s+(r-1)\del)}\,\dd s\leq \frac{1}{\tau_n^2 (r-1)} \int_0^{\del} s^{-\frac\al 2}\,\dd s\\ &\lec \frac{1}{r-1}.
	\end{align*}
	For $B_2$, we use the fact that 
	\beq\label{der} \lv\frac{\partial}{\partial t} G_y(t)\rv\lec (1\vee t){t^{-1-\frac d 2}} \ee^{-\la t}\lec {t^{-1-\frac d 2}} \ee^{-\frac\la 2 t},\eeq uniformly in $y\in\R^d$. Hence, by the mean value theorem, we have that $|G_z(s+(r+1)\del)-G_z(s+r\del)|\lec (s+r\del)^{-1- d/ 2}\ee^{- \la  (s+r\del)/2}\del$. So together with \eqref{int-bound-2}, we obtain
	\begin{align*}
	B_2&\lec \frac{1}{\tau^2_n}\int_0^\infty \frac{\del}{s+r\del} (s+r\del)^{-\frac d 2} \ee^{-\frac\la 2(s+r\del)}\\
	&\quad\times\int_{\R^d} |G_y(s+\del)-G_y(s)|\,\dd y \int_{|z|\leq \sqrt{\del}} {|z|^{\al}}\,\dd z \,\dd s\\
	&\lec \frac{1}{\tau^2_nr}\del^{\frac{d-\al}{2}} \lp \int_0^{\del} s^{-\frac d 2}\,\dd s + \int_{\del}^\infty s^{-\frac d 2}\ee^{- \la  s} \frac{\del}{s}\,\dd s\rp \lec \frac 1 r.
	\end{align*}
	Finally, again by \eqref{int-bound-2}, 
	\begin{align*}
	B_3&\lec \frac{\del^{-\frac \al 2}}{\tau_n^2}\int_0^\infty \int_{\R^d} |G_y(s+\del)-G_y(s)|\,\dd y\\ &\quad\times\int_{\R^d} |G_z(s+(r+1)\del)-G_z(s+r\del)|\,\dd z\,\dd s\\
	&\lec\frac{\del^{-\frac \al 2}}{\tau_n^2} \left( \int_0^{\del} \frac{\del}{s+r\del} \ee^{-\frac\la 2 (s+r\del)}\,\dd s \right.\\
	&\quad+\left. \int_{\del}^\infty \frac{\del}s \ee^{-\frac\la 2 s}\frac{\del}{s+r\del} \ee^{-\frac\la 2 (s+r\del)}\,\dd s\right)\\
	&\lec \frac{\del^{-\frac \al 2}}{\tau_n^2 } \lp \frac{\del}{r} + \int_{\del}^1 \frac{\del}s \frac{\del}{s+r\del} \,\dd s+\del^2\int_1^\infty \ee^{-\la s}\,\dd s \rp.
	\end{align*}
	By a change of variables, we have 
	\begin{align*}\int_{\del}^1 \frac{\del}s \frac{\del}{s+r\del} \,\dd s &= \del\int_1^{1/\del} \frac{1}{s(s+r)}\,\dd s \leq \del \int_1^\infty \frac{1}{s(s+r)}\,\dd s\\ &= \del \frac{\log(r+1)}{r}.\end{align*}
	It follows that $B_3 \lec \frac{\log(r+1)}{r}$, which proves statement (2) for $h=0$.
	\epr

	\blem\label{y2} For every $r\in\N_0$, we have
	\beq\label{pin2} \iiint_0^\infty (|y|^2+|z|^2)\,|\Pi^n_{r,0}|(\dd s,\dd y,\dd z) = o(\delh)\as. \eeq
	\elem  
	\bpr Since the arguments involved are completely analogous, we only consider the integral against $|y|^2$. We first examine the case where $F\neq\delta_0$, and choose numbers $\beta$, $\ga$, and $\delta$ satisfying the constraints
	\beq\label{constraints} \frac 2 5<\beta<\frac 12-\frac{2\delta-1}{2\al},\qquad \frac14<\ga<\frac\delta 2, \qquad \frac12<\delta<1,  \eeq
	which is clearly possible if $\delta$ is sufficiently close to $\frac12$.
	Because we have $\iiint_0^\infty |y|^2\bone_{|y|\leq \del^\ga}\,|\Pi^n_{r,0}|(\dd s,\dd y,\dd z) \lec \del^{2\ga}$ by \eqref{cross-bound-3} and $2\ga>\frac12$ by \eqref{constraints}, we may restrict the domain of integration to $|y|>\del^\ga$ in the following. Then
	\[ \iiint |y|^2\bone_{|y|>\del^\ga}\,|\Pi^n_{r,0}|(\dd s,\dd y,\dd z)=c_\al(C_1+C_2+C_3), \]
	where
	\begin{align*}
	C_1&=\iiint_0^{\del} \frac{G_y(s)|G_z(s+r\del)-G_z(s+(r-1)\del)|}{\tau_n^2|z-y|^{\al}}|y|^2\bone_{|y|>\del^\ga}\,\dd s\,\dd y\,\dd z,\\
	C_2&=\iiint  \frac{|G_y(s+\del)-G_y(s)||G_z(s+(r+1)\del)-G_z(s+r\del)|}{\tau_n^2|z-y|^{\al}}|y|^2\\
	&\quad\times\bone_{|y|>\del^\ga,|z-y|>\del^\beta}\,\dd s\,\dd y\,\dd z,\\
	C_3&=\iiint \frac{|G_y(s+\del)-G_y(s)||G_z(s+(r+1)\del)-G_z(s+r\del)|}{\tau_n^2|z-y|^{\al}}|y|^2\\
	&\quad\times\bone_{|y|>\del^\ga,|z-y|\leq\del^\beta} \,\dd s\,\dd y\,\dd z.
	\end{align*}
	We have $\int_{\R^d} |G_z(s+r\del)-G_z(s+(r-1)\del)||z-y|^{-\al}\,\dd z \lec s^{-\al/ 2}$ by \eqref{convolution}, and furthermore, 
	\beq\label{dens}\int_{\R^d} G_y(s)|y|^2\,\dd y \leq s\eeq 
	because $G_y(s)$ is a Gaussian density  (up to an exponential factor). Hence,
	\[ C_1 \lec \frac{1}{\tau_n^2} \int_0^{\del} s^{1-\frac\al 2}\,\dd s\lec \Del^{2-\frac\al 2 -1+\frac\al 2} = \del = o(\del^{\frac 12}). \]
	
	Furthermore, $C_2\leq C_{21}+C_{22}+C_{23}$ where
	\begin{align*} C_{21}&=  \iiint_0^{\del^\delta} \frac{|G_y(s+\del)-G_y(s)||G_z(s+(r+1)\del)-G_z(s+r\del)|}{\tau_n^2\del^{\al\beta}}|y|^2\\
	&\quad\times\bone_{|y|>\del^\ga}\,\dd s\,\dd y\,\dd z,\\
	C_{22}&= \iiint_{\del^\delta}^\infty \frac{ |G_y(s+\del)-G_y(s)||G_z(s+(r+1)\del)-G_z(s+r\del)|}{\tau_n^2\del^{\al\beta}}|y|^2\\
	&\quad\times\bone_{\del^\ga<|y|\leq 1}\,\dd s\,\dd y\,\dd z,\\
	C_{23}&= \iiint_{\del^\delta}^\infty \frac{ |G_y(s+\del)-G_y(s)||G_z(s+(r+1)\del)-G_z(s+r\del)|}{\tau_n^2\del^{\al\beta}}|y|^2\\
	&\quad\times\bone_{|y|> 1}\,\dd s\,\dd y\,\dd z. \end{align*}
	Bounding the $\dd z$-integral by $2$ and substituting $r^2/(2\kappa s)\mapsto u$, we obtain
	\begin{align*}
	C_{21}&\lec \frac{\del^{-\al\beta}}{\tau_n^2} \int_0^{\del^\delta+\del}   \int_{\R^d} G_y(s)|y|^2\bone_{|y|>\del^\ga}\,\dd y\,\dd s\\
	&\lec\frac{\del^{-\al\beta}}{\tau_n^2} \int_0^{2\del^\delta}  s^{-\frac{d}{2}} \int_{\del^\ga}^\infty \ee^{-\frac{r^2}{2\kappa s}}r^{d+1}\,\dd r\,\dd s\\
	&=C\frac{\del^{-\al\beta}}{\tau_n^2} \int_0^{2\del^\delta}  s \int_{\frac{\del^{2\ga}}{2\kappa s}}^\infty \ee^{-u}u^{\frac{d}{2}}\,\dd u\,\dd s\\
	&=C \frac{\del^{-\al\beta}}{\tau_n^2} \int_{\frac{\del^{2\ga-\delta}}{4\kappa}}^\infty \lp \int_{\frac{\del^{2\ga}}{2\kappa u}}^{2\del^\delta}s \,\dd s\rp \ee^{-u} u^{\frac d 2}\,\dd u\\
	&\lec \frac{\del^{-\al\beta}}{\tau_n^2} \del^{2\delta} \int_{\frac{\del^{2\ga-\delta}}{4\kappa}}^\infty\ee^{-u} u^{\frac d 2}\,\dd u.
	\end{align*}
	Since $\ga<\frac\delta 2$ by \eqref{constraints}, the last integral goes to $0$ as $n\to\infty$ at an exponential speed, and $C_{21}=o(\Del^{1/2})$ follows. Regarding $C_{22}$, we use \eqref{int-bound-2} on both the $\dd y$- and $\dd z$-integral and derive
	\[ C_{22} \lec\frac{\del^{-\al\beta}}{\tau_n^2}  \int_{\del^\delta}^\infty  \frac{\del^2}{s^2}\ee^{-\la s}\,\dd s \lec \del^{-\al\beta-1+\frac\al 2+2-\delta}.\]
	This exponent is larger than $\frac12$ if and only if $\beta<\frac12 -\frac{2\delta-1}{2\al}$, which is true by \eqref{constraints} and yields $C_{22}=o(\Del^{1/2})$ as well. For $C_{23}$, because $|y|>1$, we have \beq\label{ylarge} |G_y(s+\del)-G_y(s)| \lec \ov G(s,y) \del,\eeq where $\ov G$ has the same properties as in the proof of Lemma~\ref{delta}. Together with \eqref{int-bound-2} on the $\dd z$-integral, this shows
	\[ C_{23}\lec \frac{\del^{1-\al\beta}}{\tau_n^2} \int_{\del^{\delta}}^\infty \frac{\del}{s}\ee^{-\frac\la 2 s}\,\dd s \lec \del^{2-\al\beta-1+\frac\al 2}|\log\del|. \]
	Comparing this exponent with that of $C_{22}$, we see that $C_{23}=o(\Del^{1/2})$.
	
	Next, we have $C_3\leq C_{31}+C_{32}$ where
	\begin{align*}
	C_{31}&= \iiint_0^{\del^\eta} \frac{|G_y(s+\del)-G_y(s)||G_z(s+(r+1)\del)-G_z(s+r\del)|}{\tau_n^2|z-y|^{\al}}|y|^2\\
	&\quad\times\bone_{|z-y|\leq \del^\beta}\,\dd s\,\dd y\,\dd z,\\
	C_{32}&= \iiint_{\del^\eta}^\infty \frac{|G_y(s+\del)-G_y(s)||G_z(s+(r+1)\del)-G_z(s+r\del)|}{\tau_n^2|z-y|^{\al}}|y|^2\\
	&\quad\times\bone_{|z-y|\leq \del^\beta}\,\dd s\,\dd y\,\dd z,
	\end{align*}
	and where $\eta=\frac 12$ in dimension $d=1$ and $\eta=\frac 4 5$ in dimensions $d\geq2$. 
	By an elementary calculation, \beq\label{elementary}\int_{|z-y|\leq \del^\beta} |z-y|^{-\al}\,\dd z = \int_{|z|\leq \del^\beta} |z|^{-\al}\,\dd z \lec \del^{\beta(d-\al)}.\eeq Thus, if we apply the mean value theorem to $|G_z(s+(r+1)\del)-G_z(s+r\del)|$ and use \eqref{der} and \eqref{dens}, then in the case $d=1$, we have
	\begin{align*}
	C_{31}&\lec \frac 1{\tau_n^2} \int_0^{\delh} \lp\int_{\R} |G_y(s+\del)-G_y(s)||y|^2\,\dd y\rp\frac{\del}{s^{\frac 1 2+1}}\del^{\beta(1-\al)}\,\dd s\\ 
	&\lec \frac{\del^{1+\beta(1-\al)}}{\tau_n^2} \int_0^{2\delh} s^{-\frac 1 2}\,\dd s\lec \del^{1+\beta(1-\al)-1+\frac\al 2+\frac14}.
	\end{align*}
	The exponent is larger than $\frac12$ if and only if $(\beta-\frac14)+\al(\frac12-\beta)>0$, which is true by our choice of $\beta$. For $d\geq2$, observe that \eqref{convolution} yields the bound $\int_{\R^d} |G_z(s+(r+1)\del)-G_z(s+r\del)||z-y|^{-\al}\,\dd z\lec s^{-\al/ 2}$. Hence, for $d\geq2$, 
	\begin{align*}
	C_{31}&\lec \frac {1}{\tau_n^2} \int_0^{\del^{\frac 4 5}} \lp\int_{\R^d} |G_y(s+\del)-G_y(s)||y|^2\,\dd y\rp s^{-\frac \al 2}\,\dd s\\ 
	&\lec \frac{1}{\tau_n^2} \int_0^{2\del^{\frac 4 5}} s^{1-\frac \al 2}\,\dd s\lec \del^{-1+\frac\al 2+\frac45(2-\frac\al 2)}=\del^{\frac 3 5+\frac \al {10}} = o(\delh).
	\end{align*}
	Concerning $C_{32}$ for $d=1$, we estimate similarly as for $C_{31}$ but use \eqref{int-bound-2} and \eqref{ylarge} on the $\dd y$-integral, which leads to 
	\begin{align*}
	C_{32}&\lec \frac {\del^{1+\beta(1-\al)}}{\tau_n^2} \int_{\del^{\frac 12}}^\infty \lp\int_{\R} |G_y(s+\del)-G_y(s)||y|^2\,\dd y\rp \frac{\ee^{-\frac\la 2 s}}{s^{\frac 1 2 +1}}\,\dd s \\
	&\lec \frac{\del^{1+\beta(1-\al)}}{\tau_n^2} \int_{\delh}^\infty \lp \frac{\del}{s}+\del\rp \frac{\ee^{-\frac\la 2 s}}{s^{\frac 1 2 +1}}\,\dd s\\ 
	&\lec \del^{1+\beta(1-\al)-1+\frac\al 2+1-\frac34}\lec\del^{\beta+\frac 14 +\al(\frac12 -\beta)} = o(\delh)
	\end{align*}
	because $\frac14<\beta<\frac12$. In dimensions $d\geq2$, we use the same estimate on the $\dd z$-integral but bound the $\dd y$-integral as for $C_{31}$, which yields
	\begin{align*}
	C_{32}&\lec \frac{\del^{1+\beta(d-\al)}}{\tau_n^2} \int_{\del^{\frac 4 5}}^\infty \lp \int_{\R^d} |G_y(s+\del)-G_y(s)||y|^2\,\dd y\rp \frac{\ee^{-\frac\la 2 s}}{s^{\frac d 2 +1}}\,\dd s\\
	&\lec \frac{\del^{1+\beta(d-\al)}}{\tau_n^2} \int_{\del^{\frac 4 5}}^\infty  \frac{\ee^{-\frac\la 2 s}}{s^{\frac d 2}}\,\dd s\lec \del^{1+\beta d-\al\beta-1+\frac\al 2+\frac 45 (1-\frac d 2)}
	\end{align*}
	($\times |\log \del|$ if $d=2$).
	The exponent of $\del$ is $d(\beta-\frac 2 5)+\al(\frac 12-\beta)+\frac 4 5>\frac 12$ because $\frac 2 5<\beta<\frac 12$.
	
	We are left to consider the case $d=1$ with $F=\delta_0$. Then the left-hand side of \eqref{pin2} equals $D_1+D_2$ where $D_1$ and $D_2$ are the $\dd s\,\dd y$-integrals of 
	\begin{align*}
	\frac{|G_y(s)-G_y(s-\del)||G_y(s+r\del)-G_y(s+(r-1)\del)|}{\tau_n^2}|y|^2
	\end{align*}	
	over $[0,\sqrt{\Del}]\times\R$ and $[\sqrt{\Del},\infty)\times\R$, respectively.
	Using \eqref{der} and \eqref{dens}, we obtain
	\[ D_1\lec \frac1{\tau_n^2}\int_0^{\delh} \frac{\del}{s^{\frac12+1}}s\,\dd s \lec \del^{1+\frac14-\frac12} = o(\delh). \]
	Furthermore, using Plancherel's theorem and \eqref{fourierG}, we have for all $s>0$,
	\beq\label{square} \begin{split}
		&\int_{\R} |G_y(s+\del)-G_y(s)|^2\,\dd y\\
		&\quad= \ee^{-2\la s} \int_\R \lp\ee^{-\la \del-2\pi^2\kappa (s+\del) \xi^2}-\ee^{-2\pi^2\kappa s \xi^2}\rp^2 \,\dd \xi \\
		&\quad=\ee^{-2\la s} \int_\R \ee^{-4\pi^2\kappa s \xi^2}\lp \ee^{-\la \del-2\pi^2\kappa \del \xi^2}-1\rp^2 \,\dd \xi\\
		&\quad\lec \ee^{-2\la s}\del^2 \int_\R \ee^{-4\pi^2\kappa s \xi^2}\lp \la +2\pi^2\kappa  \xi^2\rp^2 \,\dd \xi\\
		&\quad\lec \del^2\ee^{-2\la s} (s^{-\frac12}\vee s^{-\frac 52})\lec\del^2\ee^{-\la s} s^{-\frac 52}.
	\end{split}\eeq
	Hence, if we apply  the Cauchy--Schwarz inequality on the $\dd y$-integral in $D_2$ and use the formula that we have just derived for $|y|\leq1$ and \eqref{ylarge} for $|y|>1$, we can show that 
	\begin{align*}
	D_2\lec \frac{1}{\tau_n^2} \lp \int_{\delh}^\infty \del^2\ee^{-\la s} s^{-\frac 52}\,\dd s +\del^2\rp= \del^{-\frac12+2-\frac34}+\del^{-\frac12+2} = o(\delh),
	\end{align*}
	which completes the proof of the lemma.
	\epr
	
	\blem\label{nu} For all $\theta\in(0,1)$, we have
	\beq\label{pinrate} |\Pi^n_{0,0}|((\Delta_n^{1-\theta},\infty)\times\R^d\times\R^d) \lec \Delta_n^{\nu \theta} \eeq
	with $\nu=1+\frac\al 2$.
	\elem
	\bpr First, consider the case $F\neq \delta_0$. Then the left-hand side of \eqref{pinrate} equals $E_1+E_2$ where
	\begin{align*}
	E_1&=\frac{c_\al}{\tau^2_n}\iiint_{\del^{1-\theta}-\del}^\infty \frac{|G_y(s+\del)-G_y(s)||G_z(s+\del)-G_z(s)|}{|z-y|^{\al}}\\
	&\quad\times\bone_{|z-y|\leq \del^{\frac{1-\theta}{2}}}\,\dd s\,\dd y\,\dd z,\\
	E_2&=\frac{c_\al}{\tau^2_n}\iiint_{\del^{1-\theta}-\del}^\infty \frac{|G_y(s+\del)-G_y(s)||G_z(s+\del)-G_z(s)|}{|z-y|^{\al}}\\
	&\quad\times\bone_{|z-y|> \del^{\frac{1-\theta}{2}}}\,\dd s\,\dd y\,\dd z.
	\end{align*}
	From \eqref{int-bound-2}, \eqref{der}, and \eqref{elementary}, it follows that
	\begin{align*}
	E_1&\lec \frac{1}{\tau_n^2} \int_{\frac12 \del^{1-\theta}}^\infty \frac{\del}{s}\ee^{-\la s} \frac{\del}{s^{1+\frac d 2}} \del^{\frac{1-\theta}{2}(d-\al)}\,\dd s \lec \del^{2+\frac{1-\theta}{2}(d-\al)-1+\frac\al 2-(1-\theta)(1+\frac d 2)}\\ 
	&=\del^{\nu \theta}.
	\end{align*}
	With similar considerations as for the term $A_3$ in the proof of Lemma~\ref{delta}, we obtain
	\begin{align*}
	E_2&\lec\frac{\del^{-\frac{1-\theta}{2}\al}}{\tau_n^2}\int_{\frac12 \del^{1-\theta}}^\infty \frac{\ee^{-\la s}\del^2}{s^2}\,\dd s \lec \del^{-\frac{1-\theta}{2}\al-1+\frac\al 2 +2-(1-\theta)} = \del^{\nu a},
	\end{align*}

	It remains to consider $d=1$ and $F=\delta_0$. Using \eqref{square}, we obtain 
	\begin{align*}
	|\Pi^n_{0,0}|((\Delta_n^{1-\theta},\infty)\times\R\times\R) &= \frac1{\tau_n^2} \int_{\del^{1-\theta}}^\infty \int_{\R^d} |G_y(s)-G_y(s-\del)|^2\,\dd y\,\dd s\\
	&\lec  \frac{\del^2}{\tau_n^2} \int_{\frac12\del^{1-\theta}}^\infty \frac{\ee^{-\la s}}{s^{\frac52}}\,\dd s \lec \del^{2-\frac12-\frac32(1-\theta)} = \del^{\frac32 \theta},
	\end{align*}
	which is the claim with $\nu=1+\frac\al2 = \frac32$.
	\epr

	\section{Details for the proof of Theorem~\ref{LLN}}\label{details-LLN}
	
	\bpr[Proof of Lemma~\ref{sibounded}]  Let $\si_m(\om,t,x) = (\si(\om,t,x) \vee (-m))\wedge m$ for $m\in\N$ and define $Y^m_x$ like $Y_x$ but with volatility process  $\si_m$ instead of $\si$. It is clear that $\si_m$ also satisfies Assumptions~\ref{AssLLN}2 and \ref{AssLLN}3 as soon as $\si$ does. We decompose $V^n_f(Y_x,t) - V_f(Y_x,t)$ as
	\beq\label{split3}\begin{split}  &\Big( V^n_f(Y_x,t)-V^n_f(Y_x^m,t)\Big) + \Big( V^n_f(Y_x^m,t) - V_f(Y_x^m,t) \Big) \\
		&\quad+ \Big( V_f(Y_x^m,t) - V_f(Y_x,t)\Big), 
	\end{split}\eeq
	and assume that Theorem~\ref{LLN} holds whenever $\si$ is bounded. Then for every $m\in\N$, the second term in \eqref{split3} converges in $L^1$ to $0$, uniformly on $[0,T]$.
	
	Next, since $f$ is continuous with $f(z)\lec 1+|z|^p$, the function $\mu_f$ is continuous on $\R^N$, so by Assumption~\ref{AssLLN}2, $\bbe[|\mu_f(\un\si^2(t,x))|]\lec 1+\bbe[|\un\si(t,x)|^p]$, which is uniformly bounded in $t$ and $x$. Hence, the dominated convergence theorem shows that also $V_f(Y_x^m,t) - V_f(Y_x,t) \to0$ in $L^1$, uniformly in $t\in[0,T]$.

	Finally, by Assumption~\ref{AssLLN}1,  for $0<\delta<1<A$, there are finite numbers $\Phi(A)$, $\Phi^\prime_A(\delta)$, and $\Phi^{\prime\prime}(A)$ such that $\Phi^\prime_A(\delta) \to 0$ as $\delta\to0$, $\Phi^{\prime\prime}(A)\to 0$ as $A\to \infty$, and
	\begin{align*}
	|z|\leq 2A &\implies |f(z)|\leq \Phi(A),\\
	|z|\leq 2A,~|z'|\leq \delta &\implies |f(z+z')-f(z)|\leq \Phi^\prime_A(\delta), \\
	|z|>A &\implies |f(z)|\leq \Phi^{\prime\prime}(A) |z|^p
	\end{align*}
	hold for all $z,z'\in\R^{N\times L}$;
	cf.\ the proof of Lemma~3.4.6 in \citeB{Jacod12-2}. This implies 
	\beq\label{fdiff} |f(z+z')-f(z)|\leq C\lp\Phi^\prime_A(\delta) + \frac{\Phi(A)|z'|}{\delta}+\Phi^{\prime\prime}(A)(|z|^p + |z'|^p)\rp \eeq
	for $z,z'\in\bbr^{N\times L}$ and some  $C>0$ only depending on $p$. Thus, we obtain
	\begin{align*} \left|f\lp \frac{\un \Delta^n_i Y_x}{\tau_n} \rp - f\lp \frac{\un\Delta^n_i Y^m_x}{\tau_n} \rp\right| &\lec \Phi^\prime_A(\delta) + \frac{\Phi(A)|\un\Delta^n_i Y_x-\un\Delta^n_i Y^m_x|}{\delta\tau_n}\\
	&\quad+\frac{\Phi^{\prime\prime}(A)(|\un\Delta^n_i Y_x|^p + |\un\Delta^n_i Y^m_x|^p)}{\tau_n^p}.
	\end{align*}
	Using \eqref{cross-bound-3} and Assumption~\ref{AssLLN}2, we deduce the estimate
	\begin{align*} \bbe[|\un\Delta^n_i Y_x-\un\Delta^n_i Y^m_x|^2] &=\bbe\lb\lv\iint \un\Delta^n_i G_{x,y}(s) (\si(s,y)-\si_m(s,y)) \,W(\dd s,\dd y)\rv^2\rb\\
	&\lec \sup_{(s,y)\in\bbr\times\bbr^d} \bbe[(\si-\si_m)^2(s,y)]\tau_n^2\\
	& = \sup_{(s,y)\in\bbr\times\bbr^d} \bbe[\si^2(s,y)\bone_{\{|\si(s,y)|>m\}}]\tau_n^2\\
	&\leq m^{-\eps} \sup_{(s,y)\in\R\times\bbr^d} \bbe[|\si(s,y)|^{2+\eps}]\tau_n^2\lec \tau_n^2 m^{-\eps},
	\end{align*}
	and similarly,
	$
	\bbe[|\un\Delta^n_i Y_x|^p] + \bbe[|\un\Delta^n_i Y^m_x|^p] \lec \sup_{(s,y)\in\bbr\times\bbr^d} \bbe[|\si(s,y)|^p]\tau_n^p \lec 
	\tau_n^p$.
	Thus,
	\begin{align*} \bbe\lb\supt \left|V^n_f(Y_x,t)-V^n_f(Y^m_x,t)\right|\rb&\leq \Del\sumT \bbe\lb\left|f\lp \frac{\un\Delta^n_i Y_x}{\tau_n} \rp - f\lp \frac{\un\Delta^n_i Y^m_x}{\tau_n} \rp\right|\rb\\
	& \lec  \Phi'_A(\delta)+\frac{\Phi(A)}{\delta m^{\eps/2}}+\Phi''(A). \end{align*}
	The right-hand side does not depend on $n$. Furthermore, it converges to $0$ if we let $m\to\infty$, $\delta\to0$, and $A\to\infty$ in this order, which completes the proof of Lemma~\ref{sibounded}.
	\epr 
	
	\bpr[Proof of Lemma~\ref{lemma-eps}] Similarly to the proof of Lemma~\ref{sibounded}, the assertion is proved once we can show that
	\beq\label{first}\begin{split}
		\lim_{n\to\infty} \sup_{i=1,\ldots,T^*_n} \frac{\bbe[|\un\Delta^n_i Y_x-\alpha^{n,i,\eps}_x|]}{\tau_n} &= 0,\\
		\limsup_{n\to\infty} \sup_{i=1,\ldots,T^*_n} \frac{\bbe[|\alpha^{n,i,\eps}_x|^p]}{\tau_n^p} &< \infty. 
	\end{split}\eeq
	The second part in \eqref{first} is an immediate consequence of the BDG inequality and \eqref{cross-bound-3}.
	Regarding the first statement in \eqref{first}, we have by definition
	\[ \un\Delta^n_i Y_x-\alpha^{n,i,\eps}_x= \iint  \un\Delta^n_i G_{x,y}(s)\bone_{s\leq i\del-\eps} \si(s,y)\,W(\dd s,\dd y).\]
	Therefore, another application of \eqref{cross-bound-3} yields
	\beq\label{expr-1} \frac{\bbe[|\un\Delta^n_i Y_x-\al^{n,i,\eps}_x|^2]}{\tau_n^2} \lec \sup_{(s,y)\in\bbr\times\bbr^d} \bbe[|\si(s,y)|^2] |\Pi^n_{0,0}|([\eps,\infty)\times\R^d\times\R^d), \eeq
	which tends to $0$ as $n\to\infty$ by Lemma~\ref{delta}, proving \eqref{first}. 
	\epr
	
	\bpr[Proof of Lemma~\ref{LLN-disc}]
	Using again the estimate \eqref{fdiff}, the claim follows from
	\begin{align}
	\lim_{\eps\to0}	\limsup_{n\to\infty} \sup_{i=1,\ldots,T^*_n} \frac{\bbe[|\alpha^{n,i,\eps}_x-\wh\alpha^{n,i,\eps}_x|]}{\tau_n} &= 0, \label{firstterm-2} \\
	\limsup_{\eps\to0}	\limsup_{n\to\infty} \sup_{i=1,\ldots,T^*_n} \frac{\bbe[|\alpha^{n,i,\eps}_x|^p]+\bbe[|\wh\alpha^{n,i,\eps}_x|^p]}{\tau_n^p} &< \infty. \label{secondterm-2}
	\end{align}
	While \eqref{secondterm-2} follows from the BDG inequality and \eqref{cross-bound-3}, \eqref{firstterm-2} holds because
	\begin{align*}
	\frac{\bbe[|\alpha^{n,i,\eps}_x-\wh\alpha^{n,i,\eps}_x|^2]}{\tau_n^2}\lec \sup\left\{\bbe[|\si(t,y)-\si(s,y)|^2]\colon |t-s|\leq \eps,~y\in\R^d\right\},
	\end{align*}
	which converges to $0$ as $\eps\to0$, uniformly in $n$ and $i$, by Assumption~\ref{AssLLN}3.
	\epr
	
	\bpr[Proof of Lemma~\ref{centering}] The expression in braces has zero expectation and by Assumption~\ref{AssLLN}1 and the boundedness of $\sigma$, a variance that is uniformly bounded in $i$, $n$, and $\eps$. Furthermore, as soon as $\del\leq \eps$, two such expressions, with indices $i\neq j$, are uncorrelated as soon as $|i-j|>\eps/\del+L-1$. So the second moment in \eqref{expr-2} is $\lec \del^2[t/\Delta_n](\eps/\Delta_n+L-1)\lec \eps$, and \eqref{expr-2} follows.
	\epr
	
	\bpr[Proof of Lemma~\ref{cond-exp-conv}] Conditionally on $\calf_{i\del-\eps}$, the matrix $\wh \al^{n,i,\eps}_x$ has a multivariate normal distribution with mean $0$ and covariances
	\beq\label{vni}\begin{split} (\beta^{n,i,\eps}_x)_{j_1k_1,j_2k_2}&=\cov((\wh \al^{n,i,\eps}_x)_{j_1 k_1}, (\wh \al^{n,i,\eps}_x)_{j_2 k_2})\\
		&=\iiint_0^\infty \si(i\del-\eps,x_{j_1}-y)\si(i\del-\eps,x_{j_1}-z)\\
		&\quad\times\bone_{s<\eps+(k_1-1)\del} \,\Pi^n_{k_2-k_1,x_{j_2}-x_{j_1}}(\dd s,\dd y,\dd z) \end{split}\eeq
	for $j_1,j_2\in\{1,\ldots,N\}$ and $k_1,k_2\in\{1,\ldots,L\}$ with $k_1\leq k_2$. It follows that \beq\label{muf}\bbe\left[ f\left( \frac{\wh\alpha^{n,i,\eps}_x}{\tau_n} \right)\,\Big|\, \calf_{i\del-\eps} \right]=\un\mu_f(\beta^{n,i,\eps}_{x})\eeq where $\un\mu_f$ is the function that maps $(\beta_{j_1k_1,j_2k_2})_{j_1,j_2,k_1,k_2=1}^{N,N,L,L}$ to $\bbe[f(\un Z)]$, and where $\un Z=(\un Z_{jk})_{j,k=1}^{N,L}$ has a multivariate normal distribution with mean $0$ and $\cov(\un Z_{j_1 k_1},\un Z_{j_2 k_2})=\beta_{j_1 k_1,j_2 k_2}$.
	
	Hence, the difference in \eqref{conv} equals $A^{n,\eps,1}_x(t)+A^{n,\eps,2}_x(t)+A^{n,\eps,3}_x(t)+A^{n,4}_x(t)$
	where
	\begin{align*}
	A^{n,\eps,1}_x(t)&= \sumt \int_{(i-1)\Del}^{i\Del} \lp\un\mu_f(\beta^{n,i,\eps}_x)-\un\mu_f(\beta^{n,\eps}_x(r))\rp \,\dd r,\\
	A^{n,\eps,2}_x(t)&=\int_0^{\del t^*_n} \lp\un \mu_f(\beta^{n,\eps}_x(r)) - \mu_f(\un\si^2(r-\eps,x))\rp \,\dd r,\\
	A^{n,\eps,3}_x(t)&=\int_0^{\del t^*_n} \lp \mu_f(\un\si^2(r-\eps,x)) -\mu_f(\un\si^2(r,x))\rp \,\dd r,\\
	A^{n,4}_x(t)&=\int_{\del t^*_n}^t \mu_f(\un\si^2(r,x))\,\dd r,
	\end{align*}
	and where $(\beta^{n,\eps}_x(r))_{j_1k_1,j_2k_2}$ is defined by the right-hand side of \eqref{vni} with $i\del$ replaced by $r$. Now, by Assumption~\ref{AssLLN}1, $\beta\mapsto\un\mu_f(\beta)$ is a continuous function with $\un\mu_f(\beta)=o(|\beta|^{p/ 2})$. Hence, if we apply the estimate \eqref{fdiff} on $\un\mu_f$, we can bound $\bbe[|A^{n,\eps,1}_x(t)|]$ by a constant times
	\begin{align*}
	& \Phi'_A(\delta)\del[t/\del]+\frac{\Phi(A)}{\delta}\sumt \int_{(i-1)\Del}^{i\Del} \bbe[|\beta^{n,i,\eps}_x-\beta^{n,\eps}_x(r)|] \,\dd r\\
	&\quad+\Phi^{\prime\prime}(A)\sumt \int_{(i-1)\Del}^{i\Del} (\bbe[|\beta^{n,i,\eps}_x|^{\frac p2}] + \bbe[|\beta^{n,\eps}_x(r)|^{\frac p 2}]) \,\dd r,
	\end{align*}
	where $\Phi(A)$, $\Phi'_A(\delta)$, and $\Phi''(A)$ have the same properties as in Lemma~\ref{sibounded} (but may take different values, of course).
	As a consequence of the elementary inequality 
	\beq\label{elem} |a_1a_2-a_3a_4|\leq (|a_1|+\cdots+|a_4|)(|a_1-a_3|+|a_2-a_4|),\eeq which holds for all $a_1,\ldots,a_4\in\R$, we deduce using the Cauchy--Schwarz inequality,
	\beq\label{expr-3} \begin{split}
		&\bbe\left[\left|(\beta^{n,i,\eps}_x)_{j_1k_1,j_2k_2}-(\beta^{n,\eps}_x(r)))_{j_1k_1,j_2k_2}\right|\right]\\
		&\quad\lec \sup_{(s,y)\in\R\times\R^d}  \bbe[\si^2(s,y)]^{\frac 12}     \\
		&\quad\quad\times \iiint\left(\bbe[|\si(i\del-\eps,x_{j_1}-y)-\si(r-\eps,x_{j_1}-y)|^2]^{\frac 12}\right.\\
		&\qquad \quad+\left.\bbe[|\si(i\del-\eps,x_{j_1}-z)-\si(r-\eps,x_{j_1}-z)|^2]^{\frac 12}\right)\\
		&\qquad\times\bone_{s<\eps+(k_1-1)\del}|\,\Pi^n_{k_2-k_1,x_{j_2}-x_{j_1}}|(\dd s,\dd y,\dd z) \\
		&\quad\lec \sup\left\{\bbe[|\si(t,y)-\si(s,y)|^2]^{\frac12}\colon |t-s|\leq \del,~y\in\R^d\right\}
	\end{split}\eeq 
	for all $r\in[(i-1)\del,i\del)$. Moreover, by the BDG inequality and \eqref{cross-bound-3}, we also have that $\bbe[|\beta^{n,i,\eps}_x|^{p/2}]$ and $\bbe[|\beta^{n,\eps}_x(r)|^{p /2}]$ are uniformly bounded in $n$, $i$, $\eps$, and $r$. So altogether, we have 
	\begin{align*}
	\bbe[|A^{n,\eps,1}_x(t)|]&\lec \Phi'_A(\delta)t+\frac{\Phi(A)w(\del)t}{\delta}+\Phi^{\prime\prime}(A)t,
	\end{align*}
	which vanishes if we let $n\to\infty$, $\delta\to0$, and $A\to\infty$.
	
	Next, since $\si$ is bounded, observe from Lemma~\ref{delta} that $(\beta^{n,\eps}_x(r)))_{j_1k_1,j_2k_2}$ converges to $0$ in $L^1$ as $n\to\infty$ whenever $j_1\neq j_2$. If $j_1=j_2 = j$ and $k_1\leq k_2$, we  consider the decomposition
	\[ (\beta^{n,\eps}_x(r)))_{jk_1,jk_2}-\Ga_{k_2-k_1}\si^2(r-\eps,x_j) = B^{n,\eps,R,1}_x + B^{n,\eps,R,2}_x- B^{n,\eps,3}_x,\]
	where $\Ga_r$ are the numbers from \eqref{Ga-formula}, and where we define for $R>0$, 
	\begin{align*} B^{n,\eps,R,1}_x&= \iiint_0^\infty \lp\si(r-\eps,x_{j_1}-y)\si(r-\eps,x_{j_1}-z)-\si^2(r-\eps,x_{j_1})\rp\\
	&\quad\times \bone_{s<\eps+(k_1-1)\del}\bone_{|y|+|z|< R}\, \Pi^n_{k_2-k_1,0}(\dd s,\dd y,\dd z),\\
	B^{n,\eps,R,2}_x &=
	\iiint_0^\infty \lp\si(r-\eps,x_j-y)\si(r-\eps,x_j-z)-\si^2(r-\eps,x_j)\rp\\
	&\quad\times \bone_{s<\eps+(k_1-1)\del} \bone_{|y|+|z|\geq R}\,\Pi^n_{k_2-k_1,0}(\dd s,\dd y,\dd z),\\
	B^{n,\eps,3}_x&=   \si^2(r-\eps,x_j)(\Ga_{k_2-k_1}-\Pi^n_{k_2-k_1,0}([0,\eps+(k_1-1)\del)\times\R^d\times\R^d)).\end{align*}
	With similar arguments as in \eqref{expr-3}, we can bound
	\begin{align*} \bbe[|B^{n,\eps,R,1}_x|]&\lec \iiint_0^{\infty} (w(|y|)+w(|z|)) \bone_{|y|+|z|< R}  \,|\Pi^n_{k_2-k_1,0}|(\dd s,\dd y,\dd z)\\ &\lec w(R).\end{align*}
	Moreover, since $\si$ is bounded, we have
	\begin{align*} \bbe[|B^{n,\eps,R,2}_x|]&\lec |\Pi^n_{k_2-k_1,0}|(\R\times\{(y,z)\colon |y|+|z|>R\}),\\
	\bbe[|B^{n,\eps,3}_x|]&\lec \Ga_{k_2-k_1}-\Pi^n_{k_2-k_1,0}([0,\eps+(k_1-1)\del)\times\R^d\times\R^d).\end{align*}
	By Lemma~\ref{delta}, these two terms vanish as $n\to\infty$. Since $R$ was arbitrary, Assumption~\ref{AssLLN}3 implies that $\beta^{n,\eps}_x(r) \to \beta^{\eps}_x(r)$ in $L^1$ where $(\beta^\eps_x(r))_{j_1k_1,j_2k_2} = \Ga_{|k_2-k_1|}\si^2(r-\eps,x_j)\bone_{j_1=j_2=j}$. By the continuity of $\un\mu_f$, the boundedness of $\si$, and dominated convergence, it follows that $\un\mu_f(\beta^{n,\eps}_x(r))\to \un\mu_f(\beta^{\eps}_x(r))$ in $L^1$. Upon realizing that $\un\mu_f(\beta^\eps_x(r))=\mu_f(\un\si^2(r-\eps,x))$ by definition, the dominated convergence theorem shows that $A^{n,\eps,2}_x(t)\to0$ in $L^1$ if $n\to\infty$ and $\eps\to0$. 
	Using Assumption~\ref{AssLLN}3, we can apply similar bounds as in \eqref{expr-3} to show that
	$\lim_{\eps\to0}\limsup_{n\to\infty} \bbe[|A^{n,\eps,3}_x(t)|]=0$.
	And
	finally, since $\si$ is bounded, we also have $\bbe[|A^{n,4}_x(t)|]\lec t-\del t^*_n\leq L\del\to0$. 
	\epr
	
	\section{Details for the proof of Theorem~\ref{CLT}}\label{details}
	
	In the sequel, we will frequently use a \emph{(standard) size estimate} or \emph{argument} to determine the asymptotic behavior 
	of expressions of the following or similar form [recall \eqref{abbr}]:
	\beq\label{u-expr} \begin{split} U_n(t) &= \delh\sumt h(z^n_i)\\
		&\quad\times\iint \frac{\Delta^n_i G_{x_0-y}(s)}{\tau_n} (\si(s,y)-\si(t^{n,i}_{\ell_n},y))\bone^{n,i}_{\ell_n,\ell'_n}(s)\,W(\dd s,\dd y).\end{split} \eeq
	Here, the point $x_0\in\R^d$ is arbitrary, $\ell^{(\prime)}_n = [\del^{-l^{(\prime)}}]$ with some $-\infty\leq l'<l\leq\infty$, $h$ is a function with  $|h(z)|\lec 1+|z|^{p-1}$ for some $p\geq2$, and $z^n_i$ are random variables such that 
	\beq\label{size} \sup_{n\in \N} \sup_{i=1,\ldots,T^*_n} \bbe[|z^n_i|^p]<\infty. \eeq
	In most cases, the variables $z^n_i$ are normalized increments, possibly truncated or with modified $\si$ [e.g., $z^n_i=\iint \frac{\Delta^n_i G_{x_0-y}(s)}{\tau_n} \si(s,y)\,W(\dd s,\dd y)$ or $z^n_i=\iint \frac{\Delta^n_i G_{x_0-y}(s)}{\tau_n} \si(t^{n,i}_{\ell_n},y)\bone^{n,i}_{\ell_n,\ell'_n}(s)\,W(\dd s,\dd y)$],
	or combinations thereof.
	
	In this setting,  Hölder's inequality with exponents $\frac{p}{p-1}$ and $p$ implies
	\begin{align*}
	&\bbe\lb \supt |U_n(t)|\rb\\
	&\quad \leq \delh\sumT \bbe[|h(z^n_i)|^{\frac p {p-1}}]^{\frac{p-1}{p}} \\
	&\qquad\times\bbe\lb \lv\iint \frac{\Delta^n_i G_{x_0-y}(s)}{\tau_n} (\si(s,y)-\si(t^{n,i}_{\ell_n},y))\bone^{n,i}_{\ell_n,\ell'_n}(s)\,W(\dd s,\dd y)\rv^p\rb^{\frac1p}.
	\end{align*}
	By \eqref{size} and the growth assumptions on $h$, we have $\bbe[|h(z^n_i)|^{p/(p-1)}]\lec 1$, uniformly in $i$ and $n$. Thus, by the BDG inequality, Minkowski's integral inequality, and the Cauchy--Schwarz inequality, as well as Lemma~\ref{mom-ex} and Lemma~\ref{nu}, the left-hand side of the previous display is  
	\begin{align*}
	&\lec \delh\sumT\bbe\Bigg[\Bigg( \iiint \frac{\Delta^n_i G_{x_0-y}(s) \Delta^n_i G_{x_0-z}(s)}{\tau_n^2} (\si(s,y)-\si(t^{n,i}_{\ell_n},y))\\
	&\quad\times (\si(s,z)-\si(t^{n,i}_{\ell_n},z))\bone^{n,i}_{\ell_n,\ell'_n}(s)\,\dd s\,\La(\dd y,\dd z)\Bigg)^{\frac p2}\Bigg]^{\frac1p}\\
	&\leq \delh\sumT\Bigg( \iiint \frac{|\Delta^n_i G_{x_0-y}(s) \Delta^n_i G_{x_0-z}(s)|}{\tau_n^2} \bbe\Big[|\si(s,y)-\si(t^{n,i}_{\ell_n},y)|^{\frac p2}\\
	&\quad\times |\si(s,z)-\si(t^{n,i}_{\ell_n},z)|^{\frac p 2}\Big]^{\frac 2p}\bone^{n,i}_{\ell_n,\ell'_n}(s)\,\dd s\,\La(\dd y,\dd z)\Bigg)^{\frac12}\\
	&\leq \delh\sumT\Bigg( \iiint \frac{|\Delta^n_i G_{x_0-y}(s) \Delta^n_i G_{x_0-z}(s)|}{\tau_n^2} \bbe[|\si(s,y)-\si(t^{n,i}_{\ell_n},y)|^p]^{\frac1p}\\
	&\quad\times \bbe[|\si(s,z)-\si(t^{n,i}_{\ell_n},z)|^p]^{\frac1p}\bone^{n,i}_{\ell_n,\ell'_n}(s)\,\dd s\,\La(\dd y,\dd z)\Bigg)^{\frac12}\\
	&\lec  \delh[T/\del] (\ell_n\del)^{\frac 12}|\Pi^n_{0,0}|((\ell'_n\del,\ell_n\del)\times\R^d\times\R^d)^{\frac12}\\
	&\lec \del^{-\frac12} (\ell_n\del)^{\frac12}\del^{\frac\nu 2 \ell'}.
	\end{align*}
	The individual factors in this final bound can be attributed to the components in \eqref{u-expr}. The part $\Del^{1/2}\sumt$ contributes $\del^{-1/2}$, the variables $h(z^n_i)$ have bounded moments and hence do not contribute, the difference $\si(s,y)-\si(t^{n,i}_{\ell_n,y})$ contributes $(\ell_n\del)^{1/2}$, and the stochastic integral of $\bone^{n,i}_{\ell_n,\ell'_n}$ contributes $\del^{\nu  \ell'/2}$. 
	Therefore, in the following proofs, if we encounter a term like $U_n(t)$, we will \emph{directly} conclude from \eqref{u-expr} that $$\bbe\lb\supt |U_n(t)|\rb\lec \del^{-\frac12} (\ell_n\del)^{\frac12}\del^{\frac\nu 2 \ell'},$$ without going through similar arguments again.
	
	We will also apply this type of estimates to expressions that are more complicated than \eqref{u-expr}, for example, if the stochastic integral in \eqref{u-expr} is squared, or if it is replaced by the product of two stochastic integrals and the function $h$ satisfies $|h(z)|\lec 1+|z|^{p-2}$. Then the generalized Hölder's inequality with exponent $\frac{p}{p-2}$ for $h(z^n_i)$ and exponent $p$ for each appearing stochastic integral can be used to factorize $\bbe[\supt |U_n(t)|]$ in the same manner as before. The key observation is that the total size of such expressions can  always be  determined  component by component, which, of course,  applies to vector- or matrix-valued versions of \eqref{u-expr} as well.
	
	In all proofs below, except for the proof of \eqref{Cmt} for Proposition~\ref{CLT-core}, we may and will assume that $M=1$.

	\bpr[Proof of Lemma~\ref{int-trunc}] Since $|f(z+z')-f(z)| \lec (1+|z|^{p-1}+|z'|^{p-1})|z'|$ by Assumption~\ref{AssCLT}1, a standard size argument implies that 
	\begin{align*}
	&\Delh \sumT \bbe\left[\left|f\bigg( \frac{\un\Delta^n_i Y_x}{\tau_n}\bigg) - f\bigg( \frac{\ga^{n,i,0}_x}{\tau_n}\bigg)\right|\right]\\
	& \quad\lec \Delh \sumT \bbe\Bigg[\lp 1+\bv \frac{\un\Delta^n_i Y_x}{\tau_n}\bv^{p-1} + \bv \frac{\ga^{n,i,0}_x}{\tau_n}\bv^{p-1}\rp \\
	&\qquad\times \lv \iint \frac{ \un \Delta^n_i G_{x,y}(s)}{\tau_n}\si(s,y)\bone_{s\leq (i-\la^0_n)\del}\,W(\dd s,\dd y) \rv\Bigg]\\
	& \quad\lec \Del^{-\frac12+\frac{\nu}{2}a_0}, 
	\end{align*}
	which converges to $0$ as $n\to\infty$ because $a_0>\frac 1 \nu$.	
	\epr

	\bpr[Proof of Lemma~\ref{knknprime}] Fix $r=1,\ldots,R$. We rearrange the sum of the $\ov \delta^{n,r}_i$-terms as
	\[ \sumt \ov\delta^{n,r}_i = \sum_{k=1}^{\la^{r-1}_n+L-1} C^{n}_k(t)   \]
	with
	\[ C^{n}_k(t)=\sum_{j=1}^{N^{n,r}_k(t)} \ov\delta^{n,r}_{k+(j-1)(\la^{r-1}_n+L-1)}\qquad\text{and}\qquad N^{n,r}_k(t)=\left[\frac{t^*_n-k}{\la^{r-1}_n+L-1}\right]. \]
	For a fixed value of $k=1,\ldots,\la^{r-1}_n+L-1$, the term $C^{n}_k(t)$ sums up only every $(\la^{r-1}_n+L-1)$th member of the sequence $(\ov\delta^{n,r}_i\colon i=1,\ldots, t^*_n)$, starting with $\ov\delta^{n,r}_k$. The important observation is now that $\ov\delta^{n,r}_{k+(j-1)(\la^{r-1}_n+L-1)}$ is $\calf^n_{k+(j-1)(\la^{r-1}_n+L-1)+L-1}$-measurable with vanishing conditional expectation given $\calf^n_{k+(j-1)(\la^{r-1}_n+L-1)-\la_n^{r-1}}$. 
	
	Hence, for each $k$, $C^{n}_k(t)$ is a sum of martingale differences with respect to the discrete-time filtration $(\calf^n_{k+(j-1)\la_n^{r-1}+j(L-1)}\colon j=0,\ldots,N^{n,r}_k(t))$, whose quadratic variation process is given by $$\sum_{j=1}^{N^{n,r}_k(t)} (\ov\delta^{n,r}_{k+(j-1)(\la^{r-1}_n+L-1)})^2.$$  By Assumption~\ref{AssCLT}1, we have $|f(z+z')-f(z)| \lec (1+|z|^{p-1}+|z'|^{p-1})|z'|$, so using the BDG inequality and a standard size estimate, we obtain
	\beq\label{BDG-standard}\begin{split} \bbe\lb \sup_{t\in[0,T]} |C^{n}_k(t)|^2\rb &\lec \sum_{j=1}^{N^{n,r}_k(T)} \bbe[(\ov\delta^{n,r}_{k+(j-1)(\la^{r-1}_n+L-1)})^2]\\
		& \lec \sum_{j=1}^{N^{n,r}_k(T)} \bbe[(\delta^{n,r}_{k+(j-1)(\la^{r-1}_n+L-1)})^2]\\
		&\lec \Delta_n N^{n,r}_k(T) |\Pi^n_{0,0}|(\{\la_n^r\Delta_n<s<\la_n^{r-1} \Delta_n\})\\
		&\lec \frac{\Delta_n\del^{\nu a_r} }{\del \la^{r-1}_n} =  \frac{\del^{\nu a_r}}{\la^{r-1}_n}. \end{split}\eeq
	Consequently,
	\begin{align*} \bbe\lb \sup_{t\in[0,T]} \lv \sumt \ov\delta^{n,r}_i\rv\rb &\leq \sum_{k=1}^{\la_n^{r-1}+L-1} \bbe\lb \sup_{t\in[0,T]} |C^{n}_k(t)|^2\rb^{\frac 1 2} \lec \sqrt{\la_n^{r-1}}  \del^{\frac{\nu}{2} a_r} \\ &= \del^{\frac{\nu a_r-a_{r-1}}{2}} \to 0\end{align*}
	as $n\to\infty$ since $a_r>\frac{a_{r-1}}{\nu}$.
	\epr

	\bpr[Proof of Lemma~\ref{condex-trunc}] Recalling the notation from \eqref{abbr}, we define the variables
	\beq\label{eta}\begin{split}
		\zeta^n_i&=\frac{\ga^{n,i,r-1}_x-\ga^{n,i,r}_x}{\tau_n}= \iint \frac{\un\Delta^n_i G_{x,y}(s)}{\tau_n}\bone^{n,i}_{\la_n^{r-1},\la_n^r}(s)\si(s,y)\,W(\dd s,\dd y),\\
		\eta^{n}_i &= \iint \frac{\un\Delta^n_i G_{x,y}(s)}{\tau_n}\si\left(t^{n,i}_{\la_n^r},y\right)\bone_{s>t^{n,i}_{\la_n^r}}\,W(\dd s,\dd y). \end{split}\eeq
	Then by Taylor's theorem and Assumption~\ref{AssCLT}1, we have
	\beq\label{etani} \begin{split} \delta^{n,r}_i &=\delh \sum_{\al} \frac{\partial}{\partial z_\al} f\lp\frac{\ga^{n,i,r}_x}{\tau_n} \rp (\zeta^n_i)_\al + \frac{1}{2}\delh \sum_{\al,\beta} \frac{\partial^2}{\partial z_\al\,\partial z_\beta} f(\eta^{n,1}_i)(\zeta^n_i)_\al(\zeta^n_i)_\beta\\
		&= \delh \sum_{\al} \frac{\partial}{\partial z_\al} f(\eta^n_i ) (\zeta^n_i)_\al+
		\delh \sum_{\al,\beta} \frac{\partial^2}{\partial z_\al\,\partial z_\beta} f(\eta^{n,2}_i)(\zeta^n_i )_\al\lp \frac{\ga^{n,i,r}_x}{\tau_n}- \eta^n_i\rp_\beta \\
		&\quad+\frac12\delh \sum_{\al,\beta} \frac{\partial^2}{\partial z_\al\, \partial z_\beta} f(\eta^{n,1}_i)(\zeta^n_i)_\al(\zeta^n_i)_\beta,\end{split}\eeq
	where $\al$ and $\beta$ run through the set of indices in $\{1,\ldots, N\}\times\{1,\ldots,L\}$, and where $\eta^{n,1}_i$ (resp., $\eta^{n,2}_i$) is some point on the line between $\frac{\ga^{n,i,r-1}_x}{\tau_n}$ and $\frac{\ga^{n,i,r}_x}{\tau_n}$ (resp., between $\frac{\ga^{n,i,r}_x}{\tau_n}$ and $\eta^n_i$). Let us denote the last three terms on the right-hand side of \eqref{etani} by $\delta^{n,r,1}_i$, $\delta^{n,r,2}_i$ and $\delta^{n,r,3}_i$.
	Then
	\[
	\sumt \bbe[\delta^{n,r} _i\,|\, \calf^n_{i-\la_n^{r-1}}] = \sum_{j=1}^3 D^n_j(t),\qquad D^n_j(t)=\sumt \bbe[\delta^{n,r,j}_i\,|\, \calf^n_{i-\la_n^{r-1}}].\]

	Since $|\frac{\partial^2}{\partial x_\al\,\partial x_\beta} f(z)|\lec 1+|z|^{p-2}$ by Assumption~\ref{AssCLT}1 and
	\beq\label{diff}
	\frac{\ga^{n,i,r}_x}{\tau_n}- \eta^n_i= 
	\iint \frac{\un\Delta^n_i G_{x,y}(s)}{\tau_n}\left(\si(s,y)-\si\left(t^{n,i}_{\la_n^r},y\right)\right)\bone_{s>t^{n,i}_{\la_n^r}}\,W(\dd s,\dd y),
	\eeq
	a standard size estimate and the fact that $\nu>1$ (see Lemma~\ref{nu}) show that
	\beq\label{Cn1}  \bbe\lb\sup_{t\in[0,T]} |D^n_2(t)|\rb \lec \Delta_n^{-\frac12}   \Delta_n^{\frac{\nu}{2} a_r}  (\la_n^r\Delta_n)^{\frac12}=  \Delta_n^{(\frac \nu 2 -\frac12)a_r}\to0 \eeq
	as $n\to\infty$. Similarly, we have
	$\bbe\lb\sup_{t\in[0,T]} |D^n_3(t)|\rb \lec \del^{-\frac12}\del^{\nu a_r}\to0$ because $a_r>\frac{1}{2\nu}$.
	Finally, concerning $D^n_1(t)$, notice that
	\begin{align*}
	\bbe\left[ \frac{\partial}{\partial z_\al} f(\eta^{n}_i)(\zeta^n_i)_\al \,\Big|\,\calf^n_{i-\la_n^{r-1}}\right]&= \bbe\left[ \bbe\left[ \frac{\partial}{\partial z_\al}  f(\eta^{n}_i )(\zeta^n_i)_\al\,\Big|\,\calf^n_{i-\la_n^{r}}\right] \,\Big|\, \calf^n_{i-\la_n^{r-1}}\right] \\
	&= \bbe\left[(\zeta^n_i)_\al  \bbe\left[ \frac{\partial}{\partial z_\al} f(\eta^{n}_i )\,\Big|\,\calf^n_{i-\la_n^{r}}\right] \,\Big|\, \calf^n_{i-\la_n^{r-1}}\right]\\
	&=0
	\end{align*}
	because $\eta^{n}_i$ has a centered normal distribution conditionally on $\calf^n_{i-\la_n^{r}}$ and $\frac{\partial}{\partial z_\al} f$ is an odd function as the derivative of an even function. This implies that $D^n_1(t)= 0$ for all $t\in[0,T]$.
	\epr

	\bpr[Proof of Lemma~\ref{approx}] Since $V^{n,m,3}(t)$ contains at most $(m+1)\la_n+L-1$ terms, it follows from a size estimate that
	\[ \bbe\lb \supt | V^{n,m,3}(t)|\rb \lec (m+1)\la_n\delh\leq (m+1)\del^{\frac12-a}\to0 \]
	if $a$ is close  to $\frac 1{2\nu}$.
	Next, we recall \eqref{hatpsi}
	and approximate $V^{n,m,2}(t)$ by 
	\beq\label{Vhat-2} \wh V^{n,m,2}(t) =\sum_{j=1}^{J^{n,m}(t)}\sum_{k=1}^{\la_n+L-1} \wh\psi^n_{(j-1)((m+1)\la_n+L-1)+m\la_n+k,k}.\eeq
	As in the proof of Lemma~\ref{knknprime}, $V^{n,m,2}(t)-\wh V^{n,m,2}(t) = \sum_{k=1}^{\la_n+L-1} E^{n,m}_k(t)$, where $E^{n,m}_k(t)$ is a martingale sum of the form
	\[ E^{n,m}_k(t)=\!\sum_{j=1}^{J^{n,m}(t)} [\psi^n_{(j-1)((m+1)\la_n+L-1)+m\la_n+k}-\wh\psi^n_{(j-1)((m+1)\la_n+L-1)+m\la_n+k,k}]\]
	relative to the filtration $(\calf^n_{(j-1)((m+1)\la_n+L-1)+m\la_n+k}\colon j=0,\ldots,J^{n,m}(t))$. Moreover, since $\bbe[ f( \wh\xi^{n}_{i,k}) \,|\, \calf^n_{i-\la_n-k} ]=\bbe[ f( \wh\xi^{n}_{i,k}) \,|\, \calf^n_{i-\la_n} ]$, we obtain for arbitrary $i$ and $k$,
	\begin{align*}
	\bbe[|\psi^n_i -\wh\psi^n_{i,k}|^2] 
	&\lec \del \bbe\lb \lv f\lp \frac{\ga_x^{n,i}}{\tau_n} \rp - f\lp \wh\xi^{n}_{i,k}\rp \rv^2\rb \lec \del((\la_n\del)^{\frac12})^2= \del^2\la_n,
	\end{align*}
	where the penultimate step follows from a  standard size estimate together with Lemma~\ref{mom-ex}. Hence, as $n\to\infty$,
	\beq\label{approx-2}\begin{split} \bbe\lb\sup_{t\in[0,T]} |V^{n,m,2}(t)-\wh V^{n,m,2}(t)|\rb &\lec \la_n(J^{n,m}(T)\del^2\la_n)^{\frac12}\\
		&\lec m^{-\frac12}\Del^{\frac12- a}\to0. \end{split}\eeq
	
	It remains to prove that $\wh V^{n,m,2}$ is asymptotically negligible. Observe again that the sum over $j$ in \eqref{Vhat-2} has a martingale structure. Moreover, conditionally on $\calf^n_{(j-1)((m+1)\la_n+L-1)+m\la_n}$, the variables $\wh\xi^n_{(j-1)((m+1)\la_n+L-1)+m\la_n+k,k}$ and $\wh\xi^n_{(j-1)((m+1)\la_n+L-1)+m\la_n+l,l}$ are normally distributed. Thus,
	\beq\label{Vhat-3}\begin{split}
		\bbe\lb\sup_{t\in[0,T]} |\wh V^{n,m,2}(t)|^2\rb &\lec \sum_{j=1}^{J^{n,m}(T)}\bbe\lb \lv \sum_{k=1}^{\la_n+L-1} \wh\psi^n_{(j-1)((m+1)\la_n+L-1)+m\la_n+k,k}\rv^2\rb\\
		&=\sum_{j=1}^{J^{n,m}(T)}\sum_{k,l=1}^{\la_n+L-1} \bbe[\wh\psi^n_{(j-1)((m+1)\la_n+L-1)+m\la_n+k,k}\\
		&\quad\times\wh\psi^n_{(j-1)((m+1)\la_n+L-1)+m\la_n+l,l}]\\
		&=\del\sum_{j=1}^{J^{n,m}(T)}\sum_{k,l=1}^{\la_n+L-1} \bbe[ \un\rho_f(v^{n,m,j}_k,v^{n,m,j}_l, c^{n,m,j}_{k,l}) ],
	\end{split}\raisetag{-5\baselineskip}\eeq
	where for $v^{(1|2)} = (v^{(1|2)}_{pq,p'q'})_{p,p',q,q'=1}^{N,N,L,L}$ and $c = (c_{pq,p'q'})_{p,p',q,q'=1}^{N,N,L,L}$, we define
	$$\un\rho_f(v^{(1)},v^{(2)},c)=\cov(f(\un Z^{(1)}), f(\un Z^{(2)})),$$ 
	and where $\un Z^{(1)},\un Z^{(2)}\in\R^{N\times L}$ are jointly Gaussian with zero mean, covariances $\cov(\un Z^{(1|2)}_{pq}, \un Z^{(1|2)}_{p'q'}) = v^{(1|2)}_{pq,p'q'}$, and cross-covariances $\cov(\un Z^{(1)}_{pq}, \un Z^{(2)}_{p'q'}) = c_{pq,p'q'}$. If $v^{(1)}=v^{(2)}=v$, we use the notation $\un\rho_f(v,c)$.
	Furthermore, writing $i(j-1,n,m)$ for $(j-1)((m+1)\la_n+L-1)+m\la_n$, we have by a change of variables, for $q'\geq q$,
	\beq\label{variances}\begin{split}
		\quad&(v^{n,m,j}_k)_{pq,p'q'}\\
		&\quad = \iiint_0^\infty \frac{\Delta^n_{i(j-1,n,m)+k+q-1} G_{x_p-y}(s)\Delta^n_{i(j-1,n,m)+k+q'-1} G_{x_{p'}-z}(s)}{\tau^2_n}\\
		&\qquad\times\si\lp t^{n,i(j-1,n,m)}_{\la_n},y\rp\si\lp t^{n,i(j-1,n,m)}_{\la_n},z\rp \bone_{s>t^{n,i(j-1,n,m)+k}_{\la_n}} \,\dd s\,\La(\dd y,\dd z)\\
		&\quad=\iiint_{0}^{(\la_n+q-1)\del} \frac{G_y(s)-G_y(s-\del)}{\tau_n} \\
		&\qquad\times\frac{G_{z+x_{p'}-x_p}(s+(q'-q)\del)-G_{z+x_{p'}-x_p}(s+(q'-q-1)\del)}{\tau_n}\\
		&\qquad \times \si\lp t^{n,i(j-1,n,m)}_{\la_n},x_p-y\rp\si\lp t^{n,i(j-1,n,m)}_{\la_n},x_p-z\rp\,\dd s\,\La(\dd y,\dd z)\\
		& \quad=\iiint_{0}^{(\la_n+q-1)\del} \si\lp t^{n,i(j-1,n,m)}_{\la_n},x_p-y\rp\\
		&\qquad\times\si\lp t^{n,i(j-1,n,m)}_{\la_n},x_p-z\rp\,\Pi^n_{q'-q,x_{p'}-x_p}(\dd s,\dd y,\dd z).
	\end{split}\raisetag{-8\baselineskip}\eeq
	
	Similarly, we have for $k\geq l$ the following two cases: if $q'\geq q+(k-l)$, then
	\beq\label{covariances}\begin{split}
		&(c^{n,m,j}_{k,l})_{pq,p'q'}\\
		&\quad= \iiint_0^\infty \frac{\Delta^n_{i(j-1,n,m)+k+q-1} G_{x_p-y}(s)\Delta^n_{i(j-1,n,m)+l+q'-1} G_{x_{p'}-z}(s)}{\tau^2_n} \\
		&\qquad\times \si\lp t^{n,i(j-1,n,m)}_{\la_n},y\rp \si\lp t^{n,i(j-1,n,m)}_{\la_n},z\rp\\
		&\qquad\times\bone_{((i(j-1,n,m)+k-\la_n)\del,(i(j-1,n,m)+k+q-1)\del)}(s)\,\dd s\,\La(\dd y,\dd z)\\
		&\quad=\iiint_0^{(\la_n+q-1)\del}  \si\lp t^{n,i(j-1,n,m)}_{\la_n},x_p-y\rp\\
		&\qquad\times\si\lp t^{n,i(j-1,n,m)}_{\la_n},x_p-z\rp \,\Pi^n_{q'-q-(k-l),x_{p'}-x_p}(\dd s,\dd y,\dd z);
	\end{split}\raisetag{-5\baselineskip}\eeq
	and if $q'\leq q+(k-l)$, then
	\beq\label{covariances-2}\begin{split}
		&(c^{n,m,j}_{k,l})_{pq,p'q'}\\
		&\quad= \iiint_0^{(\la_n+q'-(k-l)-1)\del}  \si\lp t^{n,i(j-1,n,m)}_{\la_n},x_{p'}-y\rp\\
		&\qquad\times\si\lp t^{n,i(j-1,n,m)}_{\la_n},x_{p'}-z\rp\,\Pi^n_{q-q'+(k-l),x_p-x_{p'}}(\dd s,\dd y,\dd z).
	\end{split}\eeq
	Note that $v^{n,m,j}_k$ does not depend on $k$, so we write $v^{n,m,j}= v^{n,m,j}_k=v^{n,m,j}_l$ in the following.
	
	From \eqref{variances},  Lemma~\ref{mom-ex}, and \eqref{cross-bound-3}, we obtain $\bbe[|v^{n,m,j}_k|^{p/ 2}]\lec 1$. Furthermore, $c^{n,m,j}_{k,l}$ only depends on the difference $|k-l|$, which we shall therefore denote by $c^{n,m,j}_{|k-l|}$. And since $q,q'\leq L$, the case $q'\geq q +|k-l|$ analyzed in \eqref{covariances} can only occur for small values of $|k-l|$, while in \eqref{covariances-2}, we have $q-q'+|k-l|\geq (|k-l|-L)\vee 0$. Consequently, \eqref{covariances}, \eqref{covariances-2}, and \eqref{cross-bound-3} imply that $\bbe[|c^{n,m,j}_{|k-l|}|^{p}]^{1/ p}\lec \ov \Ga'_{|k-l|}$ for all $k$ and $l$, where 
	\beq\label{Gaprime}\ov \Ga'_{|k-l|} = \ov \Ga_{(|k-l|-L)\vee0}.\eeq 
	
	Next, observe that if we view $\un Z^{(1|2)}$ as an $(NL)$-dimensional vector, we can apply Lemma~\ref{cov-add} to  $\un\rho_f(v^{n,m,j},c^{n,m,j}_{|k-l|})$. Together with Assumption~\ref{AssCLT}1, Hölder's inequality, and Lemma~\ref{mom-ex}, we have from \eqref{Vhat-3},
	\beq\label{Vnm2}\begin{split}
		\bbe\lb\sup_{t\in[0,T]} |\wh V^{n,m,2}(t)|^2\rb &\lec\del
		\sum_{j=1}^{J^{n,m}(T)}\sum_{k,l=1}^{\la_n+L-1} \lp 1+ \sum_{p'=1}^N\sum_{q'=1}^L \bbe[|v^{n,m,j}_{p'q',p'q'}|^p]^{\frac{p-2}{p}} \rp \\
		&\quad\times\sum_{p',p''=1}^N \sum_{q',q''=1}^L \bbe[|(c^{n,m,j}_{|k-l|})_{p'q',p''q''}|^p]^{\frac 2 p}\\
		&\lec \del \sum_{j=1}^{J^{n,m}(T)}\sum_{k,l=1}^{\la_n+L-1} (\ov\Ga'_{|k-l|})^2\\
		&\lec  \del J^{n,m}(T)(\la_n+L-1)\sum_{r=0}^\infty (\ov\Ga'_r)^2.
	\end{split}\raisetag{-5\baselineskip}\eeq
	Since $\ov\Ga_r$---and hence $\ov\Ga'_r$---is square-summable by Lemma~\ref{Ga}, the last term is bounded by a multiple of $\del J^{n,m}(T)\la_n\lec m^{-1}$, which converges to $0$ as $m\to\infty$, uniformly in $n\in\N$. 
	\epr
	
	\bpr[Proof of Lemma~\ref{Vhat}] The proof is very similar to \eqref{approx-2}, so we omit the details.\epr 
	
	\bpr[Proof of Proposition~\ref{CLT-core}] It remains to prove properties (1)--(3) described after the statement of Proposition~\ref{CLT-core}.
	Note also that $t\in[0,T]$ is fixed, so for simplicity, we suppress $t$ in the variables appearing below. Furthermore, we use the abbreviation $i[j,n,m]=j((m+1)\la_n+L-1)$ [note the difference to $i(j,m,n)$ as defined in the proof of Lemma~\ref{approx}].  
	
	\vspace{0.5\baselineskip}
	\noindent{\itt{Property (1): }} Since the proof for $m_1\neq m_2$ is completely analogous, we only consider the case $m_1=m_2$ in \eqref{Cmt}, where we may assume without loss of generality that $M=1$. Then we have
	\begin{align*}
	\calc^{n,m} &= \sum_{j=1}^{J^{n,m}} \bbe\lb \lp\sum_{k=1}^{m\la_n} \wh \psi^{n}_{i[j-1,n,m]+k,k}\right)^2 \ \Bigg\vert\ \calf^n_{i[j-1,n,m]-\la_n}\rb\\
	&=\del \sum_{j=1}^{J^{n,m}}\sum_{k,l=1}^{m\la_n} \left\{\bbe\lb f(\wh\xi^{n}_{i[j-1,n,m]+k,k})f(\wh\xi^{n}_{i[j-1,n,m]+l,l}) \,|\, \calf^n_{i[j-1,n,m]-\la_n} \rb\right.\\
	&\quad -\bbe\lb f(\wh\xi^{n}_{i[j-1,n,m]+k,k}) \,|\, \calf^n_{i[j-1,n,m]-\la_n} \rb\\
	&\qquad\times \left.\bbe\lb f(\wh\xi^{n}_{i[j-1,n,m]+l,l}) \,|\, \calf^n_{i[j-1,n,m]-\la_n} \rb\right\}.
	\end{align*}
	
	Conditionally on $\calf^n_{i[j-1,n,m]-\la_n}$, each $\wh\xi$ variable has a centered normal distribution, and two of them, one indexed by $k$ and one indexed by $l$, are conditionally independent as soon as $|k-l|\geq\la_n+L-1$. Thus, many terms above cancel, and $C^{n,m}$ becomes
	\beq\label{Cnm}
	\begin{split}
		&\del \sum_{j=1}^{J^{n,m}}\sum_{|k-l|<\la_n+L-1} \left\{\bbe\lb f(\wh\xi^{n}_{i[j-1,n,m]+k,k})f(\wh\xi^{n}_{i[j-1,n,m]+l,l}) \,|\, \calf^n_{i[j-1,n,m]-\la_n} \rb\right.\\
		&\qquad - \bbe\lb f(\wh\xi^{n}_{i[j-1,n,m]+k,k}) \,|\, \calf^n_{i[j-1,n,m]-\la_n} \rb\\
		&\quad\qquad\times\left.\bbe\lb f(\wh\xi^{n}_{i[j-1,n,m]+l,l}) \,|\, \calf^n_{i[j-1,n,m]-\la_n} \rb\right\}\\
		&\quad=\del \sum_{j=1}^{J^{n,m}}\sum_{|k-l|<\la_n+L-1} \un\rho_{f}(\wh v^{n,m,j}, \wh c^{n,m,j}_{|k-l|}).
	\end{split}\raisetag{-4.5\baselineskip}\eeq
	This has a similar structure as \eqref{Vhat-3}, and $\wh v^{n,m,j}$ and $\wh c^{n,m,j}_{|k-l|}$ are defined as in \eqref{variances} and \eqref{covariances} but with $i[j-1,n,m]$ instead of $i(j-1,n,m)$.
	Hence, $\calc^{n,m}=\calc^{n,m}_1+\calc^{n,m}_2$ where
	\beq\label{c-dec} \begin{split}
		\calc^{n,m}_1&=m\la_n\del \sum_{j=1}^{J^{n,m}} \un \rho_{f}(\wh v^{n,m,j},\wh c^{n,m,j}_{0}), \\
		\calc^{n,m}_2&=2\del \sum_{j=1}^{J^{n,m}}\sum_{r=1}^{\la_n+L-2} (m\la_n-r) \un\rho_{f}(\wh v^{n,m,j},\wh c^{n,m,j}_{r}). \end{split}
	\eeq
	
	We now claim that for every $m\in\N$ and $t\in[0,T]$, we have as $n\to\infty$,
	\beq\label{C12m}\begin{split}
		\calc^{n,m}_1\stackrel{\bbp}{\longrightarrow} \calc^m_1&=\frac{m}{m+1}\int_0^t \rho_{f,f}(0;\un\si^2(s,x))\,\dd s,\\
		\calc^{n,m}_2\stackrel{\bbp}{\longrightarrow} \calc^m_2&=2\frac{m}{m+1}\sum_{r=1}^\infty \int_0^t \rho_{f,f}(r;\un\si^2(s,x))\,\dd s,
	\end{split}\eeq
	from which \eqref{Cmt} would immediately follow. 
	We only show $\calc^{n,m}_2\stackrel{\bbp}{\longrightarrow} \calc^m_2$; the  convergence of $\calc^{n,m}_1$ is simpler and proved in a similar way. To begin with, we show that the ``$-r$''-part in $\calc^{n,m}_2$ is negligible. Indeed, as in \eqref{Vnm2},
	\beq\label{term-3}\begin{split}
		&\bbe\lb \lv 2\del\sum_{j=1}^{J^{n,m}} \sum_{r=1}^{\la_n+L-2} r \un\rho_{f}(\wh v^{n,m,j},\wh c^{n,m,j}_{r})\rv\rb \\ &\quad\lec  \del J^{n,m}\sum_{r=1}^{\la_n+L-2} r  (\ov \Ga'_r)^2 \lec \frac 1 m \sum_{r=1}^{\la_n+L-2} \frac{r}{\la_n}(\ov \Ga'_r)^2\to0 \end{split} \eeq
	as $n\to\infty$ by  dominated convergence and the square-summability of $\ov\Ga'_r$. 
	
	Next, we multiply the remaining part of $\calc^{n,m}_2$ by $\frac{m+1}{m}$ and decompose the difference between the resulting term 
	and $2\sum_{r=1}^\infty \int_0^t \rho_{f,f}(r;\un\si^2(s,x))\,\dd s$ as $$\sum_{l=1}^{5} F^{n,m}_l,\qquad F^{n,m}_l=2\sum_{r=1}^{\la_n+L-2} F^{n,m}_{l,r},$$ where $F^{n,m}_{l,r}$ is given as follows: writing $\un\rho_f(\wh v^{n,m,j},\wh c^{n,m,j}_r\,|\, A \to B)$ for the expression $\un\rho_f(\wh v^{n,m,j},\wh c^{n,m,j}_r)$, where in the definitions \eqref{variances}, \eqref{covariances}, and \eqref{covariances-2} of $\wh v^{n,m,j}$ and $\wh c^{n,m,j}_r$ (with $i[j-1,n,m]$ instead of $i(j-1,n,m)$), we also change the term $A$ to $B$, and further using the notation $t^{n,m}=J^{n,m}((m+1)\la_n+L-1)\del$, we define 
	\begin{align*}
	F^{n,m}_{1,r}&=(m+1)\la_n\del \sum_{j=1}^{J^{n,m}}  \left\{\un\rho_f\lp \wh v^{n,m,j},\wh c^{n,m,j}_r \rp \right.\\
	&\quad-\left.\un\rho_f\lp \wh v^{n,m,j},\wh c^{n,m,j}_r  \,|\, t^{n,i[j-1,n,m]}_{\la_n} \to t^{n,i[j-1,n,m]}_0 \rp\right\},\\
	F^{n,m}_{2,r}&= \sum_{j=1}^{J^{n,m}} \int_0^t\left\{ \un\rho_f\lp \wh v^{n,m,j},\wh c^{n,m,j}_r  \,|\, t^{n,i[j-1,n,m]}_{\la_n} \to t^{n,i[j-1,n,m]}_0 \rp\right.\\
	&\qquad -\left.\rho_f\lp \wh v^{n,m,j},\wh c^{n,m,j}_r  \,|\, t^{n,i[j-1,n,m]}_{\la_n} \to u-s \rp\right\}\\
	&\quad\times\bone_{(i[j-1,n,m]\del,i[j,n,m]\del)}(u)\,\dd u,\\
	F^{n,m}_{3,r}&=\int_0^{t^{n,m}} \left\{ \un\rho_f\lp \wh v^{n,m,j},\wh c^{n,m,j}_r  \,|\, t^{n,i[j-1,n,m]}_{\la_n} \to u-s \rp\right.\\
	&\quad	-\un\rho_f\lp \wh v^{n,m,j},\wh c^{n,m,j}_r  \,|\, t^{n,i[j-1,n,m]}_{\la_n} \to u-s,\right.\\
	&\qquad\quad\left.\left.\text{integral upper bounds} \to \infty \rp\right\}\,\dd u,   \\ 
	F^{n,m}_{4,r}&=\int_0^{t^{n,m}} \left\{ \un\rho_f\lp \wh v^{n,m,j},\wh c^{n,m,j}_r  \,|\, t^{n,i[j-1,n,m]}_{\la_n} \to u-s,\right.\right.\\
	&\qquad\left.\left.\vphantom{ t^{n,i[j-1,n,m]}_{\la_n}}\text{integral upper bounds} \to \infty \rp - \rho_{f,f}(r;\un\si^2(u,x))\vphantom{t^{n,i[j-1,n,m]}_{\la_n}} \right\}\,\dd u,  
	\end{align*}
	while for $l=5$, we define
	\begin{align*}
	F^{n,m}_{5}&= 2\sum_{r=1}^{\la_n+L-2}\int_0^{t^{n,m}}\rho_{f,f}(r;\un\si^2(u,x))\,\dd u - 2\sum_{r=1}^\infty \int_0^t \rho_{f,f}(r;\un\si^2(u,x))\,\dd u.
	\end{align*}
	
	The proof is complete once we show that $\bbe[|F^{n,m}_l|]\to0$ as $n\to\infty$ for all $l=1,\ldots,5$.
	For $F^{n,m}_1$, 
	if we apply Lemma~\ref{lem-rhoineq}, which is shown after the current proof, to the expression in braces in the definition of $F^{n,m}_{1,r}$ and then Hölder's inequality on the resulting product (with exponents $\frac{2p}{2p-1}$ and $2p$), we obtain from the elementary identity
	\begin{align*} \si(r,y)\si(r,z)-\si(u,y)\si(u,z)&=\si(u,z)(\si(r,y)-\si(u,y))\\ &\quad+\si(u,y)(\si(r,z)-\si(u,z))\\
	&\quad+(\si(r,y)-\si(u,y))(\si(r,z)-\si(u,z)) \end{align*}
	and Lemma~\ref{mom-ex} that
	\begin{align*}
	\bbe[|F^{n,m}_{1,r}|] &\lec \sup\left\{ \bbe[ |\si(r,y)-\si(u,y)|^{2p}]^{\frac{1}{4p}}\colon |r-u|\leq \la_n\del, y\in\R^d\right\}\\ &\lec (\la_n\del)^{\frac1 4} \to0\as.
	\end{align*}
	Furthermore, by similar considerations as in \eqref{term-3}, we have the upper bound $\bbe[|F^{n,m}_{1,r}|] \lec (\ov \Ga'_r)^2$, which is summable and independent of $n$. Hence, the dominated convergence theorem implies that $\bbe[|F^{n,m}_1|]\to0$ as $n\to\infty$.
	
	Similar arguments can be employed to show that $F^{n,m}_2$, $F^{n,m}_3$, and $F^{n,m}_4$ are negligible. Indeed, we only need to show respectively that the difference of the terms appearing in the argument of $\un\rho_f$ converges to $0$ in $L^1$. For $l=2$, this follows from $|t^{n,i[j-1,n,m]}_{0}-(u-s)|\lec m\la_n\del$ for $s\in(0,(\la_n+q-1)\del)$ or $s\in (0,(\la_n-r+q'-1)\del)$, and from $u\in (i[j-1,n,m]\del,i[j,n,m]\del)$; for $l=3$, it follows from Lemma~\ref{delta}; and for $l=4$, upon realizing that \begin{align*}&\un\rho_f\left((v_j\Ga_{|k_1-k_2|}\bone_{j_1=j_2=j})_{j_1,j_2,k_1,k_2=1}^{N,N,L,L},(v_j\Ga_{|k_1-k_2+r|}\bone_{j_1=j_2=j})_{j_1,j_2,k_1,k_2=1}^{N,N,L,L}\right)\\ &\quad=\rho_{f,f}(r;v_1,\ldots,v_N),\end{align*}
	one can use a similar argument to that given in the proof of Lemma~\ref{cond-exp-conv} for the term $A^{n,\eps,2}_x(t)$.
	
	Finally, since $|t-t^{n,m}|\leq L\del+((m+1)\la_n+L-1)\del\to0$ as $n\to\infty$ and the expectation of the $\rho_{f,f}(r;\cdot)$-term in $F^{n,m}_5$ is bounded by a constant times $\ov \Ga_r^2$ [cf.\ \eqref{term-3}], we have $\bbe[F^{n,m}_5|]\to0$ as $n\to\infty$ by dominated convergence.

	\vspace{0.5\baselineskip}
	\noindent{\itt{Property (2): }}
	We show \eqref{conv-r} for $q=4$. Then the left-hand side of \eqref{conv-r} is given by
	\begin{align*}
	\sum_{j=1}^{J^{n,m}} \sum_{k_1,\ldots,k_4=1}^{m\la_n} \bbe\lb \prod_{\iota=1}^4 \wh \psi^{n}_{i[j-1,n,m]+k_\iota,k_\iota} \,\Big|\,   \calf^n_{i[j-1,n,m]-\la_n}\rb.
	\end{align*}
	Since each $\wh \psi$-term has conditional expectation zero and two of them, one with $k_{\iota_1}$ and one with $k_{\iota_2}$, are conditionally independent if $|k_{\iota_1}-k_{\iota_2}|\geq \la_n+L-1$, the expectation in the last display vanishes unless the four indices $k_1$, $k_2$, $k_3$, and $k_4$ can be divided into two pairs $(k_{\iota_1},k_{\iota_2})$ and $(k_{\iota_3},k_{\iota_4})$ such that $|k_{\iota_1}-k_{\iota_2}|\leq\la_n+L-2$ and $|k_{\iota_3}-k_{\iota_4}|\leq\la_n+L-2$. Thus, the expression in the last display is bounded by a constant times
	\beq\label{allclose}\begin{split}
		&\sum_{j=1}^{J^{n,m}} \sum_{|k_1-k_2|\leq \la_n+L-2} \sum_{|k_3-k_4|\leq \la_n+L-2} \bbe\lb \prod_{\iota=1}^4 \wh \psi^{n}_{i[j-1,n,m]+k_\iota,k_\iota} \,\Big|\,   \calf^n_{i[j-1,n,m]-\la_n}\rb.
	\end{split}\eeq
	
	We distinguish two further cases. First, suppose that $(k_1,k_2)$ and $(k_3,k_4)$ are far from each other, that is, $\min\{|k_1-k_3|,|k_1-k_4|,|k_2-k_3|,|k_2-k_4|\} \geq \la_n+L-1$. Then the expectation in \eqref{allclose} factorizes into two parts, and similarly to \eqref{Vhat-3} and \eqref{c-dec}, the resulting expression equals
	\beq\label{eq-3}\begin{split}
		&\sum_{j=1}^{J^{n,m}} \lp\sum_{|k_1-k_2|\leq \la_n+L-2}  \bbe\lb  \wh \psi^{n}_{i[j-1,n,m]+k_1,k_1}\wh\psi^{n}_{i[j-1,n,m]+k_2,k_2} \,\Big|\,   \calf^n_{i[j-1,n,m]-\la_n}\rb\rp^2\\
		&\quad = \sum_{j=1}^{J^{n,m}} \lp m\la_n\del \un\rho_{f}(\wh v^{n,m,j},\wh c^{n,m,j}_{0}) \vphantom{\sum_{r=1}^{\la_n+L-2}}  \right.\\
		&\quad\quad\left.+2\del\sum_{r=1}^{\la_n+L-2} (m\la_n-r) \un\rho_{f}(\wh v^{n,m,j},\wh c^{n,m,j}_{r})\rp^2.
	\end{split}\eeq
	By Lemma~\ref{si-bounded}, we may assume that $\si$ is a bounded random field. Hence, by Lemma~\ref{cov-add}, we have $\rho_{f}(\wh v^{n,m,j},\wh c^{n,m,j}_{r})\leq C(\ov\Ga'_r)^2$, where $C$ is a deterministic constant that is independent of all indices, and where $\ov\Ga'_r$ was defined in \eqref{Gaprime}. Consequently, the first summand on the right-hand side of \eqref{eq-3} gives rise to a term of order $\Delta_n^2 {J^{n,m}}(m\la_n)^2\lec m\la_n\del$, which converges to $0$ as $n\to\infty$. The second term in \eqref{eq-3} is also negligible as it is of order $(m\la_n\Delta_n)^2{J^{n,m}}(\sum_{r=1}^{\la_n} (\ov\Ga'_r)^2)^2\lec  m\la_n\del$.
	
	The second case  $\min\{|k_1-k_3|,|k_1-k_4|,|k_2-k_3|,|k_2-k_4|\} \leq  \la_n+L-2$  in \eqref{allclose} yields
	\beq\label{termleft}
	\sum_{j=1}^{J^{n,m}} \sum_{k_1,\ldots,k_4} \bbe\lb \prod_{\iota=1}^4 \wh \psi^{n}_{i[j-1,n,m]+k_\iota,k_\iota} \,\Big|\,   \calf^n_{i[j-1,n,m]-\la_n}\rb, 
	\eeq
	where $\sum_{k_1,\ldots,k_4}$ denotes the sum over all quadruples $(k_1,\ldots,k_4)$ as just described. By Lemma~\ref{f-poly}, we may assume without loss of generality that $f$, which is hidden in the $\wh \psi$-variables, is an even polynomial, say, of degree $D\in\N_0$. 
	
	We will prove that
	\beq\label{toprove-2} \sum_{k_1,\ldots,k_4} \bbe\lb \prod_{\iota=1}^4 \wh \psi^{n}_{i[j-1,n,m]+k_\iota,k_\iota} \,\Big|\,   \calf^n_{i[j-1,n,m]-\la_n}\rb \leq Cm\la_n^3\del^2 \eeq
	for some $C>0$ that can be chosen independently of $j$ and uniformly for all polynomials $f$ with degree $D$ for which the sum of the absolute values of its coefficients is bounded by a given finite number, say, $A$. Then \eqref{conv-r} is evident from the first part of the proof and the fact that $J^{n,m}m\la_n^3\del^2\lec \del^{1-2a} \to 0$ for values of $a$ that are sufficiently close to $\frac 1 {2\nu}$.
	
	The variables $\wh \xi^n_{i[j-1,n,m]+k_1,k_1},\ldots,\wh \xi^n_{i[j-1,n,m]+k_4,k_4}$, which determine the $\wh\psi$-variables via \eqref{hatpsi}, are jointly Gaussian with identical law, conditionally on $\calf^n_{i[j-1,n,m]-\la_n}$. Let $T\colon \R^{N\times L}\to \bbr^{N\times L}$ be an orthogonal linear mapping such that  for each $k$, $T(\wh \xi^n_{i[j-1,n,m]+k,k})$ is a matrix of independent entries. Since $f\circ T^{-1}$ is still an even polynomial of degree $D$ with the sum of its coefficients bounded by $(NL)^{\frac D 2}A$, we may assume that each $\wh \xi$-matrix has independent entries. Moreover, as \eqref{toprove-2} is stable under approximation, we may assume that the variances of these components are strictly positive and, because of Lemma~\ref{si-bounded}, even all equal to $1$ after appropriate scaling. Then, 
	with notation from the proof of Lemma~\ref{cov-add}, we can expand $f(z)=\sum_{|l|\,\text{even},\,|l|\leq D} a_l(f) H_l(z)$ where $|a_l(f)|\leq l! A'$ for all $l\in\N_0^{N\times L}$ and $A'$ only depends on $A$ and the uniform bound on $\si$. Since the zeroth-order term is removed by the conditional expectation in \eqref{hatpsi}, we have
	\begin{align*}
	&\sum_{k_1,\ldots,k_4} \bbe\lb \prod_{\iota=1}^4 \wh \psi^{n}_{i[j-1,n,m]+k_\iota,k_\iota} \,\Big|\,   \calf^n_{i[j-1,n,m]-\la_n}\rb\\
	&\quad = \Delta_n^2 \sum_{k_1,\ldots,k_4} \sum_{\begin{smallmatrix}|l^{(1)}|,\ldots,|l^{(4)}|\,\text{even},\\ 2\leq|l^{(1)}|,\ldots,|l^{(4)}|\leq D\end{smallmatrix}}   a_{l^{(1)}}(f)\cdots a_{l^{(4)}}(f)\\
	&\qquad\times \bbe\lb \prod_{\iota=1}^4 H_{l^{(\iota)}}(\wh \xi^{n}_{i[j-1,n,m]+k_\iota,k_\iota}) \,\bigg|\, \calf^n_{i[j-1,n,m]-\la_n} \rb.
	\end{align*}
	
	Under our hypotheses and
	conditionally on $\calf^n_{i[j-1,n,m]-\la_n}$, the four matrices inside the arguments of the Hermite polynomials are jointly Gaussian, each with independent $N(0,1)$ entries, and with cross-covariances given by $\wh c^{n,m,j}_{|k_\iota-k_{\iota'}|}$ between any two of them. By the definition of the multivariate Hermite polynomials, the final conditional expectation in the previous display is taken for a product of altogether $\nu= \sum_{\iota=1}^4\sum_{p=1}^N \sum_{q=1}^L l^{(\iota)}_{pq}$ univariate Hermite polynomials. So we can apply Proposition~4.3.22 in \citeB{Pipiras17} to evaluate this conditional expectation. Noticing a different prefactor in the definition of the Hermite polynomials, we obtain  
	\beq\label{term}\begin{split}
		&\Delta_n^2\sum_{k_1,\ldots,k_4} \sum_{\begin{smallmatrix}|l^{(1)}|,\ldots,|l^{(4)}|\,\text{even},\\ 2\leq |l^{(1)}|,\ldots,|l^{(4)}|\leq D\end{smallmatrix}}   \frac{a_{l^{(1)}}(f)\cdots a_{l^{(4)}}(f)}{l^{(1)}!\cdots l^{(4)}!}\\
		&\quad\times \sum_{\begin{smallmatrix} \si\in \calm^0_2([\nu],\pi^\ast)\colon \\ \ell^{\iota,\iota}_{pq,p'q'}(\si) = 0 \end{smallmatrix}} \prod_{\begin{smallmatrix} 1\leq \iota< \iota'\leq 4,\\ (p,q),\,(p',q')\end{smallmatrix}}  (\wh c^{n,m,j}_{|k_{\iota}-k_{\iota'}|})_{pq,p'q'}^{\ell^{\iota,\iota'}_{pq,p'q'}(\si)}\end{split}
	\eeq
	with the following specifications: Taking the terminology of Definitions~4.3.1 and 4.3.2 in \citeB{Pipiras17} and the notation $\calm^0_2([\nu],\pi^\ast)$ on page 242 of the same reference for granted, $\ell^{\iota,\iota'}_{pq,p'q'}(\si)$ is the number of closed curves connecting the blocks $(\iota,p,q)$ and $(\iota',p',q')$ of $\pi^\ast$ in the diagram $\Ga(\pi^\ast,\si)$, where $\si$ is a generic element of $\calm^0_2([\nu],\pi^\ast)$ (and should not be confused with the volatility process $\si$). Furthermore, the product in \eqref{term} is taken for all $1\leq \iota<\iota'\leq 4$, $p,p'=1,\ldots,N$, and $q,q'=1,\ldots,L$. The reason why we have the restriction $\ell^{\iota,\iota}_{pq,p'q'}(\si) = 0$ in \eqref{term} is because the entries of each $\hat \xi$-matrix are assumed to be conditionally independent. The only properties that we need from these objects are the following: For all $|l^{(1)}|,\ldots,|l^{(4)}|\leq D$ and $D$ fixed, we have
	\beq\label{properties1}  \#\calm^0_2([\nu],\pi^\ast)\lec 1. \eeq
	And for all $\iota,\iota'=1,\ldots,4$, $p,p'=1,\ldots,N$, $q,q'=1,\ldots,L$, and $\si\in \calm^0_2([\nu],\pi^\ast)$, we have
	\beq\label{properties2}
	\sum_{\iota'=1}^4\sum_{p'=1}^N \sum_{q'=1}^L \ell^{\iota,\iota'}_{pq,p'q'}(\si) = l^{(\iota)}_{pq},\qquad \sum_{\iota=1}^4\sum_{p=1}^N \sum_{q=1}^L \ell^{\iota,\iota'}_{pq,p'q'}(\si) = l^{(\iota')}_{p'q'}.\eeq
	
	Recalling that $|a_l(f)|\lec l!$, we now deduce that \eqref{term} is bounded by a constant multiple of
	\begin{align*}
	&\Delta_n^2 \sum_{k_1,\ldots,k_4} \sum_{\begin{smallmatrix}|l^{(1)}|,\ldots,|l^{(4)}|\,\text{even},\\ 2\leq |l^{(1)}|,\ldots,|l^{(4)}|\leq D\end{smallmatrix}}\sum_{\begin{smallmatrix} \si\in \calm^0_2([\nu],\pi^\ast)\colon \\ \ell^{\iota,\iota}_{pq,p'q'}(\si) = 0 \end{smallmatrix}} \prod_{\begin{smallmatrix} 1\leq \iota< \iota'\leq 4,\\ (p,q),\,(p',q')\end{smallmatrix}}  \lv(\wh c^{n,m,j}_{|k_{\iota}-k_{\iota'}|})_{pq,p'q'}\rv^{\ell^{\iota,\iota'}_{pq,p'q'}(\si)}.
	\end{align*}
	For the quadruplets $(k_1,\ldots,k_4)$ under consideration, if $k_1\leq k_2\leq k_3\leq k_4$, we must have $r_1=k_2-k_1\leq\la_n+L-2$, $r_2=k_3-k_2\leq\la_n+L-2$, and $r_3=k_4-k_3\leq \la_n+L-2$. Hence, by a change of variables, and because $|\wh c^{n,m,j}_r|\lec \ov\Ga'_r\lec 1$, we derive
	\begin{align*}
	&\sum_{k_1,\ldots,k_4} \bbe\lb \prod_{\iota=1}^4 \wh \psi^{n}_{i[j-1,n,m]+k_\iota,k_\iota} \,\Big|\,   \calf^n_{i[j-1,n,m]-\la_n}\rb\\
	&\quad \lec \Delta_n^2 m \la_n\sum_{r_1,r_2,r_3=0}^{\la_n+L-2}  \sum_{\begin{smallmatrix}|l^{(1)}|,\ldots,|l^{(4)}|\,\text{even},\\ 2\leq |l^{(1)}|,\ldots,|l^{(4)}|\leq D\end{smallmatrix}}\sum_{\begin{smallmatrix} \si\in \calm^0_2([\nu],\pi^\ast)\colon \\ \ell^{\iota,\iota}_{pq,p'q'}(\si) = 0 \end{smallmatrix}} \prod_{(p,q),\,(p',q')}   \\
	&\quad \quad \times \lv(\wh c^{n,m,j}_{r_1})_{pq,p'q'}\rv^{\ell^{1,2}_{pq,p'q'}(\si)}\lv(\wh c^{n,m,j}_{r_1+r_2})_{pq,p'q'}\rv^{\ell^{1,3}_{pq,p'q'}(\si)} \\
	&\quad\quad\times \lv(\wh c^{n,m,j}_{r_1+r_2+r_3})_{pq,p'q'}\rv^{\ell^{1,4}_{pq,p'q'}(\si)}  \lv(\wh c^{n,m,j}_{r_2})_{pq,p'q'}\rv^{\ell^{2,3}_{pq,p'q'}(\si)}\\
	&\quad\quad\times \lv(\wh c^{n,m,j}_{r_2+r_3})_{pq,p'q'}\rv^{\ell^{2,4}_{pq,p'q'}(\si)} \lv(\wh c^{n,m,j}_{r_3})_{pq,p'q'}\rv^{\ell^{3,4}_{pq,p'q'}(\si)}\\
	&\quad \lec \Delta_n^2 m\la_n  \sum_{\begin{smallmatrix}|l^{(1)}|,\ldots,|l^{(4)}|\,\text{even},\\ 2\leq |l^{(1)}|,\ldots,|l^{(4)}|\leq D\end{smallmatrix}}\sum_{\begin{smallmatrix} \si\in \calm^0_2([\nu],\pi^\ast)\colon \\ \ell^{\iota,\iota}_{pq,p'q'}(\si) = 0 \end{smallmatrix}}  \sum_{r_1=0}^{\la_n+L-2} (\ov\Ga'_{r_1})^{\sum \ell^{1,2}(\si)} \\
	&\quad\quad\times\sum_{r_2=0}^{\la_n+L-2} (\ov\Ga'_{r_1+r_2})^{\sum \ell^{1,3}(\si)}\sum_{r_3=0}^{\la_n+L-2} (\ov\Ga'_{r_1+r_2+r_3})^{\sum \ell^{1,4}(\si)}\\
	&\quad\leq \Delta_n^2 m\la_n  \sum_{\begin{smallmatrix}|l^{(1)}|,\ldots,|l^{(4)}|\,\text{even},\\ 2\leq |l^{(1)}|,\ldots,|l^{(4)}|\leq D\end{smallmatrix}}\sum_{\begin{smallmatrix} \si\in \calm^0_2([\nu],\pi^\ast)\colon \\ \ell^{\iota,\iota}_{pq,p'q'}(\si) = 0 \end{smallmatrix}}  \sum_{r_1=0}^{\la_n+L-2} (\ov\Ga'_{r_1})^{\sum \ell^{1,2}(\si)} \\
	&\qquad\times\sum_{r_2=0}^{2(\la_n+L-2)} (\ov\Ga'_{r_2})^{\sum\ell^{1,3}(\si)}\sum_{r_3=0}^{3(\la_n+L-2)} (\ov\Ga'_{r_3})^{\sum \ell^{1,4}(\si)},
	\end{align*}
	where $\sum \ell^{1,\iota'}(\si)$ stands for $\sum_{(p,q),\,(p',q')} \ell^{1,\iota'}_{pq,p'q'}(\si)$.
	
	By \eqref{properties2}, the three exponents $\sum \ell^{1,2}(\si)$, $\sum \ell^{1,3}(\si)$, and $\sum \ell^{1,4}(\si)$ add up to
	\begin{align*} \sum_{1<\iota'\leq 4} \sum_{(p,q),\,(p',q')} \ell^{1,\iota'}_{pq,p'q'}(\si) &= \sum_{1<\iota'\leq 4} \sum_{p=1}^N\sum_{q=1}^L\sum_{p'=1}^N\sum_{q'=1}^L \ell^{1,\iota'}_{pq,p'q'}(\si)=\sum_{p=1}^N\sum_{q=1}^L l^{(1)}_{pq} \\ &= |l^{(1)}|\geq2. \end{align*}
	So exactly two different situations may occur: either at least one of these three exponents is $\geq2$, or there are at least two of them equal to $1$. In the first case, because of \eqref{properties1} and because $\ov\Ga'_r$ is square-summable, the resulting term is of magnitude $\Delta_n^2 m\la_n \la_n^2$, as desired. In the second case, Hölder's inequality implies that $\sum_{r=0}^{\la_n} \ov\Ga'_r \leq (\sum_{r=0}^{\la_n} (\ov \Ga'_r)^2)^{1/2} (\la_n+1)^{1/2} \lec \la_n^{1/2}$, so the resulting term is of order $\Delta_n^2m\la_n (\la_n^{1/2})^2\la_n$, which proves \eqref{toprove-2}.

	\vspace{0.5\baselineskip}
	\noindent{\itt{Property (3): }} It suffices to consider the case where $f$ is real-valued.
	Then the desired result \eqref{conv-mart} obviously follows if we can show that
	\beq\label{toshow} \bbe\lb \wh \psi^n_{i[j-1,n,m]+k,k}\left(M({\tau^{n}_j})-M({\tau^{n}_{j-1}})\right) \,|\, \calf^n_{i[j-1,n,m]-\la_n}\rb = 0  \eeq
	for all  $n$, $m$, $j$,  and $k$. 
	By an orthogonal transformation, as explained in the paragraph after \eqref{toprove-2}, we may assume that for each $k$, the matrix $\wh\xi^n_{i[j-1,n,m]+k,k}$ consists of independent entries.
	Then, since $f$ is even, conditionally on $\calf^n_{i[j-1,n,m]-\la_n}$, the variable $\wh\psi^n_{i[j-1,n,m]+k,k}$ belongs to the direct sum of Wiener chaoses of even orders $\geq 2$, that is, with the notation from Lemma~\ref{cov-add},
	\beq\label{exp-psi} \wh\psi^n_{i[j-1,n,m]+k,k}=\sum_{\begin{smallmatrix} l\in\N_0^{N\times L},\\ |l|\geq2\,\text{even}\end{smallmatrix}} a_{l}(\wh v^{n,m,j},f) H_{l}(\wh v^{n,m,j},\wh\xi^n_{i[j-1,n,m]+k,k}),\eeq
	where we identify $\wh v^{n,m,j}$ with $(\wh v^{n,m,j}_{pq,pq})_{p,q=1}^{N,L}$. 
	
	By contrast, the Brownian increment $W^{\iota}(\tau^{n}_j)-W^{\iota}(\tau^{n}_{j-1})$ belongs to the first-order Wiener chaos. In other words, if we let $Y \in \R^{N\times L}$ be a matrix where $Y_{11}= \sqrt{{\wh v^{n,m,j}_{11}}}(W^{\iota}({\tau^{n}_j})-W^{\iota}({\tau^{n}_{j-1}}))/\sqrt{{\tau^n_j-\tau^n_{j-1}}}$ and all other entries are independent of $W^\iota$ with the same law as the corresponding entries of $\wh \xi^n_{i[j-1,n,m]+k,k}$, then $W^\iota(\tau^n_j)-W^\iota(\tau^n_{j-1}) = \sum_{|l|=1} a_l(\wh v^{n,m,j},g)H_l(\wh v^{n,m,j},Y)$ where $g(y)=\sqrt{{\tau^n_j-\tau^n_{j-1}}}y_{11}/\sqrt{{\wh v^{n,m,j}_{11}}}$. Hence, by \eqref{cross}, $\wh \psi^n_{i[j-1,n,m]+k,k}$ and $W^{\iota}(\tau^{n}_j)-W^{\iota}(\tau^{n}_{j-1})$ are conditionally uncorrelated, from which \eqref{toshow} follows for $M=W^\iota$. 
	
	Now suppose that $M$ is a bounded martingale that is orthogonal to all $W^\iota$. By the defining properties of $W^\iota$, $M$ is actually orthogonal to all stochastic integral processes with respect to $W$. Since each variable $H_{l}(\wh v^{n,m,j},\wh\xi^n_{i[j-1,n,m]+k,k})$ can be written in such a form (see Proposition~1.1.4 in \citeB{Nualart06-2}), \eqref{toshow} follows from \eqref{exp-psi}. 
	\epr
	
	\blem\label{lem-rhoineq} Suppose that $f$ satisfies Assumption~\ref{AssCLT}1, and recall the definition of $\un\rho_f(v,c)$ introduced after \eqref{Vhat-3}. Then for all $v, v', c, c'\in (\R^{N\times L})^2$ such that $\un\rho_f(v,c)$ and $\un\rho_f(v',c')$ are well defined,
	\beq\label{rhoineq} |\un\rho_f(v,c)-\un\rho_f(v',c')|\lec (1+|v|^{p-\frac 12}+|v'|^{p-\frac12})(\sqrt{|v-v'|} + \sqrt{|c-c'|}).\eeq
	\elem
	\bpr We only give the details in the case where $N=L=1$, so we have $v,v'\in[0,\infty)$, $c\in[0,v]$, and $c'\in[0,v']$. The general case follows similarly if we keep in mind that all norms on $\R^{N\times L}$ are equivalent, and that $\sum_{j=1}^N\sum_{k=1}^L \sqrt{|x_{jk}-y_{jk}|} \lec \sqrt{|x-y|}$. 
	
	Consider the auxiliary function $\wt\rho_f(x,y)=\cov(f(xU+yW),f(xV+yW))$ for $x,y\in\R$, where $U$, $V$, and $W$ are independent standard normal random variables. Since $f'$ exists and $|f'(z)|\lec 1+|z|^{p-1}$, we can differentiate under the expectation sign by dominated convergence, so together with the Cauchy--Schwarz inequality, we obtain
	\begin{align*}
	\frac{\partial}{\partial x}\wt \rho_f(x,y)&=\frac{\partial}{\partial x}\left\{ \bbe[f(xU+yW)f(xV+yW)] - \bbe[f(xU+yW)]^2\right\}\\
	&=\bbe[Uf'(xU+yW)f(xV+yW)]\\
	&\quad+\bbe[Vf(xU+yW)f'(xV+yW)] \\
	&\quad- 2\bbe[f(xU+yW)]\bbe[Uf'(xU+yW)]\\
	&=2\cov(Uf'(xU+yW),f(xV+yW))\\
	&\leq 2\sqrt{\var(Uf'(xU+yW))\var(f(xV+yW))}\\
	&\lec \sqrt{(1+(x^2+y^2)^{p-1})(1+(x^2+y^2)^p)}\lec 1+|x|^{2p-1}+|y|^{2p-1}.
	\end{align*}
	Similarly,
	\[\frac{\partial}{\partial y}\wt \rho_f(x,y) = 2\cov(Wf'(xU+yW),f(xU+yW))\lec 1+|x|^{2p-1}+|y|^{2p-1}. \]
	Thus, by the mean value theorem and the inequality $|\sqrt{x}-\sqrt{y}|\leq \sqrt{|x-y|}$ for $x,y\in[0,\infty)$,
	\begin{align*}
	\un\rho_f(v,c)-\un\rho_f(v',c')&=\wt \rho_f(\sqrt{v-c},\sqrt{c})-\wt\rho_f(\sqrt{v'-c'},\sqrt{c'})\\
	&\lec (1+(\sqrt{v})^{2p-1}+(\sqrt{v'})^{2p-1})\\
	&\quad\times(|\sqrt{v-c}-\sqrt{v'-c'}|+|\sqrt{c}-\sqrt{c'}|)\\
	&\leq (1+|v|^{p-\frac12}+|v'|^{p-\frac12})(\sqrt{|v-v'|}+2\sqrt{|c-c'|}),
	\end{align*}
	which shows \eqref{rhoineq}.
	\epr
	
	\bpr[Proof of Lemma~\ref{si-bounded}] Let $K\in\N$ and $\si_K = (\si \vee (-K))\wedge K$. Define the variables
	\beq\label{hatVm} \wh\xi^{n,K}_{i,k}, ~\wh \psi^{n,K}_{i,k}, ~\wh V^{n,m,K,1}(t),~ \wh V^{n,m,K}_j \eeq
	as in \eqref{hatpsi} and \eqref{hatV} but with $\si$ replaced by $\si_K$. By Proposition~2.2.4 of \citeB{Jacod12-2}, it suffices to show 
	\begin{align}\label{first-eq} \lim_{K\to\infty}\limsup_{n\to\infty} \bbe\lb \supt | \wh V^{n,m,1}(t)-\wh V^{n,m,K,1}(t)|^2\rb = 0, \\
	\lim_{K\to\infty} \bbe\lb \supt | \calc^{m,K}(t) -\calc^m(t)|\rb  =0\label{second-eq}
	\end{align}
	for each $m\in\N$, where $\calc^{m,K}$ is given by $\calc^m = \frac{m}{m+1}\calc$ when $\si$ is replaced by $\si_K$. Due to Lemma~\ref{cov-add} and the fact that $\Ga_r$ in \eqref{Ga-formula} is square-summable, \eqref{second-eq} follows from the continuity of $(v_1,\ldots,v_N)\mapsto \rho_{f_{m_1},f_{m_2}}(0; v_1,\ldots,v_N)$ and the dominated convergence theorem. 
	
	Hence, we only give the details for \eqref{first-eq}. To this end, observe that $$\wh V^{n,m,1}(t)-\wh V^{n,m,K,1}(t)=\sum_{j=1}^{J^{n,m}(t)} (\wh V^{n,m}_j - \wh V^{n,m,K}_j)$$ is a sum of martingale increments, whose predictable quadratic variation process is given by $\sum_{j=1}^{J^{n,m}(t)} \bbe[(\wh V^{n,m}_j - \wh V^{n,m,K}_j)^2 \,|\, \calf^n_{i[j-1,n,m]-\la_n}]$. Hence, the expectation in \eqref{first-eq} is bounded by a multiple of 
	\begin{align*}
	\bbe \lb  \sum_{j=1}^{J^{n,m}(T)} \bbe[(\wh V^{n,m}_j - \wh V^{n,m,K}_j)^2 \,|\, \calf^n_{i[j-1,n,m]-\la_n}] \rb.
	\end{align*}
	This has a similar form as $\calc^{n,m}$ in \eqref{Cnm}, and by the same arguments, it equals
	\beq\label{cov-dec}\begin{split}
		&\del\! \sum_{j=1}^{J^{n,m}(T)}\sum_{|k-l|<\la_n+L-1} \cov[ f(\wh \xi^n_{i[j-1,n,m]+k,k})- f(\wh \xi^{n,K}_{i[j-1,n,m]+k,k}),\\
		&\qquad\quad  f(\wh \xi^n_{i[j-1,n,m]+l,l})- f(\wh \xi^{n,K}_{i[j-1,n,m]+l,l}) \,|\, \calf^n_{i[j-1,n,m]-\la_n} ]\\
		&\quad =\del\! \sum_{j=1}^{J^{n,m}(T)} \sum_{k=1}^{m\la_n} \var[ f(\wh \xi^n_{i[j-1,n,m]+k,k})- f(\wh \xi^{n,K}_{i[j-1,n,m]+k,k}) \,|\, \calf^n_{i[j-1,n,m]-\la_n} ]\\
		&\quad \quad + \del\! \sum_{j=1}^{J^{n,m}(T)} \sum_{r=1}^{\la_n+L-2} (m\la_n-r) \cov[ f(\wh \xi^n_{i[j-1,n,m]+1,1})- f(\wh \xi^{n,K}_{i[j-1,n,m]+1,1}),\\
		&\quad\qquad  f(\wh \xi^n_{i[j-1,n,m]+1+r,1+r})- f(\wh \xi^{n,K}_{i[j-1,n,m]+1+r,1+r}) \,|\, \calf^n_{i[j-1,n,m]-\la_n} ].
	\end{split} \eeq
	Since $|f(z)-f(z')|\lec (1+|z|^{p-1}+|z'|^{p-1})|z-z'|$ by Assumption~\ref{AssCLT}1, a size argument shows that
	\begin{align*}
	&\bbe\left[ \lp f(\wh \xi^n_{i[j-1,n,m]+k,k})- f(\wh \xi^{n,K}_{i[j-1,n,m]+k,k})\rp^2\right]\\
	&\quad\lec \bbe\lb \lp \wh \xi^n_{i[j-1,n,m]+k,k} - \wh \xi^{n,K}_{i[j-1,n,m]+k,k}\rp^{2p}\rb^{\frac 1 p}.
	\end{align*}
	Recalling that $\si$ has bounded moments of order $2p+\eps$ by Lemma~\ref{mom-ex}, we conclude that 
	\beq\label{varbound}\begin{split}&\bbe\left[ \lp f(\wh \xi^n_{i[j-1,n,m]+k,k})- f(\wh \xi^{n,K}_{i[j-1,n,m]+k,k})\rp^2\right]\\ &\quad\lec \sup_{(s,y)\in\R\times\R^d} \bbe[|\si(s,y)|^{2p}\bone_{\{|\si(s,y)|>K\}}]^{\frac 1 p} \lec K^{-\frac \eps p}.\end{split}
	\eeq
	As a consequence, the first term on the right-hand side of \eqref{cov-dec} is bounded by multiple of
	\[ \Delta_n J^{n,m}(T)m\la_n K^{-\frac \eps p} \lec K^{-\frac \eps p}, \]
	which goes to $0$ as $K\to\infty$, uniformly in $n$ and $m$.
	
	We now consider the second term on the right-hand side of \eqref{cov-dec}. With similar reasoning as for \eqref{term-3}, the ``$-r$''-portion of the sum can be neglected as $n\to\infty$ because it is of size $\lec \Delta_n J^{n,m}(T) \sum_{r=1}^{\la_n+L-2} r (\ov\Ga'_r)^2$. For the remaining part, which is given by
	\beq\label{remain}\begin{split} &\del m\la_n \sum_{j=1}^{J^{n,m}(T)} \sum_{r=1}^{\la_n+L-2}  \cov\lp f(\wh \xi^n_{i[j-1,n,m]+1,1})- f(\wh \xi^{n,K}_{i[j-1,n,m]+1,1}),\right.\\
		&\quad \left. f(\wh \xi^n_{i[j-1,n,m]+1+r,1+r})- f(\wh \xi^{n,K}_{i[j-1,n,m]+1+r,1+r}) \,|\, \calf^n_{i[j-1,n,m]-\la_n} \rp,
	\end{split}\eeq
	notice from \eqref{varbound} that 
	\beq\label{dom-2}\begin{split}
		&\limsup_{n\to\infty} \del m\la_n \sum_{j=1}^{J^{n,m}(T)} \bbe\lb \lv \cov\lp f(\wh \xi^n_{i[j-1,n,m]+1,1})- f(\wh \xi^{n,K}_{i[j-1,n,m]+1,1}),\right.\right.\right.\\
		&\quad\quad \left.\left. \left.f(\wh \xi^n_{i[j-1,n,m]+1+r,1+r})- f(\wh \xi^{n,K}_{i[j-1,n,m]+1+r,1+r}) \,|\, \calf^n_{i[j-1,n,m]-\la_n} \rp\rv\rb\\
		&\quad \lec \del m\la_nJ^{n,m}(T)K^{-\frac{\eps}{p}} \lec K^{-\frac{\eps}{p}} \to 0\qquad\text{as}\qquad K\to\infty.
	\end{split} \raisetag{-3.4\baselineskip}\eeq
	By the dominated convergence theorem, the lemma is proved once we can bound the left-hand side of \eqref{dom-2} by a quantity that is independent of $K$ and summable in $r$.
	
	To this end, observe that conditionally on $\calf^n_{i[j-1,n,m]-\la_n}$, the variables $$(\wh \xi^{n,K}_{i[j-1,n,m]+1,1},\wh \xi^{n,K}_{i[j-1,n,m]+1+r,1+r})$$  and $$(\wh \xi^{n}_{i[j-1,n,m]+1,1}-\wh \xi^{n,K}_{i[j-1,n,m]+1,1},\wh \xi^{n}_{i[j-1,n,m]+1+r,1+r}-\wh \xi^{n,K}_{i[j-1,n,m]+1+r,1+r})$$
	consist of Wiener  integrals over the disjoint sets $\{|\si(t^{n,i[j-1,n,m]}_{\la_n},y)|\leq K\}$ and $\{|\si(t^{n,i[j-1,n,m]}_{\la_n},y)|> K\}$, respectively, and are hence independent. If we define $g(z,z') = f(z+z')-f(z)$ for $z,z'\in\R^{N\times L}$, then the covariance term in \eqref{remain} can be written as $\cov'(g(X_1,Y_1),g(X_2,Y_2))$
	where $(X_1,X_2,Y_1,Y_2)$ is centered Gaussian, defined on another probability space $(\Om',\calf',\bbp')$, such that $(X_1,X_2)$ is independent of $(Y_1,Y_2)$ and 
	\begin{align*}
	\cov'((X_{1|2})_{pq},(X_{1|2})_{p'q'})&=(v^{n,m,j}_{\leq K})_{pq,p'q'}, \\ \cov'((Y_{1|2})_{pq},(Y_{1|2})_{p'q'})&=(v^{n,m,j}_{> K})_{pq,p'q'},\\
	\cov'((X_{1})_{pq},(X_{2})_{p'q'})&=(c^{n,m,j}_{\leq K,r})_{pq,p'q'},\\ \cov'((Y_{1})_{pq},(Y_{2})_{p'q'})&=(c^{n,m,j}_{> K,r})_{pq,p'q'}.
	\end{align*}
	Here, $v^{n,m,j}_{\leq K}$, $c^{n,m,j}_{\leq K,r}$, $v^{n,m,j}_{>K}$, and $c^{n,m,j}_{>K,r}$ are defined as in \eqref{variances} and \eqref{covariances} but with $i(j-1,n,m)$ replaced by $i[j-1,n,m]$, $(k-l)$ replaced by $r$, and, respectively, $\si$ replaced by $\si_K$ for the variables with index $\leq K$ and by $\si-\si_K$ for the variables with index $>K$. 
	Consequently, if we apply Lemma~\ref{cov-add} to $g$, we obtain
	\begin{align*} &|\cov'(g(X_1,Y_1),g(X_2,Y_2))|\\ &\quad\lec (1+|v^{n,m,j}_{\leq K}|^{\frac{p-2}2}+|v^{n,m,j}_{> K}|^{\frac{p-2}2})(|c^{n,m,j}_{\leq K,r}|^2+|c^{n,m,j}_{> K,r}|^2),
	\end{align*}
	so by Hölder's inequality, and since $\bbe[|c^{n,m,j}_{\leq K,r}|^p]^{1 /p}+\bbe[|c^{n,m,j}_{> K,r}|^p]^{1/ p} \lec \ov\Ga'_r$,
	\beq\label{dom-1}\begin{split}
		&\limsup_{n\to\infty} \del m\la_n \sum_{j=1}^{J^{n,m}(T)}  \bbe\lb  \lv\cov\lp f(\wh \xi^n_{i[j-1,n,m]+1,1})- f(\wh \xi^{n,K}_{i[j-1,n,m]+1,1}),\right.\right.\right.\\
		&\quad \quad\left.\left. \left.f(\wh \xi^n_{i[j-1,n,m]+1+r,1+r})- f(\wh \xi^{n,K}_{i[j-1,n,m]+1+r,1+r}) \,|\, \calf^n_{i[j-1,n,m]-\la_n} \rp\rv\rb\\
		&\quad\lec \limsup_{n\to\infty} \del m\la_n J^{n,m}(T)(\ov \Ga'_r)^2 \lec (\ov \Ga'_r)^2,
	\end{split}\raisetag{-3.4\baselineskip}\eeq
	which is independent of $m$ and $K$ and is summable in $r$.
	\epr

	\bpr[Proof of Lemma~\ref{f-poly}]  Since $f$  is even and continuous, we have by the Stone--Weierstrass theorem  that for every $K\in\N$ and $m=1,\ldots,M$, there is an even polynomial $f_m^{(K)}$ in $NL$ variables such that $|f_m(z)-f_m^{(K)}(z)|\leq K^{-1}$ for $|z|\leq K$. Then we need to check \eqref{first-eq} and \eqref{second-eq} once more, where we redefine the variables in \eqref{hatVm} by taking $f_m^{(K)}$ instead of $f$ in \eqref{hatpsi} and \eqref{hatV} (we do not change $\si$ to $\si_K$, but $\si$ is nevertheless bounded by Lemma~\ref{si-bounded}). The arguments are analogous to those in Lemma~\ref{si-bounded}, so we do not present the details here. 
	\epr
	
	\bpr[Proof of Lemma~\ref{discretization}] Let us denote the left-hand side of \eqref{term-2} by $H^n(t)$. Then $H^n(t)=H^n_1(t)+H^n_1(t)$ with
	\begin{align*}
	H^n_1(t)&=\delh\sumt \sum_\al \bbe\lb \frac{\partial}{\partial z_\al} f(\theta^n_i)(\kappa^n_i)_\al \,\Big|\, \calf^n_{i-\la_n}\rb,\\
	H^n_2(t)&= \frac12\delh\sumt\sum_{\al,\beta} \bbe\lb\frac{\partial}{\partial z_\al\,\partial z_\beta} f(\ov\theta^n_i)(\kappa^n_i)_\al(\kappa^n_i)_\beta\,\Big|\, \calf^n_{i-\la_n}\rb,
	\end{align*}
	where the sums are taken over $\al,\beta\in\{1,\ldots,N\}\times\{1,\ldots,L\}$, $\kappa^n_i$ is defined by
	\[ \kappa^n_i = \sum_{q=1}^Q \iint  \frac{\un\Delta^n_i G_{x,y}(s)}{\tau_n} \lp\si(s,y)-\si\lp t^{n,i}_{\la_n^{(q-1)}},y\rp\rp\bone^{n,i}_{\la_n^{(q-1)},\la_n^{(q)}}(s) \,W(\dd s,\dd y), \]
	and $\ov\theta^n_i$ is some value between $\frac{\ga^{n,i}_x}{\tau_n}$ and $\theta^n_i$. Regardless of the exact choice of the numbers $a^{(q)}$, by Assumption~\ref{AssCLT}1, Lemma~\ref{mom-ex}, and a standard size argument, we have that
	\beq\label{Hn2} \bbe\lb \supt |H^n_2(t)|\rb \lec \sum_{q=1}^Q \del^{-\frac12}((\la^{(q-1)}_n\del)^{\frac12})^2 \lec \del^{\frac12 -a},  \eeq
	which converges to $0$ as $n\to\infty$ if $a$ is close enough to $\frac 1 {2\nu}$.
	
	For the term $H^n_1(t)$, since $\kappa^n_i$ depends linearly on $\si$, it is no restriction to assume that $\si=\si^{(0)}$ or $\si=\si^{(1)}=\si-\si^{(0)}$ in the definition of $\kappa^n_i$. In the former case, we have by a size argument,
	\beq\label{eq-4} \begin{split} \bbe\lb \supt |H^n_1(t)|\rb  &\lec \sum_{q=1}^Q \Delta_n^{-\frac 12} (\la_n^{(q-1)}\del)^\ga \del^{\frac \nu 2 a^{(q)}}\\ &\leq  \sum_{q=1}^Q \Delta_n^{-\frac12 +(1-a^{(q-1)})\ga +\frac\nu2 a^{(q)}},\end{split}\eeq
	with the convention that $a^{(Q)}=0$. Then \eqref{eq-4} goes to $0$ as $n\to\infty$ if 
	\beq\label{cond-a}-\frac12 +(1-a^{(q-1)})\ga +\frac\nu2 a^{(q)}>0\eeq
	for every $q=1,\ldots,Q$.
	Putting equality in \eqref{cond-a} gives rise to the recurrence equation
	\beq\label{recurrence}-\frac12 +(1-b^{(q-1)})\ga +\frac\nu2 b^{(q)}=0. \eeq 
	
	Then we could take the starting point $b^{(0)}=a$, and since $\nu>1$, $\ga>\frac12$, and $a$ is close to $\frac 1 {2\nu}$, it would be straightforward to show that this recurrence equation decreases to a strictly negative fixed point or $-\infty$. Hence, if $Q$ is such that $b^{(Q)}$ is the first strictly negative value of \eqref{recurrence}, then it would suffice to choose $a^{(q)}$ slightly  larger than $b^{(q)}$ for $q=1,\ldots,Q-1$ to meet the requirements of this lemma. But for reasons that will become clear in the proof of Lemma~\ref{iter-condexp}, we do \emph{not} take this construction. Instead, we shall solve \eqref{recurrence} backwards by specifying a terminal value $b^{(Q)} = b$, where $b$ is a strictly negative number, and then compute $b^{(Q-1)},\ldots,b^{(0)}$ iteratively. In explicit terms, we have
	\beq\label{bQ} b^{(Q-r)} = \frac{(\textstyle\frac{\nu}{2\ga})^r((2\ga-\nu)b-(2\ga-1))+(2\ga-1)}{2\ga-\nu},\qquad r=0,\ldots,Q. \eeq
	If $\nu=2\ga$, we have $b^{(Q-r)}=b+(1-\frac1 \nu) r$. Since $\nu>1$ and $\ga>\frac12$, as $r\to\infty$, the right-hand side of \eqref{bQ} tends to $\frac{2\ga-1}{2\ga-\nu}>1>\frac{1}{2\nu}$ if $\nu<2\ga$; and it tends to $+\infty$ if $\nu\geq 2\ga$ and the absolute value of $b$ is sufficiently small. Hence, if $a$ is close enough to $\frac1 {2\nu}$, we can always guarantee that $b^{(0)}>a$ by taking a large number $Q$ of iterations. By further decreasing the absolute value of $b$ if necessary, we can also make sure that  $b^{(Q-1)}>0$. Thus, we obtain \eqref{cond-a} by choosing $a^{(q)}$, $q=1,\ldots,Q-1$, slightly larger than $b^{(q)}$. 
	
	Next assume that $\kappa^n_i$ is defined using $\si=\si^{(1)}$. Then we can split $\kappa^n_i = \kappa^{n,1}_i+\kappa^{n,2}_i$ where $\kappa^{n,1}_i$ and $\kappa^{n,2}_i$ are defined in the same way as $\kappa^n_i$ but with $\si=\si^{(1)}$ replaced by the second and third term on the right-hand side of \eqref{si-fixed-x}, respectively. Accordingly, we also obtain the decomposition $H^n_1(t)=H^{n,1}_1(t)+H^{n,2}_1(t)$. Since the second term on the right-hand side of \eqref{si-fixed-x} is a Lebesgue integral, $H^{n,1}_1$ is asymptotically of size at most $\del^{-1/2}\la_n\del = \del^{1/2 -a}$ and hence negligible for small $a$. Next, we further decompose $H^{n,2}_1$ into $H^{n,2}_1(t)=H^{n,2}_{11}(t)+H^{n,2}_{12}(t)+H^{n,2}_{13}(t)$ where
	\begin{align*}
	H^{n,2}_{11}(t)&= \delh\sumt \sum_{q=1}^Q \sum_{j=1}^N\sum_{k=1}^L \bbe\lb \frac\partial{\partial z_{jk}} f(\theta^n_i) \iint\frac{\Delta^n_{i+k-1} G_{x_j-y}(s)}{\tau_n} \right.\\
	&\quad \times \iint_{(i-\la^{(q-1)}_n)\del}^{s} (\si^{(12)}_y(r,z)-\si^{(12)}_y((i-\la^{(q-1)}_n)\del,z))\,W'(\dd r,\dd z)\\
	&\quad\times\left.\bone^{n,i}_{\la_n^{(q-1)},\la_n^{(q)}}(s) \,W(\dd s,\dd y) \,\bigg|\,\calf^n_{i-\la_n}\rb,\\
	H^{n,2}_{12}(t)&=\delh\sumt\sum_{q=1}^Q\sum_{j=1}^N\sum_{k=1}^L\bbe\lb \lp\frac\partial{\partial z_{jk}} f(\theta^n_i)-\frac\partial{\partial z_{jk}} f(\theta^{n,q}_i)\rp  \right.\\
	& \quad\times \iint\frac{\Delta^n_{i+k-1} G_{x_j-y}(s)}{\tau_n} \\
	&\quad\times\iint_{(i-\la^{(q-1)}_n)\del}^{s} \si^{(12)}_y((i-\la^{(q-1)}_n)\del,z)\,W'(\dd r,\dd z)\\
	&\quad\times\left.\bone^{n,i}_{\la_n^{(q-1)},\la_n^{(q)}}(s)\vphantom{\lp\frac\partial{\partial z_{jk}} f(\theta^n_i)-\frac\partial{\partial z_{jk}} f(\theta^{n,q}_i)\rp}\,W(\dd s,\dd y) \,\bigg|\,\calf^n_{i-\la_n}\right],\\
	H^{n,2}_{13}(t)&=\delh\sumt\sum_{q=1}^Q\sum_{j=1}^N\sum_{k=1}^L\bbe\lb \frac\partial{\partial z_{jk}} f(\theta^{n,q}_i) \iint\frac{\Delta^n_{i+k-1} G_{x_j-y}(s)}{\tau_n} \right.\\
	& \quad\times  \iint_{(i-\la^{(q-1)}_n)\del}^{s} \si^{(12)}_y((i-\la^{(q-1)}_n)\del,z)\,W'(\dd r,\dd z)\\
	&\quad\times\left.\bone^{n,i}_{\la_n^{(q-1)},\la_n^{(q)}}(s)\vphantom{\frac\partial{\partial z_{jk}}}\,W(\dd s,\dd y) \,\bigg|\,\calf^n_{i-\la_n}\right].
	\end{align*}
	Here, for $r=1,\ldots,Q$, the variable $\theta^{n,r}_i$ is defined in the same way as $\theta^n_i$ in \eqref{thetani} but with all $t^{n,i}_{\la^{(q-1)}_n}$ for $q=r+1,\ldots, Q$ replaced by $t^{n,i}_{\la^{(r-1)}_n}$.
	
	Since $\si^{(12)}_x(s,y)=K(0,x-y)\rho(s,y)$, a standard size argument yields
	\beq\label{cond-a-2}\begin{split} \bbe\left[ \supt |H^{n,2}_{11}(t)|\right] &\lec \sum_{q=1}^Q \del^{-\frac12} (\la^{(q-1)}_n\del)^{\frac12}\del^{\frac\nu2 a^{(q)}}(\la^{(q-1)}_n\del)^{\eps'}\\ &\lec \sum_{q=1}^Q \del^{-\frac12 + (\frac12+\eps')(1-a^{(q-1)})+\frac\nu2 a^{(q)}}.  \end{split}\eeq
	The last term is \eqref{eq-4} with $\frac12+\eps'$ instead of $\ga$, so we can apply the procedure described after \eqref{cond-a} to find $a^{(q)}$ such that \eqref{cond-a-2} goes to $0$ as $n\to\infty$. [In the general case where $\si=\si^{(0)}+\si^{(1)}$, we may assume without loss of generality that $\frac12+\eps'<\ga$. Then the constructed sequence $a^{(q)}$ also satisfies \eqref{cond-a}.]
	
	For $H^{n,2}_{12}(t)$, since $\theta^n_i-\theta^{n,q}_i$ is a term of magnitude $(\la_n^{(q-1)}\del)^{1/2}$ according to Lemma~\ref{mom-ex}, a size estimate yields
	\begin{align*} \bbe\lb\supt |H^{n,2}_{12}(t)|\rb &\lec \sum_{q=1}^{Q}\del^{-\frac12} (\la_n^{(q-1)}\del)^{\frac12}\del^{\frac\nu2 a^{(q)}} (\la_n^{(q-1)}\del)^{\frac12}\\ 
	&\leq \sum_{q=1}^{Q}\del^{-\frac12+(1-a^{(q-1)})+\frac\nu2 a^{(q)}}.  \end{align*}
	This amounts to \eqref{eq-4} with $\ga=1$ and is negligible by our choice of $a^{(q)}$. 
	
	Finally, for $H^{n,2}_{13}(t)$, we compute the conditional expectation in the definition  by first conditioning on $\calf_{i-\la_n^{(q-1)}}^n$. Since $f$ is even, and $\theta^{n,q}_i$ is normally distributed given  $\calf^n_{i-\la^{(q-1)}_n}$, it follows that  $\frac{\partial}{\partial z_{jk}}f(\theta^{n,q}_i)$ belongs to the direct sum of all Wiener chaoses of odd orders. At the same time, the double stochastic integral in $H^{n,2}_{13}(t)$ belongs to the second Wiener chaos; see Proposition~1.1.4 in \citeB{Nualart06-2}. Since Wiener chaoses are mutually orthogonal, conditioning on $\calf^n_{i-\la^{(q-1)}_n}$ shows that $H^{n,2}_{13}(t)=0$.
	\epr

	\bpr[Proof of Lemma~\ref{iter-condexp}] Let $\un\mu_f(m,v) = \un\mu_{f(m+\cdot)}(v)$ for $m\in\bbr^{N\times L}$ and $v\in(\R^{N\times L})^2$, where $\un\mu_f$ is defined after \eqref{muf}. 
	Then
	\begin{align*} \bbe[f(\theta^n_i)\,|\, \calf^n_{i-\la_n}]&=\bbe\left[\bbe\left[f(m^{n,i}_Q)\,\Big|\, \calf^n_{i-\la^{(Q-1)}_n}\right]\,\Big|\, \calf^n_{i-\la_n}\right]\\ &=\bbe\lb \un\mu_f(m^{n,i}_{Q-1},v^{n,i}_{Q-1})\,\Big|\,\calf^n_{i-\la_n} \rb\\ &=\bbe\lb \bbe\lb \un\mu_f(m^{n,i}_{Q-1},v^{n,i}_{Q-1}) \,\Big|\,\calf^n_{i-\la_n^{(Q-2)}}\rb\,\Big|\,\calf^n_{i-\la_n} \rb.
	\end{align*}
	Suppose that we can replace $v^{n,i}_{Q-1}$ by $\bbe[v^{n,i}_{Q-1} \,|\, \calf^n_{i-\la_n^{(Q-2)}}]$ in the second argument of $\un\mu_f$. Then, as $\bbe[\un\mu_f(m+\un Z,v_1)]=\un\mu_f(m,v_1+v_2)$ for all $m\in\bbr^{N\times L}$ and $v_1,v_2\in(\R^{N\times L})^2$, where $\un Z$ is Gaussian with mean $0$ and $\cov(\un Z_{jk},\un Z_{j' k'}) = (v_2)_{j k,j' k'}$, we obtain
	\begin{align*} &\bbe\lb \bbe\lb \un\mu_f\left(m^{n,i}_{Q-1},\bbe\lb v^{n,i}_{Q-1}\,\Big|\,\calf^n_{i-\la_n^{(Q-2)}}\rb \right)\,\Big|\,\calf^n_{i-\la_n^{(Q-2)}}\rb  \,\Big|\,\calf^n_{i-\la_n} \rb\\
	&\quad = \bbe\lb  \un \mu_f\left(m^{n,i}_{Q-2},\bbe\lb v^{n,i}_{Q-2}\,\Big|\,\calf^n_{i-\la_n^{(Q-2)}}\rb \right)  \,\Big|\,\calf^n_{i-\la_n} \rb\\
	&\quad =\bbe\lb\bbe\lb  \un \mu_f\left(m^{n,i}_{Q-2},\bbe\lb v^{n,i}_{Q-2}\,\Big|\,\calf^n_{i-\la_n^{(Q-2)}}\rb \right)  \,\Big|\,\calf^n_{i-\la_n^{(Q-3)}}\rb\,\Big|\,\calf^n_{i-\la_n} \rb.
	\end{align*}
	If we iterate this argument and replace $\bbe[v^{n,i}_{r}\,|\, \calf^n_{i-\la_n^{(r)}}]$ by $\bbe[v^{n,i}_{r}\,|\, \calf^n_{i-\la_n^{(r-1)}}]$ in each step, we get $\bbe[\un\mu_f(m^{n,i}_0,\bbe[v^{n,i}_0\,|\, \calf^n_{i-\la_n^{(0)}}])\,|\, \calf^n_{i-\la_n}] = \un\mu_f(\bbe[v^{n,i}_0\,|\, \calf^n_{i-\la_n}])$, which is the desired expression. 
	
	Hence, it remains to prove that modifying the conditional expectations in the argument above only leads to asymptotically negligible errors, that is, for each $r=1,\dots,Q-1$,
	\beq\label{toshow-2}\begin{split}&\bbe\lb \supt \left|\Delh\sumt \bbe\lb \un\mu_f( m^{n,i}_r, \bbe[v^{n,i}_r \,|\, \calf^n_{i-\la_n^{(r)}} ] )\right.\right.\right.\\
		&\quad-\left.\left.\left.\un\mu_f( m^{n,i}_r, \bbe[v^{n,i}_r\,|\, \calf^n_{i-\la_n^{(r-1)}} ] ) \,\Big|\, \calf^{n}_{i-\la_n}\vphantom{\un\mu_f( m^{n,i}_r, \bbe[v^{n,i}_r \,|\, \calf^n_{i-\la_n^{(r)}} ] )} \rb\vphantom{\sumt}\right|\rb \to 0\end{split}
	\eeq
	as $n\to\infty$.
	For this purpose, we claim that
	$\un\mu_f(\cdot,\cdot)$ is twice continuously differentiable with 
	\beq\label{itoformula}\begin{split}\left|\frac{\partial^2}{\partial v_{\al,\beta}\,\partial v_{\ga,\delta}} \un\mu_f(m,v)\right|&= \frac1{2^{\bone_{\al=\beta}+\bone_{\ga=\delta}}}\lv \un \mu_{\partial_{\al\beta\ga\delta} f}(m,v)\rv\\ &\lec 1+|m|^{p-2}+|v|^{\frac{p-2}{2}}\end{split}\eeq
	for all $\al,\beta,\ga,\delta\in\{1,\ldots,N\}\times\{1,\ldots,L\}$, where $\partial_{\al\beta\ga\delta} = \frac{\partial^4}{\partial z_\al\,\partial z_\beta\,\partial z_\ga\,\partial z_\delta}$.
	
	Indeed, the inequality in \eqref{itoformula} holds by Assumption~\ref{AssCLT}1, while the equation follows from applying the identity 
	\beq\label{id} \frac{\partial}{\partial v_{\al,\beta}} \un\mu_f(m,v) = \frac{1}{2^{\bone_{\al=\beta}}} \un\mu_{\partial_{\al\beta} f}(m,v)\eeq 
	twice. In order to prove this, we may clearly assume that $f$ only depends on the arguments $z_\al$ and $z_\beta$ and that $m=0$. If $\al=\beta$, then for any $h\geq0$, $\un\mu_f(v_{\al,\al}+h)=\bbe[f(Z_\al+B_h)]$ where $Z_\al$ is $N(0, v_{\al,\al})$-distributed and $(B_h)_{h\geq0}$ is an independent standard Brownian motion. Since $f$ is twice continuously differentiable, Itô's formula yields
	\begin{align*} \un\mu_f(v_{\al,\al}+h)&=\bbe[f(Z_\al)]+\frac12 \int_0^h \bbe[f''(Z_\al+B_u)]\,\dd u \\ &= \un\mu_f(v_{\al,\al}) + \frac12 \int_0^h \un\mu_{f''}(v_{\al,\al}+u)\,\dd u, \end{align*}
	which readily gives \eqref{id} for $\al=\beta$. 
	
	If $\al\neq\beta$, then $\un\mu_f\Big(\begin{smallmatrix} v_{\al,\al}+h & v_{\al,\beta}+h\\ v_{\al,\beta}+h &v_{\beta,\beta}+h \end{smallmatrix}\Big)= \bbe[f(Z_\al+B_h,Z'_\beta+B_h)]$ where $Z'_\beta$ is independent of $B$, jointly Gaussian with $Z_\al$, and has mean $0$, variance $v_{\beta,\beta}$, and $\cov(Z_\al,Z'_\beta)=v_{\al,\beta}$. Again with Itô's formula, we get 
	\[ \un\mu_f\Big(\begin{smallmatrix} v_{\al,\al}+h & v_{\al,\beta}+h\\ v_{\al,\beta}+h &v_{\beta,\beta}+h \end{smallmatrix}\Big) = \un\mu_f\Big(\begin{smallmatrix} v_{\al,\al} & v_{\al,\beta}\\ v_{\al,\beta} &v_{\beta,\beta} \end{smallmatrix}\Big) + \frac12 \sum_{i,j\in\{\al,\beta\}} \int_0^h\un\mu_{\partial_{ij}f}\Big(\begin{smallmatrix} v_{\al,\al}+u & v_{\al,\beta}+u\\ v_{\al,\beta}+u &v_{\beta,\beta}+u \end{smallmatrix}\Big)\,\dd u, \]
	which, together with the proven identity for $\al=\beta$, shows \eqref{id}.
	
	Returning to the main line of the proof, Taylor's theorem now yields some $\rho^{n,i}_r$ between $\bbe[v^{n,i}_r\,|\, \calf^n_{i-\la_n^{(r-1)}}]$ and $\bbe[v^{n,i}_r\,|\, \calf^n_{i-\la_n^{(r)}}]$ such that
	\[\begin{split} &\un\mu_f( m^{n,i}_r, \bbe[v^{n,i}_r \,|\, \calf^n_{i-\la_n^{(r)}} ] )-\un\mu_f( m^{n,i}_r, \bbe[v^{n,i}_r\,|\, \calf^n_{i-\la_n^{(r-1)}} ] )\\
	&\quad = \sum_{\al,\beta} \frac{\partial}{\partial v_{\al,\beta}}\un\mu_f( m^{n,i}_r, \bbe[v^{n,i}_r\,|\, \calf^n_{i-\la_n^{(r-1)}} ] )\\
	&\qquad\quad\times\lp \bbe[(v^{n,i}_r)_{\al,\beta} \,|\, \calf^n_{i-\la_n^{(r)}} ]- \bbe[(v^{n,i}_r)_{\al,\beta}\,|\, \calf^n_{i-\la_n^{(r-1)}} ] \rp\\
	&\quad\quad + \frac12 \sum_{\al,\beta,\ga,\delta}\frac{\partial^2}{\partial v_{\al,\beta}\,\partial v_{\ga,\delta}}\un\mu_f( m^{n,i}_r,\rho^{n,i}_r)\\
	&\quad\qquad\times\lp \bbe[(v^{n,i}_r)_{\al,\beta} \,|\, \calf^n_{i-\la_n^{(r)}} ]- \bbe[(v^{n,i}_r)_{\al,\beta}\,|\, \calf^n_{i-\la_n^{(r-1)}} ] \rp\\
	&\quad\qquad\times\lp \bbe[(v^{n,i}_r)_{\ga,\delta} \,|\, \calf^n_{i-\la_n^{(r)}} ]- \bbe[(v^{n,i}_r)_{\ga,\delta}\,|\, \calf^n_{i-\la_n^{(r-1)}} ] \rp. \end{split}\]
	We split the left-hand side of \eqref{toshow-2} into two parts $I^{n,r}_1$ and $I^{n,r}_2$ according to this decomposition. Since 
	\beq \label{condexp-dec} \begin{split} &\bbe\lb \si( t^{n,i}_{\la_n^{(q-1)}},y)\si( t^{n,i}_{\la_n^{(q-1)}},z) \,\Big|\, \calf^n_{i-\la_n^{(r)}}\rb - \bbe\lb \si( t^{n,i}_{\la_n^{(q-1)}},y)\si( t^{n,i}_{\la_n^{(q-1)}},z) \,\Big|\, \calf^n_{i-\la_n^{(r-1)}}\rb\\
		&\quad=  \bbe\lb \si( t^{n,i}_{\la_n^{(q-1)}},y)\si( t^{n,i}_{\la_n^{(q-1)}},z)-\si( t^{n,i}_{\la_n^{(r-1)}},y)\si( t^{n,i}_{\la_n^{(r-1)}},z)\,\Big|\, \calf^n_{i-\la_n^{(r)}}\rb\\
		&\quad\quad-\bbe\lb \si( t^{n,i}_{\la_n^{(q-1)}},y)\si( t^{n,i}_{\la_n^{(q-1)}},z)-\si( t^{n,i}_{\la_n^{(r-1)}},y)\si( t^{n,i}_{\la_n^{(r-1)}},z)\,\Big|\, \calf^n_{i-\la_n^{(r-1)}}\rb, \end{split} \eeq
	we obtain from \eqref{elem}, \eqref{mv}, and a standard size estimate that for small $a$,
	\begin{align*}
	I^{n,r}_2 &\lec \Delh\sumT\bbe\left[\left| \bbe[v^{n,i}_r \,|\, \calf^n_{i-\la_n^{(r)}} ]- \bbe[v^{n,i}_r\,|\, \calf^n_{i-\la_n^{(r-1)}} ] \right|^p\right]^{\frac 2 p}\\
	&\lec \sum_{q=r+1}^{Q} \sum_{j,j'=1}^N \sum_{k,k'=1}^L \Delh\sumT  \left( \iiint  \frac{|\Delta^n_{i+k-1} G_{x_j-y}(s)||\Delta^n_{i+k'-1} G_{x_{j'}-z}(s)|}{\tau^2_n}\right.\\
	&\quad\times   \bbe\lb \left|\si( t^{n,i}_{\la_n^{(q-1)}},y)\si( t^{n,i}_{\la_n^{(q-1)}},z)-\si( t^{n,i}_{\la_n^{(r-1)}},y)\si( t^{n,i}_{\la_n^{(r-1)}},z)\right|^p\rb^{\frac 1 p}\\
	&\quad\times \bone^{n,i}_{\la_n^{(q-1)},\la_n^{(q)}}(s)\,\dd s\left.\vphantom{\frac{|\Delta^n_{i+k-1} G_{x_j-y}(s)||\Delta^n_{i+k'-1} G_{x_{j'}-z}(s)|}{\tau^2_n}}\,\La(\dd y,\dd z) \right)^{2}\\
	&\lec \Del^{-\frac12} ((\la_n^{(r-1)}\Delta_n)^{\frac12})^2 \lec \Delta_n^{\frac 12 -a^{(r-1)}} \leq \del^{\frac12-a} \to 0\as. \end{align*}
	
	Next, we use again Taylor's theorem to write $I^{n,r}_1=I^{n,r}_{11}+I^{n,r}_{12}+I^{n,r}_{13}$ where 
	\beq\label{In1j} I^{n,r}_{1l} = \bbe\lb \supt \left|\Delh\sumt \bbe\lb I^{n,r}_{i,l} \,\Big|\, \calf^{n}_{i-\la_n} \rb\right|\rb, \qquad l=1,2,3,\eeq
	and
	\begin{align*}
	I^{n,r}_{i,1} &=\sum_{\al,\beta}\frac{\partial}{\partial v_{\al,\beta}}\un\mu_f( 0, \bbe[v^{n,i}_r\,|\, \calf^n_{i-\la_n^{(r-1)}} ] )\\
	&\quad\times\lp \bbe[(v^{n,i}_r)_{\al,\beta} \,|\, \calf^n_{i-\la_n^{(r)}} ]- \bbe[(v^{n,i}_r)_{\al,\beta}\,|\, \calf^n_{i-\la_n^{(r-1)}} ] \rp,\\
	I^{n,r}_{i,2} &= \sum_{\al,\beta,\ga}\frac{\partial^2}{\partial m_\ga\, \partial v_{\al,\beta}}\un\mu_f( 0, \bbe[v^{n,i}_r\,|\, \calf^n_{i-\la_n^{(r-1)}} ] )\\ 
	&\quad\times\lp \bbe[(v^{n,i}_r)_{\al,\beta} \,|\, \calf^n_{i-\la_n^{(r)}} ]- \bbe[(v^{n,i}_r)_{\al,\beta}\,|\, \calf^n_{i-\la_n^{(r-1)}} ] \rp (m^{n,i}_r)_\ga,\\
	I^{n,r}_{i,3} &= \frac12 \sum_{\al,\beta,\ga,\delta} \frac{\partial^3}{\partial m_\ga\,\partial m_\delta\, \partial v_{\al,\beta}}\un\mu_f(\wh \rho^{n,i}_r, \bbe[v^{n,i}_r\,|\, \calf^n_{i-\la_n^{(r-1)}} ] )\\
	&\quad\times\lp \bbe[(v^{n,i}_r)_{\al,\beta} \,|\, \calf^n_{i-\la_n^{(r)}} ]- \bbe[(v^{n,i}_r)_{\al,\beta}\,|\, \calf^n_{i-\la_n^{(r-1)}} ] \rp (m^{n,i}_r)_\ga(m^{n,i}_r)_\delta,
	\end{align*}
	with some $\wh \rho^{n,i}_r$ between $0$ and $m^{n,i}_{r}$. 
	
	Because $\frac{\partial}{\partial v_{\al,\beta}}\un\mu_f( 0, \bbe[v^{n,i}_r\,|\, \calf^n_{i-\la_n^{(r-1)}} ])$ is $\calf^n_{i-\la_n^{(r-1)}}$-measurable and the difference $\bbe[v^{n,i}_r \,|\, \calf^n_{i-\la_n^{(r)}} ]- \bbe[v^{n,i}_r\,|\, \calf^n_{i-\la_n^{(r-1)}} ]$ has a vanishing conditional expectation given $\calf^n_{i-\la_n^{(r-1)}}$, it follows that $I^{n,r}_{11} = 0$. Next, as in \eqref{itoformula},
	\beq\label{itoformula-2}\begin{split} \frac{\partial^2}{\partial m_\ga\, \partial v_{\al,\beta}}\un\mu_f&=\frac1{2^{\bone_{\al=\beta}}} \un\mu_{\partial_{\al\beta\ga} f},\\ 
		\frac{\partial^3}{\partial m_\ga\,\partial m_\delta\, \partial v_{\al,\beta}}\un\mu_f&= \frac1{2^{\bone_{\al=\beta}}}\un\mu_{\partial_{\al\beta\ga\delta}f} \end{split}\eeq
	by Itô's formula.
	Since $f$ has odd third derivatives, it follows that the first term in \eqref{itoformula-2} vanishes at $m=0$, whence $I^{n,r}_{12}=0$.
	For $I^{n,r}_{13}$, the second relation in \eqref{itoformula-2}, Assumption~\ref{AssCLT}1, \eqref{elem}, \eqref{condexp-dec}, and a size estimate give
	\beq\label{In13}\begin{split}
		\qquad I^{n,r}_{13} &\lec \del^{-\frac12}  \bbe\lb \left|\bbe[v^{n,i}_r \,|\, \calf^n_{i-\la_n^{(r)}} ]- \bbe[v^{n,i}_r\,|\, \calf^n_{i-\la_n^{(r-1)}} ]\right|^p\rb^{\frac 1 p}(\bbe[|m^{n,i}_r|^{2p}]^{\frac 1 {2p}})^2\\
		&\lec \Del^{-\frac12} (\la_n^{(r-1)}\Delta_n)^{\frac12} \Del^{ \nu  a^{(r)}} = \del^{-\frac12 a^{(r-1)}+\nu a^{(r)}}.
	\end{split}\raisetag{-2.5\baselineskip}\eeq
	Without loss of generality, we may assume that $\eps'>0$ from \eqref{Holder-2} is small enough such that $1+2\eps'<\nu$. Now consider again the recurrence relation \eqref{recurrence} with $\ga$ replaced by $\frac12 +\eps'$. From \eqref{recurrence}, we obtain for any $q=1,\ldots,Q$, the equivalence 
	\beq\label{bQ-2} -\frac12 b^{(q-1)}+\nu b^{(q)}>0\iff b^{(q)}>\frac{2\eps'}{\nu(1+4\eps')}.  \eeq
	If this is true for all $q=1,\ldots, Q-1$, we can ensure that \eqref{In13} converges to $0$ for all $r=1,\ldots,Q-1$ by choosing $a^{(q)}$, $q=1,\ldots,Q-1$, sufficiently close to $b^{(q)}$, $q=1,\ldots,Q-1$. 
	
	As $b^{(q)}$ is decreasing in $q$, we only need to check the right-hand side of \eqref{bQ-2} for $q=Q-1$. Since $b^{(Q)}=b<0$, by \eqref{recurrence}, this is equivalent to
	\beq\label{bQ-3} \frac{2\eps'+\nu b}{1+2\eps'} > \frac{2\eps'}{\nu(1+4\eps')} \iff  |b|< \frac{2\eps'}{\nu^2}\lp \nu-\frac{1+2\eps'}{1+4\eps'} \rp. \eeq
	Since $\nu>1$, the last term is strictly positive. So by choosing $b$ according to the right-hand side of \eqref{bQ-3}, we can find $a^{(1)},\ldots,a^{(Q-1)}$ such that the right-hand side of \eqref{In13} tends to $0$ for all $r=1,\ldots,Q-1$.
	\epr

	\bpr[Proof of Lemma~\ref{remove-condexp}] By Taylor's theorem, we can decompose the difference $\un\mu_f(\bbe[v^{n,i}_0 \,|\, \calf^n_{i-\la_n}]) - \un\mu_f(v^{n,i}_0) = J^{n,i}_1+J^{n,i}_2$ with
	\begin{align*}
	J^{n,i}_1 &= - \sum_{\al,\beta} \frac{\partial}{\partial {v_{\al,\beta}}}\un\mu_f(\bbe[v^{n,i}_0 \,|\, \calf^n_{i-\la_n}]) ((v^{n,i}_0)_{\al,\beta}-\bbe[(v^{n,i}_0)_{\al,\beta} \,|\, \calf^n_{i-\la_n}]),\\
	J^{n,i}_2&= -\frac12\sum_{\al,\beta,\ga,\delta} \frac{\partial^2}{\partial v_{\al,\beta}\,\partial v_{\ga,\delta}}\un\mu_f(\wt\rho^n_i)((v^{n,i}_0)_{\al,\beta}-\bbe[(v^{n,i}_0)_{\al,\beta} \,|\, \calf^n_{i-\la_n}])\\
	&\quad\times((v^{n,i}_0)_{\ga,\delta}-\bbe[(v^{n,i}_0)_{\ga,\delta} \,|\, \calf^n_{i-\la_n}]),
	\end{align*}
	and some $\wt \rho^n_i$ between $v^{n,i}_0$ and $\bbe[v^{n,i}_0 \,|\, \calf^n_{i-\la_n}]$. 
	
	Observe that $J^{n,i}_1$ is $\calf^n_{i}$-measurable with a vanishing conditional expectation given $\calf^{n}_{i-\la_n}$. Hence, we can apply a martingale argument as in the proofs of Lemma~\ref{knknprime} or Lemma~\ref{approx}. Since $\bbe[|v^{n,i}_0-\bbe[v^{n,i}_0 \,|\, \calf^n_{i-\la_n}]|^2]^{1/2} \lec (\la_n\del)^{1/2}$ by \eqref{condexp-dec} and Lemma~\ref{mom-ex}, we obtain for sufficiently small $a$,
	\[ \bbe\lb \supt \lv \delh \sumt J^{n,i}_1\rv\rb\lec \sqrt{\la_n}(\la_n\del)^{\frac12}\lec \del^{\frac12 -a}\to0\as. \]
	Similarly, by a size argument, the left-hand side of the previous display with $J^{n,i}_2$ instead of $J^{n,i}_1$ is bounded by a multiple of $\del^{-1/2} ((\la_n\del)^{1/2})^2\lec \del^{1/2-a}$,
	which tends to $0$ as before and completes the proof the lemma.
	\epr

	\bpr[Proof of Lemma~\ref{remove-s}] By definition, $\mu_f(\un\si^2(s,x))=\un\mu_f(v_x(s))$ where $v_x(s)_{jk,j'k'}=\Ga_{|k-k'|}\si^2(s,x_j)\bone_{j=j'}$. Writing $t^{n,i}_1=(i-1)\Del$, we use Taylor's theorem to write the difference in braces in Lemma~\ref{remove-s} as the sum of
	\begin{align*}
	K^{n,i}_1&=\sum_{j,j'=1}^N\sum_{k,k'=1}^L \sum_{q=1}^Q \frac{\partial}{\partial v_{jk,j'k'}}\un\mu_f(v_x(t^{n,i}_1)) \iiint  \frac{\Delta^n_{i+k-1} G_{x_j-y}(s)}{\tau_n} \\
	&\quad\times\frac{\Delta^n_{i+k'-1} G_{x_{j'}-z}(s)}{\tau_n}\lp \si( t^{n,i}_{\la_n^{(q-1)}},y)\si( t^{n,i}_{\la_n^{(q-1)}},z)-\si(t^{n,i}_1,y)\si(t^{n,i}_1,z)\rp\\
	&\quad\times\bone^{n,i}_{\la_n^{(q-1)},\la_n^{(q)}}(s)\,\dd s\,\La(\dd y,\dd z),\\
	K^{n,i}_2&=\sum_{j=1}^N\sum_{k,k'=1}^L \frac{\partial}{\partial v_{jk,jk'}}\un\mu_f(v_x(t^{n,i}_1)) \iiint  \frac{\Delta^n_{i+k-1} G_{x_j-y}(s)\Delta^n_{i+k'-1} G_{x_{j}-z}(s)}{\tau_n^2} \\
	&\quad\times\left( \si(t^{n,i}_1,y)\si(t^{n,i}_1,z)-\si^2(t^{n,i}_1,x_j)\right)\bone^{n,i}_{\la_n,0}(s)\,\dd s\,\La(\dd y,\dd z),\\
	K^{n,i}_3&=\sum_{j\neq j'}\sum_{k,k'=1}^L \frac{\partial}{\partial v_{jk,j'k'}}\un\mu_f(v_x(t^{n,i}_1)) \iiint  \frac{\Delta^n_{i+k-1} G_{x_j-y}(s)}{\tau_n} \\
	&\quad\times \frac{\Delta^n_{i+k'-1} G_{x_{j'}-z}(s)}{\tau_n}\si(t^{n,i}_1,y)\si(t^{n,i}_1,z)\bone^{n,i}_{\la_n,0}(s)\,\dd s\,\La(\dd y,\dd z),\\
	K^{n,i}_4&=\sum_{j=1}^N\sum_{k,k'=1}^L  \frac{\partial}{\partial v_{jk,jk'}}\un\mu_f(v_x(t^{n,i}_1)) \si^2(t^{n,i}_1,x_j) \lb (\Ga^n_{|k-k'|}-\vphantom{\Pi^n_{|k-k'|,0}}\Ga_{|k-k'|})\right.\\
	&\quad+\left.\left(\Pi^n_{|k-k'|,0}((0, \la_n\del)\times\R^d\times\R^d)-\Pi^n_{|k-k'|,0}([0,\infty)\times\R^d\times\R^d)\right)\rb,\\
	K^{n,i}_5&=\frac12\sum_{j_1,j_2=1}^N \sum_{k_1,k_1',k_2,k_2'=1}^L \frac{\partial^2}{\partial v_{j_1k_1,j_1k_1'}\,\partial v_{j_2k_2,j_2k_2'}}\un\mu_f(\chi^n_i)\\
	&\quad\times\left((v^{n,i}_0)_{j_1k_1,j_1k_1'}-\Ga_{|k_1-k_1'|}\si^2(s,x_{j_1})\right)\\
	&\quad\times\left((v^{n,i}_0)_{j_2k_2,j_2k_2'}-\Ga_{|k_2-k_2'|}\si^2(s,x_{j_2})\right),\\
	K^{n,i}_6&=\frac12\sum_{j_1 \neq j_1'\text{ or } j_2\neq j_2'} \sum_{k_1,k_1',k_2,k_2'=1}^L \frac{\partial^2}{\partial v_{j_1k_1,j_1'k_1'}\,\partial v_{j_2k_2,j_2'k_2'}}\un\mu_f(\chi^n_i)\\
	&\quad\times\left((v^{n,i}_0)_{j_1k_1,j_1'k_1'}-\Ga_{|k_1-k_1'|}\si^2(s,x_{j_1})\bone_{j_1=j_1'}\right)\\
	&\quad\times\left((v^{n,i}_0)_{j_2k_2,j_2k_2'}-\Ga_{|k_2-k_2'|}\si^2(s,x_{j_2})\bone_{j_2=j_2'}\right),
	\end{align*}
	with some $\chi^n_i$ between $v^{n,i}_0$ and $v_x(t^{n,i}_1)$. 
	
	If we take $r=t^{n,i}_1$ and $u=t^{n,i}_{\la_n^{(q-1)}}$ in the identity
	\begin{align*}
	&\si(r,y)\si(r,z)-\si( u,y)\si( u,z)\\
	&\quad=\si( u,z)\lb (\si^{(0)}(r,y)-\si^{(0)}( u,y))+(\si^{(1)}(r,y)-\si^{(1)}( u,y))  \rb\\
	&  \qquad + \si( u,y)\lb (\si^{(0)}(r,z)-\si^{(0)}( u,z))+(\si^{(1)}(r,z)-\si^{(1)}( u,z))  \rb\\
	&\qquad + \lb (\si^{(0)}(r,y)-\si^{(0)}( u,y))+(\si^{(1)}(r,y)-\si^{(1)}( u,y))\rb\\
	&\qquad\quad\times \lb(\si^{(0)}(r,z)-\si^{(0)}(u,z))+(\si^{(1)}(r,z)-\si^{(1)}(u,z))\rb,
	\end{align*}
	with $\si^{(1)}=\si-\si^{(0)}$ and $\si^{(0)}$ as in Assumption~\ref{AssCLT}3,
	the difference in parentheses in $K^{n,i}_1$ splits into a sum of eight products. By a simple size argument, each product with a difference involving $\si^{(0)}$ at different time points as a factor contributes to $\Del^{1/2}\sumt K^{n,i}_1$ by a term of size at most $\sum_{q=1}^Q \del^{-1/2}(\la_n^{(q-1)}\del)^\ga\del^{\nu  a^{(q)}}$, which is negligible by \eqref{cond-a}. We may also discard the product of two differences of $\si^{(1)}$ as this has a total contribution of size $ \sum_{q=1}^Q \del^{-1/2}((\la_n^{(q-1)}\del)^{1/2})^2\lec \del^{1/2-a}$, which vanishes for small $a$. 
	
	So it remains to consider the two terms
	$\si( u,z) (\si^{(1)}(r,y)-\si^{(1)}( u,y))$ and
	$\si( u,y)(\si^{(1)}(r,z)-\si^{(1)}( u,z))$.
	By symmetry, it suffices to analyze the first one, whose contribution to $\Del^{1/2}\sumt K^{n,i}_1$ is
	\begin{align*}
	&-\delh\sumt\sum_{j,j'=1}^N\sum_{k,k'=1}^L \frac{\partial}{\partial v_{jk,j'k'}}\un\mu_f(v_x(t^{n,i}_1)) \iiint  \frac{\Delta^n_{i+k-1} G_{x_j-y}(s)}{\tau_n} \\
	&\quad\times \frac{\Delta^n_{i+k'-1} G_{x_{j'}-z}(s)}{\tau_n}\sum_{q=1}^Q \si( t^{n,i}_{\la_n^{(q-1)}},z) (\si^{(1)}(t^{n,i}_1,y)-\si^{(1)}( t^{n,i}_{\la_n^{(q-1)}},y))\\
	&\quad\times\bone^{n,i}_{\la_n^{(q-1)},\la_n^{(q)}}(s)\,\dd s\,\La(\dd y,\dd z).
	\end{align*}
	Now if we replace $t^{n,i}_1 $ in the argument of $v_x$ by $t^{n,i}_{\la_n^{(q-1)}}$, the error induced by this modification is $\lec\sum_{q=1}^Q \del^{-1/2}((\la^{(q-1)}_n\del)^{1/2})^2\del^{\nu a^{(q)}}\lec \del^{1/2-a}$, which vanishes for small $a$. The resulting term is $-K^{n}_{11}(t)-K^{n}_{12}(t)$ where
	\begin{align*}
	K^{n}_{11}(t)&=\delh\sumt \sum_{q=1}^Q\sum_{j,j'=1}^N\sum_{k,k'=1}^L \frac{\partial}{\partial v_{jk,j'k'}}\un\mu_f(v_x(t^{n,i}_{\la_n^{(q-1)}}))  \\
	&\quad\times\iiint  \frac{\Delta^n_{i+k-1} G_{x_j-y}(s)\Delta^n_{i+k'-1} G_{x_{j'}-z}(s)}{\tau_n^2}\si( t^{n,i}_{\la_n^{(q-1)}},z)\\
	&\quad\times \left\{ \si^{(1)}((i-1)\del,y)-\si^{(1)}( t^{n,i}_{\la_n^{(q-1)}},y) \vphantom{\bbe\lb \si^{(1)}((i-1)\del,y)-\si^{(1)}( t^{n,i}_{\la_n^{(q-1)}},y)\,\Big|\, \calf^n_{i-\la_n^{(q-1)}}\rb}\right.\\
	&\qquad- \left.\bbe\lb \si^{(1)}((i-1)\del,y)-\si^{(1)}( t^{n,i}_{\la_n^{(q-1)}},y)\,\Big|\, \calf^n_{i-\la_n^{(q-1)}}\rb\right\}\\
	&\quad\times \bone^{n,i}_{\la_n^{(q-1)},\la_n^{(q)}}(s)\,\dd s\,\La(\dd y,\dd z),\\
	K^{n}_{12}(t)&=\delh\sumt \sum_{q=1}^Q\sum_{j,j'=1}^N\sum_{k,k'=1}^L \frac{\partial}{\partial v_{jk,j'k'}}\un\mu_f(v_x(t^{n,i}_{\la_n^{(q-1)}}))  \\
	&\quad\times\iiint  \frac{\Delta^n_{i+k-1} G_{x_j-y}(s)\Delta^n_{i+k'-1} G_{x_{j'}-z}(s)}{\tau_n^2}\si( t^{n,i}_{\la_n^{(q-1)}},z)\\
	&\quad\times\bbe\lb \si^{(1)}((i-1)\del,y)-\si^{(1)}( t^{n,i}_{\la_n^{(q-1)}},y)\,\Big|\, \calf^n_{i-\la_n^{(q-1)}}\rb\\
	&\quad\times\bone^{n,i}_{\la_n^{(q-1)},\la_n^{(q)}}(s) \,\dd s\,\La(\dd y,\dd z).
	\end{align*}
	The $i$th term in the sum defining $K^n_{11}(t)$ is $\calf^n_{i-1}$-measurable and has zero conditional expectation given $\calf^n_{i-\la_n^{(q-1)}}$. By a martingale argument (see the proofs of Lemma~\ref{knknprime} and Lemma~\ref{approx}), we derive
	\[ \bbe\lb \supt |K^n_{11}(t)|\rb \lec \sum_{q=1}^Q \sqrt{\la_n^{(q-1)}}(\la_n^{(q-1)}\del)^{\frac12}\lec \del^{\frac12-a} \to0 \]
	as $n\to\infty$.
	
	Furthermore, from  \eqref{si-fixed-x}, \eqref{si11-12} and Assumption~\ref{AssCLT}3, we deduce 
	\beq\label{cond-exp-rate} \begin{split}\bbe\left[\left|\bbe[\si^{(1)}(s,y)-\si^{(1)}(r,y)\,|\, \calf_r]\right|^p\right]^{\frac1 p} &=\bbe\lb\lv \int_r^s \si^{(11)}_x(u)\,\dd u\rv^p\rb^{\frac1p} \\ &\lec s-r\end{split}
	\eeq for $s\geq r$. Thus, 
	\[ \bbe\lb \supt |K^{n}_{12}(t)|\rb \lec \sum_{q=1}^Q\del^{-\frac12}(\la_n^{(q-1)}\del)\lec \del^{\frac12-a} \to0 \]
	as $n\to\infty$.
	
	Next, by a change of variables, the triple integral in the definition of $K^{n,i}_2$ actually equals
	\begin{align*} &\iiint_0^\infty  \left( \si((i-1)\del,x_j-y)\si((i-1)\del,x_j-z)-\si^2((i-1)\del,x_j)\right)\\
	&\quad\times\bone_{s<(\la_n+(k\vee k'))\del}\,\Pi^n_{|k'-k|,0}(\dd s,\dd y,\dd z). \end{align*}
	Consequently, using Taylor's theorem and the second part of Lemma~\ref{takesform}, we obtain $$\delh\sumt K^{n,i}_2 =K^n_{21}(t)+K^n_{22}(t)$$ where, with a slight abuse of notation ($x_j$, $j=1,\ldots,N$, are the observation points in $\R^d$, while $x_l$ and $x_m$ are the $l$th and $m$th coordinate of a generic point $x\in\R^d$, respectively), 
	\begin{align*} K^n_{21}(t)&=-\delh\sumt\sum_{j=1}^N\sum_{k,k'=1}^L \frac{\partial}{\partial v_{jk,jk'}}\un\mu_f(v_x(t^{n,i}_1))  \sum_{l=1}^d \si(t^{n,i}_1,x_j)\frac{\partial}{\partial x_l}\si(t^{n,i}_1,x_j)\\
	&\quad\times\iiint_0^{(\la_n+(k\vee k'))\del} (y_l+z_l) \,\Pi^n_{|k'-k|,0}(\dd s,\dd y,\dd z), \\
	K^n_{22}(t)&= \frac{\delh}{2}\sumt \sum_{j=1}^N\sum_{k,k'=1}^L \sum_{l,m=1}^d \frac{\partial}{\partial v_{jk,jk'}}\un\mu_f(v_x(t^{n,i}_1))\\
	&\quad\times \iiint_0^{(\la_n+(k\vee k'))\del}  \bigg\{ \si(t^{n,i}_1,\ov\chi^{n,i,j}_z)\frac{\partial^2}{\partial x_l\,\partial x_m} \si(t^{n,i}_1,\wt \chi^{n,i,j}_y) y_ly_m\\
	&\quad \quad + 2\frac{\partial}{\partial x_l}\si(t^{n,i}_1,\wt \chi^{n,i,j}_y)\frac{\partial}{\partial x_m} \si(t^{n,i}_1,\ov\chi^{n,i,j}_z)y_lz_m\\
	&\qquad+\si(t^{n,i}_1,\wt \chi^{n,i,j}_y)\frac{\partial^2}{\partial x_l\,\partial x_m} \si(t^{n,i}_1,\ov\chi^{n,i,j}_z) z_l z_m
	\bigg\}\, \Pi^n_{|k'-k|,0}(\dd s,\dd y,\dd z),  \end{align*}
	and $\wt \chi^{n,i,j}_y$ (resp., $\ov\chi^{n,i,j}_z$) is some point on the line between $x_j$ and $y$ (resp., $x_j$ and $z$). By the symmetry properties of $\Pi^n_{|k'-k|,0}$, one has 
	\beq\label{important} \iiint_0^{(\la_n+(k\vee k'))\del} (y_l+z_l)\,\Pi^n_{|k'-k|,0}(\dd s,\dd y,\dd z)=0,\eeq and we immediately obtain $K^n_{21}(t)=0$.
	For $K^n_{22}(t)$, a size argument together with \eqref{si-der} and the second part of Lemma~\ref{takesform} yields (with obvious modifications if $d=1$ and $F=\delta_0$)
	\begin{align*} \bbe\lb \supt |K^n_{22}(t)|\rb &\lec \del^{-\frac12} \sum_{k,k'=1}^L \iiint_0^\infty (|y|^2+|z|^2)\,|\Pi^n_{k'-k,0}|(\dd s,\dd y,\dd z),
	\end{align*}
	which converges to $0$ by Lemma~\ref{y2}.
	
	Next, observe that $K^{n,i}_3$ only contains derivatives of $\un\mu_f$ with $j\neq j'$. Since Assumption~\ref{AssCLT}2 implies that each coordinate of $f$ only depends on the observations at one spatial point if $\al\leq 1$ and $F\neq\delta_0$, we have $K^{n,i}_3=0$ in these cases. In other words, for $K^{n,i}_3$, we only need to consider the case $\al>1$ or the case $d=1$ and $F=\delta_0$. By a size argument and \eqref{speed}, it follows that
	\begin{align*} \bbe\lb \supt \lv \delh\sumt K^{n,i}_3\rv\rb &\lec \del^{-\frac12}\sum_{j\neq j'}\sum_{k,k'=1}^L |\Pi^n_{|k'-k|,|x_j-x_{j'}|}|([0,\infty)\times\R^d\times\R^d) \\ &\lec \begin{cases} \del^{\frac{\al-1}{2}}|\log\del|,& \al>1,\\ \del|\log\del|,&F=\delta_0,  \end{cases} \end{align*}
	which vanishes in both cases.
	
	For the term $K^{n,i}_4$, a  size estimate suffices. Indeed, using $a>\frac1{2\nu}$, \eqref{Ganr-speed}, and Lemma~\ref{nu}, 
	\[ \bbe\lb\supt \lv\delh\sumt K^{n,i}_4\rv\rb \lec \del^{-\frac12}(\del^{\nu a}+\del) \to0\as. \]
	Finally, by the calculations so far and \eqref{speed}, one can quickly check that \begin{align*}&\bbe[|(v^{n,i}_0)_{jk,jk'}-\Ga_{|k'-k|}\si^2((i-1)\del,x_j)|^p]^{\frac1 p}\\
	&\quad+\bbe[|(v^{n,i}_0)_{jk,j'k'}|^p]^{\frac1p}\bone_{j\neq j'}\bone_{\al>1\text{ or }F=\delta_0}\lec (\la_n\del)^{\frac12}+\del^{\frac12}\lec (\la_n\del)^{\frac12},\end{align*}
	which implies that also $\Del^{1/2}\sumt (K^{n,i}_5+K^{n,i}_6)$ is negligible if $a$ is close to $\frac 1 {2\nu}$ because this term is of size at most $\del^{-1/2}((\la_n\del)^{1/2})^2 \lec \del^{1/2-a}$.
	\epr
	
	\brem\label{rem-discr-2} An important step in the proof of Lemma~\ref{remove-s} (and in fact, of Theorem~\ref{CLT}) is that the integral in \eqref{important} vanishes. This is the only place where we use the spatial differentiability properties of $\si$. If $\si$ were, for example, only Lipschitz continuous in space, we would obtain an integral as in \eqref{important}, but with absolute values of $y_l$ and $z_l$ instead. We are not able to show (nor do we believe) that the resulting term is $o(\sqrt{\Del})$. This is also why we prove the underlying central limit theorem (i.e., Proposition~\ref{CLT-core}) \emph{before} we approximate the volatility process $\si$ in space. Indeed, if in \eqref{important}, we still had a stochastic instead of a deterministic integral, then the resulting expression would not vanish even if $\si$ were differentiable in space.
	\erem
	
	\bpr[Proof of Lemma~\ref{lebesgue}] 
	The expression on the left-hand side of $\limL$ equals $-\sum_{i=1}^3 L^n_i(t)$
	where
	\begin{align*}
	L^n_1(t)&=\del^{-\frac12}\int_{\del t^*_n}^t \mu_f(\un\si^2(s,x))\,\dd s,\\
	L^{n}_2(t)&=\del^{-\frac12}\sumt \sum_{j=1}^N \frac{\partial}{\partial v_j}\mu_f(\un\si^2((i-1)\del,x))\\
	&\quad\times \int_{(i-1)\del}^{i\del} (\si^2(s,x_j)-\si^2((i-1)\del,x_j))\,\dd s,\\
	L^{n}_3(t)&=\del^{-\frac12}\sumt\sum_{j,j'=1}^N\int_{(i-1)\del}^{i\del}  \frac{\partial^2}{\partial v_j\,\partial v_{j'}}\mu_f(\om^{n,i}(s))\\
	&\quad \times(\si^2(s,x_j)-\si^2((i-1)\del,x_j))(\si^2(s,x_{j'})-\si^2((i-1)\del,x_{j'}))\,\dd s,
	\end{align*}
	with some $\om^{n,i}(s)$ between $\un\si^2(s,x)$ and $\un\si^2((i-1)\del,x)$.
	As $|t-\del t^*_n|\leq L\del$, it follows that  $\bbe[\supt|L^n_1(t)|]\lec\Del^{1/2}\to0$ as $n\to\infty$. Similarly, $L^n_3(t)$ is of order $\del^{-1/2}\del^{-1}\del(\del^{1/2})^2\leq\Del^{1/2}$ by a size argument. 
	
	Using the identity $a^2-b^2 = 2b(a-b)+(a-b)^2$, we can rewrite $L^n_2(t)$ as
	\begin{align*}
	&\del^{-\frac12}\sumt \sum_{j=1}^N \frac{\partial}{\partial v_j}\mu_f(\un\si^2((i-1)\del,x))\\
	&\quad\times \Bigg\{ 2\int_{(i-1)\del}^{i\del} \si((i-1)\del,x_j)(\si(s,x_j)-\si((i-1)\del,x_j))\,\dd s\\
	&\qquad + \frac{\partial}{\partial v_j}\mu_f(\un\si^2((i-1)\del,x)) \int_{(i-1)\del}^{i\del} (\si(s,x_j)-\si((i-1)\del,x_j))^2\,\dd s\Bigg\}.
	\end{align*}
	The second part is of magnitude $\lec\del^{-1/2}\del^{-1}\del(\Del^{1/2})^2=\Del^{1/2}$. The first part equals $L^n_{21}(t)+L^n_{22}(t)$ where
	\begin{align*} L^n_{21}(t)&=2\del^{-\frac12}\sumt \sum_{j=1}^N \frac{\partial}{\partial v_j}\mu_f(\un\si^2((i-1)\del,x)) \int_{(i-1)\del}^{i\del} \si((i-1)\del,x_j)\\
	&\quad\times\Big\{\si(s,x_j)-\si((i-1)\del,x_j)\\
	&\qquad-\bbe[\si(s,x_j)-\si((i-1)\del,x_j)\,|\,\calf^n_{i-1}]\Big\}\,\dd s,\\
	L^n_{22}(t)&=2\del^{-\frac12}\sumt \sum_{j=1}^N \frac{\partial}{\partial v_j}\mu_f(\un\si^2((i-1)\del,x)) \int_{(i-1)\del}^{i\del} \si((i-1)\del,x_j)\\
	&\quad \times\bbe[\si(s,x_j)-\si((i-1)\del,x_j)\,|\,\calf^n_{i-1}]\,\dd s. \end{align*}
	One readily checks that $L^n_{21}(t)$ is a martingale of size $\lec\del^{1/2}$. For $L^n_{22}(t)$, recall from  \eqref{cond-exp-rate} that the size of the conditional expectation is bounded by $\del^{\ga}+\del$, so  $L^n_{22}(t)$ is of size $\del^{\ga-1/2}\to0$. 
	\epr
	
	\section{Proofs for Sections~2.2 and 2.3}\label{proof-rest}
	\bpr[Proof of Corollary~\ref{multipower-result}]
	Under the stated hypotheses of part (1) [resp., part (2)], $\Phi$ and $\Psi$ satisfy Assumption~\ref{AssLLN}1 (resp., Assumptions~\ref{AssCLT}1 and \ref{AssCLT}2) with some $p'$ slightly larger than $\ov w$ (resp., with $p=\ov w$). It is clear that Assumptions~\ref{AssLLN}2 and \ref{AssLLN}3 would still hold with $p'$ if it is close enough to $\ov w$. Then the corollary readily follows from Theorems~\ref{LLN} and \ref{CLT}. 
	\epr
	
	\bpr[Proof of Theorem~\ref{known-alpha}] Since $\wt\tau_n/\tau_n\to 1$ by \eqref{taun-asymp}, \eqref{LLN-mpv-2} follows from Corollary~\ref{multipower-result} (1). By the same reason, and using \eqref{stable}, we obtain \eqref{stud-1-2} from Corollary~\ref{multipower-result} (2) by studentization.
	\epr
	
	\bpr[Proof of Theorem~\ref{est-alpha}] Suppose first that $N=1$ and write $x=x_1$. Furthermore, define $\Theta \colon \R^2\to \R^2$ by $\Theta_1(x_1,x_2)=\Phi(p;x_1)=|x_1|^p$ and $\Theta_2(x_1,x_2)=\Phi^{(2)}(p;x_1,x_2)=|x_1+x_2|^p$, which clearly satisfy Assumption~\ref{AssLLN}1 and Assumptions~\ref{AssCLT}1 and \ref{AssCLT}2 under the hypotheses of part (1) and (2), respectively. Thus, \eqref{alpha-conv} follows from Theorem~\ref{LLN} upon realizing that $\mu_{\Phi(p;\cdot)}=\bbe[|Z|^p]$ and $\mu_{\Phi^{(2)}(p;\cdot)}=\bbe[|Z|^p](2+2\Ga_1)^{p/2} = \bbe[|Z|^p]2^{(1-\al/ 2) p/ 2}$, where $Z$ is a standard normal variable. 
	
	Moreover, Theorem~\ref{CLT} implies that 
	\[ \begin{pmatrix} \del^{-\frac12} \lp V^n_{\Theta_1}(Y_x,t)-\mu_{\Phi(p;\cdot)} \int_0^t |\si(s,x)|^p\,\dd s\rp\\ \del^{-\frac12} \lp V^n_{\Theta_2}(Y_x,t)-\mu_{\Phi^{(2)}(p;\cdot)}\int_0^t |\si(s,x)|^p\,\dd s\rp \end{pmatrix} \limst \calz=\begin{pmatrix} \calz_1\\ \calz_2\end{pmatrix}, \]
	where $\calz$ is as in Theorem~\ref{CLT} with $\calc(t)=\calc  \int_0^t |\si(s,x)|^{2p}\,\dd s$ and $\calc = \lp\begin{smallmatrix} \calc_{11}&\calc_{12}\\ \calc_{12} & \calc_{22}\end{smallmatrix}\rp$ as defined in \eqref{Cs}.
	Now let $H(a,b)=2-\frac{4}{p}\log_2(\frac{b}{a})$ and observe that 
	\begin{align*} \wh \al^{(p)}_n &= H(V^n_{\Theta_1}(Y_x,T),V^n_{\Theta_2}(Y_x,T)),\\ \al &= H(V_{\Theta_1}(Y_x,T),V_{\Theta_2}(Y_x,T)).\end{align*} Therefore, by the mean value theorem,
	\beq \begin{split}&\del^{-\frac12}(\wh\al^{(p)}_n-\al)\\&\quad= \lp \frac{\partial}{\partial a} H ,\frac{\partial}{\partial b} H\rp (v^n_1,v^n_2) \begin{pmatrix} \del^{-\frac12} \lp V^n_{\Theta_1}(Y_x,t)-V_{\Theta_1}(Y_x,T)\rp\\ \del^{-\frac12} \lp V^n_{\Theta_2}(Y_x,t)-V_{\Theta_2}(Y_x,T)\rp \end{pmatrix}, \end{split} \eeq
	with some $(v^n_1,v^n_2)$ satisfying  $(v^n_1,v^n_2) \limL V_{\Theta}(Y_x,T)$ by Theorem~\ref{LLN}. Hence, by \eqref{stable},
	\begin{align*} &\del^{-\frac12}(\wh\al^{(p)}_n-\al)\\ &\quad\stackrel{\mathrm{st}-\call}{\longrightarrow} N\lp 0,\lp \frac{\partial}{\partial a} H ,\frac{\partial}{\partial b} H\rp (V_\Theta(Y_x,T)) \calc(T)\begin{pmatrix}\frac{\partial}{\partial a} H \\ \frac{\partial}{\partial b} H\end{pmatrix}(V_\Theta(Y_x,T))\rp .  \end{align*}
	By a direct computation, the variance term in the last line equals
	\beq\label{alpha-var} 
	\calc_0(\al)\frac{\int_0^T |\si(s,x)|^{2p}\,\dd s}{\lp \mu_{\Phi(p;\cdot)}\int_0^T {|\si(s,x)|}^p\,\dd s\rp^2}.  \eeq
	Since $\al\mapsto \Ga_r=\Ga_r(\al)$ is continuous by \eqref{Ga-formula} and $r\mapsto \rho_p(r)$ is continuous and satisfies $|\rho_p(r)|\lec r^2$ by Lemma~\ref{cov-add}, also $\al\mapsto \calc_0(\al)$ is continuous by the dominated convergence theorem. Hence, $\calc_0(\wh\al^{(p)}_n) \limp \calc_0(\al)$, and the result follows from \eqref{stable} because $V^n_{\Phi((2)p;\cdot)}(Y_x,T)\limp \bbe[|Z|^{(2)p}]\int_0^T |\si(s,x)|^{(2)p}\,\dd s$ by Theorem~\ref{LLN} and $\bbe[|Z|^{2p}] = 2^p\Ga(\frac{2p+1}{2})\pi^{-1/2}$. 
	
	The last statement holds because neither $\calc_0(\wh \al^{(p)}_n)$ nor $\frac{V^n_{\Phi(p;\cdot)}(Y_x,T)}{\sqrt{V^n_{\Phi(2p;\cdot)}(Y_x,t)}}$ depends on $\kappa$ or $\lambda$.
	
	For general $N$, \eqref{conf-int} follows from the univariate case and the fact that functionals based on increments at different locations are asymptotically independent; see \eqref{cov-mpv}.
	\epr
	
	\bpr[Proof of Theorem~\ref{correst}] By Corollary~\ref{multipower-result}, we have $\frac{V^n_{\Psi_m(1,1;\cdot)}(Y_x,T)}{V^n_{\Phi_m(2;\cdot)}(Y_x,T)} \limp \frac{\mu_{\Psi_m(1,1;\cdot)}}{\mu_{\Phi_m(2,\cdot)}} = \Ga_1 = 2^{-\al/2}-1$. So part (1) follows because $F$ is the inverse of $\al\mapsto \Ga_1(\al) = 2^{-\al/2}-1$. The proof of part (2) is completely analogous to the proof of Theorem~\ref{est-alpha} and is therefore omitted.
	\epr
	
	\bpr[Proof of Theorem~\ref{unknown-alpha}] The first part follows immediately from Theorem~\ref{LLN} and the fact that $\wh\tau_n/\tau_n \limp 1$ by \eqref{taun-asymp} and the consistency of $\al_n$ for $\al$. For the second part, we decompose
	\begin{align*} &\frac{\del^{-\frac12}}{|\log\del|} 
	(\wh V^n_{\Phi_m|\Psi_m}(Y_x,T) - V_{\Phi_m|\Psi_m}(Y_x,T))\\ &\quad=\frac{\del^{-\frac12}}{|\log\del|} (\wh V^n_{\Phi_m|\Psi_m}(Y_x,T)- V^n_{\Phi_m|\Psi_m}(Y_x,T))\\
	&\qquad+\frac{\del^{-\frac12}}{|\log\del|} (V^n_{\Phi_m|\Psi_m}(Y_x,T) - V_{\Phi_m|\Psi_m}(Y_x,T)). \end{align*}
	The second term converges to $0$ in probability by Corollary~\ref{multipower-result} (2), so we only need to analyze the first term further, which can be written as
	\beq\label{rewrite} \frac{\del^{-\frac12}}{|\log\del|} \wh V^n_{\Phi_m|\Psi_m}(Y_x,T) \frac{\tau_n^{w_m}-\wh\tau_n^{w_m}}{\tau_n^{w_m}}. \eeq
	If $g(x)=\pi^{ d/ 2-x}\Ga(\frac x 2)/((2\kappa)^{ x /2}(1-\frac x 2)\Ga(\frac d 2))$, then 
	\begin{align*} \wh\tau_n^{w_m} &= g(\al_n)^{\frac {w_m} 2}\del^{(1-\frac{\al_n} 2)\frac {w_m} 2},\\ \tau_n^{w_m} &= g(\al)^{\frac {w_m} 2}\del^{(1-\frac \al 2)\frac {w_m} 2}+o(\del^{(1-\frac\al 2)\frac {w_m} 2})\end{align*} by \eqref{taun-asymp}. So the mean value theorem shows that there is some $\ov \al_n$ between $\al_n$ and $\al$ such that $\wh \tau_n^{w_m} - \tau_n^{w_m}$ equals
	\begin{align*} &\lp\frac {w_m} 2 g(\ov\al_n)^{\frac {w_m} 2-1}\del^{(1-\frac{\ov\al_n} 2)\frac {w_m}2}-\frac {w_m} 4 g(\ov\al_n)^{\frac {w_m} 2}\del^{(1-\frac{\ov\al_n}2)\frac {w_m}2}\log\del\rp (\al_n-\al) \\ &\quad+ o(\del^{(1-\frac\al2)\frac {w_m}2}). \end{align*}
	Inserting this in \eqref{rewrite} and dividing by $\tau_n^{w_m}$, we see that the $o(\del^{(1-\al/2) w_m/2})$-part can be neglected. Thus,
	\beq\label{rewrite-2}\begin{split} \quad&-\del^{-\frac12} (\al_n-\al)\frac{\wh V^n_{\Phi_m|\Psi_m}(Y_x,T)}{|\log\del|}\\
		&\quad\times  \lp\frac {w_m} 2 g(\ov\al_n)^{\frac {w_m} 2-1}\frac{\del^{(1-\frac{\ov\al_n} 2)\frac {w_m}2}}{\tau_n^{w_m}}-\frac {w_m} 4 g(\ov\al_n)^{\frac {w_m} 2}\frac{\del^{(1-\frac{\ov\al_n}2)\frac {w_m}2}}{\tau_n^{w_m}}\log\del\rp. \end{split}\raisetag{-2.75\baselineskip}\eeq
	Since $\al_n$ converges to $\al$ at a rate of $\Del^{1/2}$, one can check that the first term in parentheses converges in probability to $\frac {w_m} 2 g(\al)^{-1}$. And because $\wh V^n_{\Phi_m|\Psi_m}(Y_x,T) \limp  V_{\Phi_m|\Psi_m}(Y_x,T)$ as seen above, the logarithmic factor cancels the first part in \eqref{rewrite-2}. 
	What remains is
	\beq\label{rewrite-3} -\wh V^n_{\Phi_m|\Psi_m}(Y_x,T)\lp \frac {w_m} 4 g(\ov\al_n)^{\frac {w_m} 2}\frac{\del^{(1-\frac{\ov\al_n}2)\frac {w_m}2}}{\tau_n^{w_m}}\rp \del^{-\frac12}(\al_n-\al). \eeq
	With similar arguments as before, the term in parentheses converges in probability to $\frac{w_m}{4}$, and we have $\wh V^n_{\Phi_m|\Psi_m}(Y_x,T)\limp V_{\Phi_m|\Psi_m}(Y_x,T)$. Moreover, by Theorem~\ref{est-alpha} (2) and Theorem~\ref{correst} (2), we have $\del^{-1/2}(\al_n-\al) \lims N(0,\calv)$ with
	\[ \calv=\begin{cases} \frac{1}{N^2} \displaystyle\sum_{j=1}^N\textstyle\calc_0(\al)\frac{\int_0^T |\si(s,x_j)|^{2p_0}\,\dd s}{\lp V_{\Phi_j(p_0,\cdot)}(Y_x,T)\rp^2} &\text{if } \al_n=\wh\al^{(p_0)}_n,\\ \frac{1}{N^2} \displaystyle\sum_{j=1}^N\textstyle\wt\calc_0(\al)\frac{\int_0^T |\si(s,x_j)|^4\,\dd s}{\lp V_{\Psi_j(1,1;\cdot)+\Phi_j(2,\cdot)}(Y_x,T)\rp^2} &\text{if } \al_n=\wt\al_n; \end{cases} \]
	see \eqref{alpha-var} in particular. As this limit is the same for all $m$, the assertion follows from \eqref{stable} and studentization.
	\epr
	
	\bibliographystyleB{plainnat}
	\bibliographyB{bib-VolaEstimation}

\end{document}